\documentclass[12pt,a4paper]{article}\usepackage[top=20mm,left=12mm]{geometry}\textwidth18.2cm\textheight25.2cm
\usepackage{amsmath,amssymb,amsthm}
\usepackage{float}\floatstyle{plaintop}\restylefloat{table}
\usepackage[tableposition=top]{caption}
\usepackage{subcaption}\usepackage{tikz}
\usepackage{graphicx,pstricks,fancyvrb,multirow,lineno,url}
\usepackage{aliascnt} 
\usepackage{enumitem}
\usepackage{listings}\usepackage{paralist, longtable}

\theoremstyle{theorem}
\newtheorem{theorem}{Theorem}[section]
\newaliascnt{definition}{theorem}
\newtheorem{definition}[definition]{Definition}
\newaliascnt{lemma}{theorem}
\newtheorem{lemma}[lemma]{Lemma}
\newtheorem{corollary}[theorem]{Corollary}
\theoremstyle{definition}
\newaliascnt{remark}{theorem}
\newtheorem{remark}[remark]{Remark}

\aliascntresetthe{lemma}
\aliascntresetthe{remark}
\aliascntresetthe{definition}
\def\bexe{\begin{exercise}}\def\eexe{\eex\end{exercise}}
\def\bsol{\begin{solution}}\def\esol{\eex\end{solution}}
\def\bexa{\begin{example}}\def\eexa{\end{example}}
\def\brem{\begin{remark}}\def\erem{\end{remark}}
\def\bthm{\begin{theorem}}\def\ethm{\end{theorem}}
\def\blem{\begin{lemma}}\def\elem{\end{lemma}}
\def\bcor{\begin{corollary}}\def\ecor{\end{corollary}}
\def\bdefi{\begin{definition}}\def\edefi{\end{definition}}

\def\bmip{\begin{minipage}{\textwidth}}\def\emip{\end{minipage}}
\def\huga#1{\begin{gather} #1 \end{gather}}
\def\hugast#1{\begin{gather*} #1 \end{gather*}}
\def\hual#1{\begin{align} #1 \end{align}}
\def\hualst#1{\begin{align*} #1 \end{align*}}
\def\hueq#1{\begin{equation} #1 \end{equation}}
\def\hueqst#1{\begin{equation*} #1 \end{equation*}}
\newcommand{\btab}[2]{\begin{tabular}{#1}#2\end{tabular}}

\newcommand{\R}{{\mathbb R}}
\newcommand{\C}{{\mathbb C}}
\newcommand{\Z}{{\mathbb Z}}

\def\CC{{\cal C}}\def\CT{{\cal T}}
\def\CV{{\cal V}}\def\CA{{\cal A}}\def\CD{{\cal D}}  
\def\CG{{\cal G}}\def\CH{{\cal H}}\def\CL{{\cal L}}
\def\CO{{\cal O}}\def\CP{{\cal P}}
\def\CM{{\cal M}}\def\CB{{\cal B}}
\def\CX{{\cal X}}

\def\per{{\rm per}}

\def\ga{\gamma}\def\om{\omega}
\def\noi{\noindent}\def\ds{\displaystyle}
\def\vt{\vartheta}\def\pa{{\partial}}\def\lam{\lambda}
\newcommand{\bi}{\begin{itemize}}\newcommand{\ei}{\end{itemize}}
\newcommand{\ben}{\begin{enumerate}}\newcommand{\een}{\end{enumerate}}
\newcommand{\bce}{\begin{center}}\newcommand{\ece}{\end{center}}
\newcommand{\bci}{\begin{compactitem}}\newcommand{\eci}{\end{compactitem}}
\newcommand{\bcen}{\begin{compactenum}}\newcommand{\ecen}{\end{compactenum}}
\newcommand{\bcena}{\begin{compactenum}[(a)]}
\newcommand{\reff}[1]{(\ref{#1})}

\newcommand{\spr}[1]{\left\langle #1 \right\rangle}
\newcommand{\hs}[1]{{\hspace{#1}}}\newcommand{\vs}[1]{{\vspace{#1}}}

\def\eps{\varepsilon}
\def\ra{\rightarrow}

\newcommand{\barr}{\begin{array}}\newcommand{\earr}{\end{array}}
\newcommand{\bpm}{\begin{pmatrix}}\newcommand{\epm}{\end{pmatrix}}
\newcommand{\bsm}{\left(\begin{smallmatrix}}
\newcommand{\esm}{\end{smallmatrix}\right)}
\newcommand{\ba}{\begin{array}}\newcommand{\ea}{\end{array}}
\def\dd{\, {\rm d}}\def\ri{{\rm i}}

\def\cc{{\rm c.c.}}\def\er{{\rm e}}

\def\om{\omega}\def\Om{\Omega}
\def\ddt{\frac{\rm d}{{\rm d}t}}
\def\bphi{\varphi}\def\del{\delta}\def\brho{\varrho}

\def\eex{\hfill\mbox{$\rfloor$}}

\def\per{{\rm per}}
\def\sign{{\rm sign}}
\def\Del{\Delta}

\def\sig{\sigma}
\def\al{\alpha}

\def\Lam{\Lambda}
\def\kap{\kappa}\def\Sig{\Sigma}

\def\eps{\varepsilon} 
  
\def\pdep{{\tt pde2path}}

\newlength{\tew}

\def\ig{\includegraphics}

\def\bsub{\begin{subequations}}
\def\esub{\end{subequations}}
\def\wit{\widetilde}

\def\ddT{\frac{{\rm d}}{{\rm d}T}}
\renewcommand{\arraystretch}{1.05}\renewcommand{\baselinestretch}{1.0}
\def\medskip{}\def\bigskip{}
\def\sm{\small }\def\rb{\raisebox }
\usepackage{setspace}
 \usepackage{xcolor}\usetikzlibrary{calc}

\def\huc#1{{\magenta #1}} 
\def\huc#1{{#1}} 

\def\Cr{\mathcal C_r} \def\Co{\mathcal C} \def\bp{\bphi} \def\sG{\mathcal G}
\def\per{\text{per}}    
\def\Psinum{\Psi_{\text{num}}}\def\tol{{\tt tol}}
\def\hra{\hookrightarrow}

\usepackage{hyperref}\usepackage{cleveref}
\def\taskip{\renewcommand{\arraystretch}{1}\renewcommand{\baselinestretch}{1}}
\def\teskip{\renewcommand{\arraystretch}{1.1}\renewcommand{\baselinestretch}{1.1}}

\def\hutab#1{\taskip\begin{\table}#1\end{table}\teskip}
\newfloat{Algorithm}{ht}{lop}[section]
\def\slfig{\opt{hu,ho}{\begin{figure}[ht]}\opt{sl}{\begin{figure}[H]}}
\def\slfigH{\opt{hu,ho}{\begin{figure}[H]}\opt{sl}{\begin{figure}[H]}}
\def\sltab{\opt{hu,ho}{\begin{table}[ht]}\opt{sl}{\begin{table}[H]}}

\def\vol{{\rm vol}}\def\gptool{{\tt gptoolbox}}\def\T{\mathbb{T}} 

\def\D{\mathcal{D}} 
\def\cc{\mathrm{c.c.}}
\def\muti{\tilde{mu}}

\allowdisplaybreaks

\begin{document}
\text{}\vspace{0mm} 
\bce\Large
Helfrich cylinders -- instabilities, bifurcations and amplitude equations
\\[4mm]
\normalsize 
Alexander Meiners, Hannes Uecker\\[2mm]
\footnotesize
Institut f\"ur Mathematik, Universit\"at Oldenburg, D26111 Oldenburg, 
alexander.meiners@uni-oldenburg.de, hannes.uecker@uni-oldenburg.de\\[3mm]
\normalsize

\today
\ece 

\begin{abstract} \noindent 
Combining local bifurcation analysis with 
numerical continuation and bifurcation methods 
we study bifurcations from cylindrical vesicles 
described by the Helfrich equation with volume and area constraints, 
with a prescribed periodicity along the cylindrical axis. 
The bifurcating solutions are in two main classes, 
axisymmetric (pearling), and non--axisymmetric (coiling, buckling, and 
wrinkling),  
and depending on the spontaneous curvature and the prescribed 
periodicity along the cylinder axis we obtain different stabilities 
of the bifurcating branches, and different secondary bifurcations. 
\end{abstract}

\noi

\noi{\bf MSC} 
35B32, 58J55, 65P30

\tableofcontents


\section{Introduction}
The (dimensionless) Canham--Helfrich model \cite{Can70,H73} 
(or spontaneous curvature model) for vesicles 
with lipid bilayer membranes is based on the energy 
\hueq{\label{helfen} E(X)= \int_X (H-c_0)^2\dd S, 
}
where $X$ is a (possibly time $t$ dependent) 
2D manifold embedded in $\R^3$ with mean curvature 
$H$, and where the parameter $c_0\in\R$ models an energetically 
favorable spontaneous (or preferred) curvature; it depends on the properties of 
the lipid bilayer membrane and can be zero, or positive, or negative. 
 Roughly speaking, 
the aim is to minimize $E$ 
for prescribed area $\CA(X)$ 
and enclosed volume $\CV(X)$, 
i.e., under the constraints 
\huga{\label{c1}
\text{
$q_1(X):=\CA(X)-\CA_0=0$ and $q_2(X):=\CV(X)-\CV_0=0$}. 
} 
There are numerous works on this, for instance 
using the minimization of $E$ to explain the shape of red blood 
cells, see, e.g., \cite{SZ89,SBL90,peter83,SL95,sei97}.%
\footnote{For $c_0=0$, $E$ is also called Willmore energy. 
More generally, a term $b\int_X K\dd S$ can be added to 
$E$, where the parameter $b\ge 0$ is called saddle--splay modulus, 
but due to the Gauss--Bonnet theorem 
$\int K dS$ is a topological invariant for closed $X$, and 
hence can be dropped from the minimization of $E$, and the same 
holds for the cylindrical topology considered here.} 
The minimization 
proceeds by introducing the Lagrangian 
\hueq{\label{Lagra}F(X,\Lam) =E(X)+\lam_1(\CA(X)-\CA_0)
+\lam_2(\CV(X)-\CV_0), }
with Lagrange multipliers $\Lam=(\lam_1,\lam_2)\in\R^2$, where 
$\lam_2$ is the pressure difference between inside and outside; since 
we choose the inner normal, here negative $\lam_2$ corresponds to 
lower pressure inside; $\lam_1$ corresponds to a surface tension, 
and can also have either sign; see, e.g., \cite{BH17} for discussion. 

\begin{figure}[ht]
\bce
\btab{lll}{{\small (a)}&{\sm (b)}&\hs{4mm} {\small (c)}  \\
\hs{-2mm}\ig[height=0.18\tew]{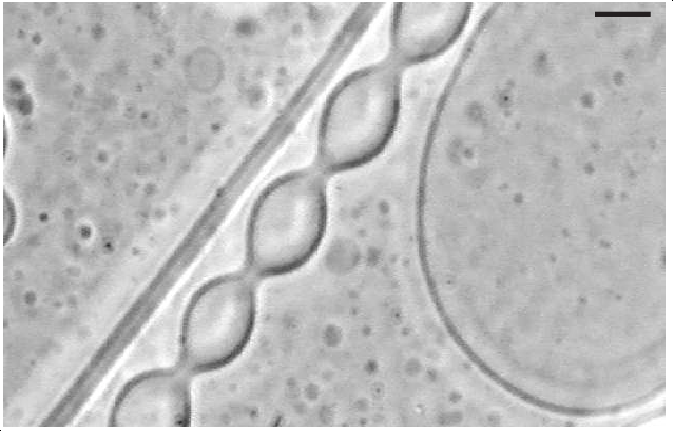}&\ig[height=0.18\tew]{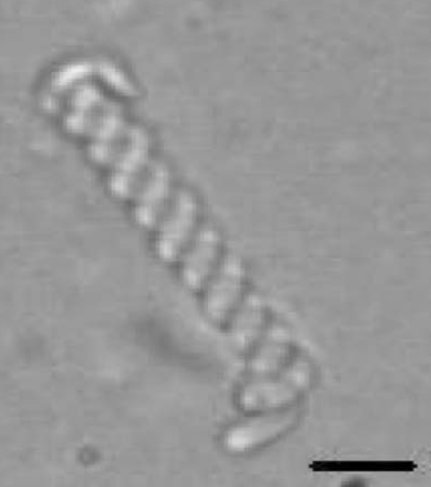}& 
\hs{1mm}\rb{-5mm}{\ig[height=0.24\tew]{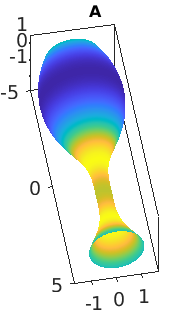}\hs{-2mm}\ig[height=0.24\tew]{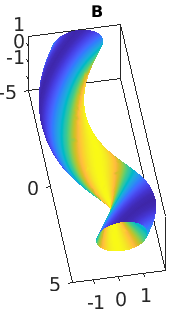}\hs{-2mm}
\ig[height=0.24\tew]{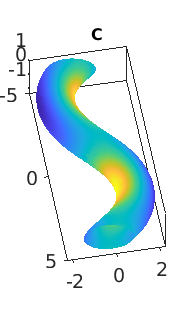}\hs{-1mm}\ig[height=0.24\tew]{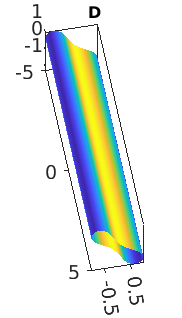}}
}
\ece 
\vs{-0mm}
\caption{\small{(a,b) Cylindrical vesicle shapes in experiments, 
pearling (a) from \cite{CR98}, and coiling (b) from \cite{HEJ14}; scale bars 
representing 10$\mu$m. 
(c) Four types of steady states bifurcating from straight cylinders 
in \reff{hflow}, with fixed length $L=10$ (and $c_0=0$): 
A pearling, B coiling, C buckling, D wrinkling; 
see \S\ref{numex} for detailed bifurcation diagrams and stability information.  
}  
\label{f0}}
\end{figure}

From \reff{Lagra} we obtain the Euler--Lagrange equation 
\hueq{\label{helf0} 
G(X,\Lam):=\Del H +2H(H^2-K) +2c_0K-2(c_0^2+\lam_1)H-\lam_2=0 
}
under variations in normal direction together with the constraints \reff{c1}, 
where  $K$ is the Gaussian curvature, 
and $\Delta$ the Laplace--Beltrami operator on $X$. 
Many works also use \reff{helfen} to study instabilities of 
tubular vesicles, e.g., \cite{OuHe89,BZM94,FTS99,TTS01}, 
in particular the so--called pearling and coiling instabilities, 
see Fig.\ref{f0}(a,b), 
and Remark \ref{irem1} for a review of the (physics) literature.

In, e.g., \cite{NY12,rupp24}, the 
constrained $L^2$--gradient flow 
associated to \reff{helfen} is studied, which can be written as 
\hueq{\label{hflow}\bpm V_t(t)\\ 0\epm=\bpm G(X(t),\Lam(t))\\ q(X(t))\epm, 
}
where $q=(q_1,q_2)$ and $V_t(t)$ is the normal velocity of $X$, and 
with now time dependent Lagrange multipliers $\Lam(t)$ which are determined 
from \reff{c1}, such that 
\huga{\label{helf1b}
\ddt E(X(t),\Lam(t))=-\|V(t)\|_{L^2(X(t))}^2,\quad \ddt \CA(X(t))=0,\quad 
\text{and }\ddt \CV(X(t))=0, 
}
as long as the flow exists. We also refer to \cite{Ga13} and 
\cite{BDGP23} for introductions to, reviews of, and extended 
bibliographies concerning geometric flows, 
both analytically and numerically, and also including mean--field 
approaches.

A cylinder $\CC(r)$ of radius $r$ is a solution to \reff{helf0} iff 
\hueq{\label{socon} 
\lam_1=\frac 1 2 \left( \frac 1 {2r^2} -2r\lam_2-2c_0^2\right).
}
However,  $E,\CA$ and $\CV$ are not well defined 
for an infinite cylinder,
and in the following we restrict to cylinders  $\CC_L(r)$ of lengths $L$. 
For these we have the scaling 
$E(\CC_L(\al r); c_0)= E(\CC_{\al L}(r); c_0/\al)$ and thus 
it is sufficient to consider $r=1$ with a variable length $L$, and we abbreviate 
$\CC_L=\CC_L(1)$. 
After showing well--posedness of \reff{hflow} for initial 
conditions $X_0$ close to $\CC_L$, 
where we assume periodic boundary conditions along the 
cylindrical axis, we study bifurcation 
diagrams of steady states of \reff{hflow}. 
In Fig.\ref{f0}(c) we preview four types of solutions from branches 
bifurcating from cylinders, obtained from numerical continuation and 
bifurcation, namely again pearling (A) and coiling (B) like in (a) and (b), 
and buckling (C), and wrinkling (D), which break the symmetry of the cylinder 
in different ways. 
\footnote{While pearling is a standard name for A, the 
nomenclature for B and in particular C varies. For instance, 
B is sometimes called writhing or twisting, and  
in \cite{NSS15}, C is called wrinkling.  We could not 
find explicit names for states like D in the ``vesicles literature'', 
and our choice of ``wrinkling'' is inspired by \cite{FVGK22}, which 
treats a related model for arteries with a restoring force.}
\footnote{The symmetry group of the problem at the cylinder is 
$O(2)\times O(2)$, translations and reflections in $x$, 
and in $\bphi$, where $x$ is along the cylindrical axis and 
$\bphi$ the polar angle 
orthogonal to it. Additionally there are the trivial symmetries 
of translation in $y$ and $z$, the euclidian coordinates orthogonal 
to $x$, but these are mostly relevant from a numerical point 
of view, \huc{i.e., must be dealt with via phase conditions in the numerical 
continuation.}}
The starting point for the bifurcation analysis is the dispersion relation 
$\mu=\mu(m,n)$ 
for the linearization of \reff{hflow} at $\CC_L$, which can be computed 
explicitly from an ansatz 
$$\text{$X=\CC_L+\er^{\mu t+\ri(k_m x+n\bphi)}N$,\,\, 
$m,n\in\Z$,\, \, $k_m=2\pi m/L$.}
$$
We study these dispersion relations with $c_0$ and 
$\lam_2$ (and $\lam_1$ given by \reff{socon}) as primary parameters, 
and for suitable $c_0$ and $L$ we find (typically narrow) 
ranges in $\lam_2$ 
of stability of the ``trivial branch'' $X=\CC_L$, with 
branch points at the ends. At one end, axisymmetric 
branches of pearling type A bifurcate, while at the other end 
the loss of stability is to B,C or D, depending on the parameters. 

After identifying these instabilities, based on the well--posedness 
results for the flow we use a centre manifold reduction to derive amplitude 
equations (AEs) for the critical modes, and from these also obtain 
approximations of the bifurcating steady state branches. The derivation 
of these AEs exploits the $O(2)\times O(2)$ symmetry of the problem, 
but is rather non--standard due to the constraints \reff{c1}, 
and in some sense 
our methods could be described as bifurcation without parameters \cite{Lie15}, 
see Remark \ref{irem0}. 

Finally, we use the package \pdep\ \cite{p2pbook,geompap} to 
numerically continue the primary bifurcating branches 
away from their bifurcation points from $\CC_L$. Close to bifurcation 
these numerical results agree well with the AE description, 
and additionally we obtain the behavior and stability at larger amplitude, 
and secondary bifurcations. In summary, this yields rather complicated 
bifurcation diagrams (BDs), 
which we discuss in detail for three exemplary 
cases, namely $c_0=0, c_0=-1$, and $c_0=0.48$, with strong qualitative 
differences, and we complement these steady state BDs with 
numerical flows for \reff{hflow} to give further indication of the 
typical behavior of \reff{hflow}. 

\brem\label{irem0} We start the linear stability analysis 
and also the numerical bifurcation analysis by continuation 
of the trivial solution (branch) $\CC_L$ in the Lagrange multiplier $\lam_2$ 
(with $\lam_1$ given by \reff{socon}). We then also often continue 
the nontrivial branches (e.g., of type A--D from Fig.\ref{f0}) in $\lam_2$ 
as a primary parameter, with fixed $\CA$ but without fixing $\CV$. 
However, the stability is always computed wrt to 
\reff{hflow}, i.e., for the area and volume conserving flow. 
This point of view of allowing the volume $\CV$ to change along the branch 
while considering stability under $\CV$ conservation is motivated by 
the closed vesicle situation discussed for example in \cite{JSL93} and 
\cite{geompap}. Roughly speaking, for closed vesicles, 
the sphere is the unique maximum of volume at fixed area. 
Hence volume and area preserving branches bifurcating 
from the sphere do not exist. The same holds for the 
cylinder if we also fix the length $L$. Cylinders 
can be destabilized preserving $\CV$ and $\CA$ when allowing $L$ to change, 
see, e.g., \cite{OuHe89}, but we found our method with fixed $L$ more 
convenient. 
\eex \erem
\brem\label{irem1} 
a) Following \cite{BZM94}, \cite{GNPS96} considers the dynamics of pearling 
based on the Helfrich energy and 
treating the hydrodynamics in the tube via a lubrication approximation. 
This leads to a 6th order Cahn--Hilliard 
type equation for axially symmetric perturbations of the cylinder, 
with $\CV$ conservation, and the critical wave numbers and growth rates 
of the pearling (including the propagation speed of pearling fronts) are 
carefully compared to \cite{BZM94}. 
The coiling instability and the dynamics of coiling are 
discussed \cite{FTS99,Spi01}, 
and compared to the dynamics of pearling in \cite{TTS01}. For both, 
it is argued that polymers anchoring at lipid bilayers {\em locally} 
generate nonzero values for the 
spontaneous curvature $c_0$, {\em and} generate significant differences 
in effective area between inside and outside of the bilayer(s), which 
points at area difference 
elasticity (ADE) terms. In particular, in \cite{TTS01}, 
measurements of diffusion/propagation speeds suggest that 
for local injection of polymers in the outer fluid, 
initially ADE terms (with bulk diffusion of polymers) give the right model, 
and spontaneous curvature terms (with surface diffusion of polymers) 
when reaching equilibrium, both with $\CA$ and $\CV$ conservation. 
On the other hand, \cite{CR98} argues that in certain 
experiments $c_0\ne 0$ is a global quantity. 
Finally, an important feature of coiling is that it rather occurs for 
multilamellar tubes, and that in equilibrium the coiling is usually ``tight'', 
as in Fig.\ref{f0}(b), and ``stretched coils'' as B in (c) only occur 
transiently.  See also \cite{St02} for an overview of these results.

Further experiments on coiling are discussed in \cite{PAM09},  
with various hypothesis for the underlying mechanism, and 
\cite{HEJ14} discusses the role of the actin cortex (not considered in \reff{hflow}) in suppressing pearling 
in actual cells.  In \cite{NSS15}, 
tubular vesicles modeled by \reff{helfen} 
with $c_0=0$ and $\CA$ and $\CV$ conservation, are subjected to external 
elongational and compressional flows, and numerical simulations are 
compared to experimental results. Depending on flow parameters 
(condensed to the capillary number), and the initial vesicle shape 
(including, e.g., the shape of the spherical caps at the ends of finite 
cylinders), pearling is observed for elongational flows, 
and buckling and coiling (called wrinkling in \cite{NSS15}) for 
compressional flows. This is also complemented by linear stability 
analysis for infinite cylinders. 

b) In our model \reff{hflow} we do not consider any interaction with 
surrounding fluids; we treat the idealized infinitely long 
cylinders, or the idealization of a fixed spatial periodicity; 
and $c_0$ and the surface tension $\lam_1$ and pressure difference $\lam_2$ 
are always global. 
Hence, mismatches to experiments (small range 
of stability for periodic pearling, long wave nature of the 
buckling and coiling, see below) are not completely surprising, but we 
believe/hope that 
our analytical (instabilities of $\CC_L$, derivation of the AEs) 
and numerical (continuation of the bifurcating branches to larger amplitude) 
discussion of \reff{hflow} will also help the further 
modeling. 
\eex \erem

\paragraph{Outline, and overview of results} 
In \S\ref{prewell} we recall the needed differential geometry and rewrite 
\reff{hflow} as an equation in the normal displacement $u$ from the cylinder. 
This gives a constrained quasilinear parabolic problem in $u$, for which 
we use maximal regularity  in Sobolev spaces 
to obtain local existence close to $\CC_L$. 
From this we also obtain 
that nonlinear stability follows from linearized stability, and since 
$\CC_L$ has constant mean and Gauss curvatures, the eigenvalues and eigenfunctions of the linearized problem can be computed explicitly 
allowing an analytical treatment of the instabilities of $\CC_L$, 
and of the bifurcating branches via amplitude equations in \S\ref{stcy}  
and \S\ref{aesec}. This is straightforward, but cumbersome 
due to the three parameters $c_0, L$ and $\lam_2$, and 
requires a careful look at the area and volume constraints \reff{c1} already 
in the linearization. 
The results of the linear stability and AE analysis can be 
summarized as follows: 
\bcen 
\item The coiling/buckling instability 
of $\CC_L$ is always of long wave type, i.e., the primary coiling/buckling 
instability has the maximal period given by the length $L$. This contradicts 
many experimental results for coiling, 
where the instability seems of finite wavelength type, 
and after onset the cylinders coils up tightly as in 
Fig.\ref{f0}(b). We could not find clear experimental results on buckling of 
vesicles. 
\item For some (positive) $c_0$ the primary pearling instability can be of 
finite wave number type, i.e., independent of $L$, corresponding 
to Fig.\ref{f0}(a). However, see 5.~below. 
\item The wrinkling instability 
(homogeneous along the cylindrical axis) naturally is independent of $L$. 
\item 
Due to the constraints, 
also the derivation and interpretation of the AEs 
is non--standard. In particular, the pearling and wrinkling 
branches bifurcating at stability 
loss of $\CC_L$ are always stable. For the coiling vs buckling at stability 
loss of $\CC_L$, always 
exactly one of these bifurcating branches is stable, and the instability 
of the unstable one is wrt perturbations into direction of the stable one. 
\ecen 

In \S\ref{numex} we complement our analysis with 
numerical continuation and bifurcation methods, and also with some 
numerical flows for \reff{hflow}, aka direct numerical simulation (DNS). 
As expected and already said, near bifurcation the numerics 
agree well with the analysis, but one remarkable result is the following: 
\bcen 
\setcounter{enumi}{4}
\item 
In case of finite wave number pearling the stability range 
of the bifurcating pearling branch rather strongly depends on $L$: 
The pearling branch becomes unstable via exchange of stability 
with a bifurcating branch of a ``localized single pearl'', and the 
associated bifurcation point moves closer to the primary bifurcation 
with increasing $L$. This suggests that finite wavelength $O(1)$ 
amplitude pearling as in Fig.\ref{f0}(a) does not arise via bifurcations 
at stability loss of $\CC_L$ in \reff{hflow}.  
\ecen 
See \S\ref{numex} for details of 
this, and also for secondary bifurcations from other primary branches. 
\huc{On the other hand: 
\bcen\setcounter{enumi}{5}
\item For the case of ``large'' $c_0$ ($c_0>1/2$, where $\CC_L$ is not stable 
for any $(\lam_1,\lam_2)$) we do find finite wavelength $O(1)$ pearling, 
but not via bifurcation from $\CC_L$. Instead we use DNS to converge 
to such solutions, which we subsequently can continue again in parameters. 
However, we only present this as an outlook in \S\ref{cpsec}, and leave 
a systematic study of this setting for future work. 
\ecen 
}
Another interesting result from the numerical flows is that cylindrical 
{\em shapes} can be more stable than linearized stability analysis 
of steady 
{\em solutions} $(X,\lam_1,\lam_2)$ at fixed $\CA_0$ and $\CV_0$ suggests. 
Namely: 
\bcen 
\setcounter{enumi}{6}
\item 
Non area and/or volume conserving perturbations, of, e.g., unstable 
coiling, may yield initial conditions which converge to nearby stable  
coiling at the new conserved $(\CA_0,\CV_0)$ by mainly 
adjusting the surface tension $\lam_1(t)$ and the pressure $\lam_2(t)$, with 
rather negligible changes of $X(t)$. 
\ecen 
In \S\ref{dsec} we summarize our analysis again, and 
in the Appendix we collect some details of the derivation of the AEs, 
and some remarks on the numerical algorithms. 

\section{Preliminaries and well-posedness of the flow} \label{prewell}
We introduce the basic differential geometric preliminaries in order
to write \reff{hflow} as a constrained parabolic 
quasi-linear equation in the normal displacement $u$ of the
cylinder. This construction is as in \cite{LeSi16}, 
and hence we refrain from extensive computational details and focus 
on the structure and mapping properties of the resulting equations. 
From these we obtain maximal regularity results and hence local existence 
for \reff{hflow}. For general background on differential geometry 
we recommend \cite{Tapp16,UY17}. 
\subsection{Preliminaries} \label{pre}
Let $\CC(r)$ be the two dimensional unbounded cylinder embedded in $\R^3$ of radius $r>0$ with the unbounded direction along the $x$ axis. $\CC(r)$ can be considered as the  embedding 
\huga{\barr{rl} \Phi :\R\times \T &\ra \R^3\\
(x,\bphi)&\mapsto \bpm x , r\cos(\bphi),r\sin(\bphi)\epm^T 
\earr \label{cylc}
} 
of the product manifold $\R\times \T$, where $\T$ is the one dimensional topological torus with length $2\pi$ endowed with the metric $d_\T(\bphi,\vt)
=\min \{ | \bphi-\vt|, 2\pi-|\bphi-\vt|\} $ for all $\bphi,\vt\in \T$. Then as the pullback metric of the $\R^3$ to $\CC(r)$ we have 
$g_0=\dd x^2+ r^2\dd \bphi^2$.
The normal vector (field) $N$ of $\CC(r)$  is oriented inwards, wrt 
$\Phi$, hence $H(\CC_r)=\frac 1 {2r}$. 

For a function $u:\CC(r)\ra \R$, for the remainder of this subsection assumed as smooth as necessary,  we consider surfaces 
\hueqst{ X = \CC(r)-u\, N,} 
assuming $u<r$, such that $X$ is still an embedding of $\R\times \T$; 
explicitly, 
\hueqst{\tilde \Phi(x,\bphi) = \Phi(x,\bphi) + \tilde u(x,\bphi) ( 0, \cos(\bphi), \sin(\bphi))^T,}
where $\tilde u= u\circ \Phi$. 
This gives the metric 
\hueqst{ (g_{ij}) =\bpm 1+u_x^2 & u_xu_\bphi \\ u_xu_\bphi &(r+u)^2+u^2_\bphi \epm \quad \text{with inverse }(g^{ij}) = \frac 1 g \bpm(r+u)^2+u^2_\bphi  & -u_xu_\bphi \\ -u_xu_\bphi & 1+u_x^2\epm  
}
with $g=\det (g_{ij}) = (r+u)^2(1+u_x^2)+u_\bphi^2$. 
In the same explicit fashion one can calculate the second fundamental form  
\hueqst{ (h_{ij}) = \frac 1 {\sqrt{g}} \bpm (r+u)u_{xx} & (r+u) u_{x\bphi}-u_xu_\bphi\\ (r+u) u_{x\bphi}-u_xu_\bphi & (r+u)(u_{\bphi\bphi} -(r+u))-2u_\bphi^2\epm, }
and this gives explicit formulas for the mean and Gaussian curvatures, 
\hual{\label{meanu}H(u)&=\frac{1}{2g^{3/2}}\left(u_\bphi^2-(r+u)\big[ ((r+u)^2+u_\bphi^2)u_{xx} +(1+u_x^2)u_{\bphi\bphi} -2u_xu_\bphi u_{x\bphi}\big]\right)+\frac 1 {2g^{1/2}},\\
\label{gaussu} K(u)&= \frac{1}{g^2}\left((r+u)u_{xx} \Big[ (r+u)\big(u_{\bphi\bphi}-(r+u)\big) -2u_\bphi^2\Big]-\big((r+u)u_{x\bphi}-u_xu_\bphi\big)^2\right).
}
 
We denote by $\Delta_g$ the negative definite Laplace--Beltrami operator 
\hueq{\label{LBu}
\Delta_u=\frac{1}{g^{1/2}}\left( \pa_x\bigg[ \frac{(r+u)^2+u_\bphi^2}{g^{1/2}}\pa_x -\frac{u_xu_\bphi}{g^{1/2}}\pa_\bphi\bigg] + \pa_\bphi\bigg[ -\frac{u_xu_\bphi}{g^{1/2}} \pa_x+\frac{1+u_x^2}{g^{1/2}} \pa_\bphi \bigg]\right),}
which is self adjoint in $L^2(X)$. 
Furthermore, the normal velocity is 
\hueq{V_t=-\spr{\ddt X,N}= \frac{(r+u)u_t}{g^{1/2}}=\eta(u)u_t
}
where $\eta(u)>0$ for $u<r$, as we assumed. 
By all of this we can express \reff{hflow} as 
a PDE in $u$, 
and with a slight 
abuse of notation we rewrite \reff{hflow} as 
\bsub\label{hflowu}
\hual{\eta(u) \pa_t u&=G(u,\Lam),\\ 
0&=q(u), 
} 
\esub 
where $G$ has the quasi-linear structure
$G(u,\Lam){=}B(u)u+C(u)$, 
where $B(u){=}\sum_{|\al|=4}b_\al \pa_x^{\al_1}\pa_\bp^{\al_2}$, 
where $b_\al=b_\al(u, \nabla u,\nabla^2u)$, see \reff{psB}, 
and where similarly $C$ is a nonlinear operator of order up to three, 
but linear in the third derivatives.

Naturally, we shall need linearizations of \reff{hflowu} at 
$\CC_L(r)$ and at $X(u)=\CC_L(r)-uN$. For this we first note the 
formulas for derivatives of area $\CA$ and volume $\CV$, where we 
  distinguish between the derivative of a functional $\CH(u)$ acting 
on $v$ in the $L^2$ scalar product, denoted by $\spr{\pa_u \CH,v}$, 
and the integral kernel $\pa_u\CH$ as a differential operator on $X$. 

\blem\label{aderlem} For $u,v\in C^\infty_\per(\CC_L(r))$ sufficiently small, 
the area of $X(u)=\CC_L(r)-u\, N$ expands as  
$
\CA(u+v)= \CA(u)+\spr{\pa_u \CA(u),v} +\frac 1 2 \spr{\pa_u^2\CA(u)v,v}
+ O(\|v\|^3)$
with 
\hueq{\label{derarea} 
\spr{\pa_u \CA(u),v}{=}2 \int_{X} H(u)\eta(u) v\dd S, 
\ \ \spr{\pa_u^2 \CA(u)v,v}{=}\int_{X} -(\Delta_u v+2(H^2(u)-K(u))v)v\dd S.
}
The volume 
\hueq{\CV(u)=\frac 1 2\int_{X} \frac{(r+u)^2}{g^{1/2}(u)} \dd S}
expands as
$\CV(u+v)= \CV(u)+\spr{\pa_u \CV(u),v}+\frac 1 2 \spr{\pa_u^2\CV(u)v,v}
+O(\|v\|^3)$ 
with 
\hueq{\label{dervol}
\spr{\pa_u \CV(u),v}= \int_{X} \eta(u) v \dd S, \qquad 
\spr{\pa_u^2 \CV(u)v,v}= \int_{X} g^{-1/2}(u) v^2\dd S.}
\elem

\brem\label{prrem}
Since $\Lam$ appears linearly in $G(u,\Lam)$, the 
first derivative $\ddt q(u)=0$ gives a linear system of equations 
\hueq{\label{eqlam}
\D\Lam=\bpm \spr{\pa_u E,\pa_u\CA(u) }\\ \spr{\pa_u E,\pa_u\CV(u)}\epm,
}
for $\Lam$, where 
\hueq{\label{ddef} 
\D:=\spr{\pa_u q(u),\pa_\Lam G(u,\lam)}=\bpm 4\int_X H^2(u)\eta^2(u) \dd S & 2\int_X H\eta^2(u) \dd S\\  2\int_X H \eta^2(u) \dd S &\int_X \eta^2(u)\dd S\epm. 
}
 If $\det \D= 4 \left(\int_X\eta^2\dd S \int_X H^2\eta^2 \dd S-\left(\int_X H\eta^2 \dd S\right)^2\right)  \neq 0$, 
then the implicit function theorem gives that in a neighborhood of 
$u$ there is a unique solution $\Lam=\Lam(u)$. 
Following \cite{NY12}, inserting $\Lam(u)$ 
in (\ref{hflowu}a) gives a parabolic integro differential equation. 
However, if $u=0$, then $\det \CD =0$, but 
the right hand side of \reff{eqlam} is orthogonal to $\ker \CD$  
(more generally, this holds at any constant mean curvature surface). 
Then a solution $\Lam$ to \reff{eqlam} exists, and  
$\spr{\Lam_\CC,\ker \CD^T}=0$ gives that 
$2\lam_1=\lam_2$ and $\lam_1=\frac 1 {3\CA_0} \spr{\pa_u E,1}$, i.e., 
$\ker \CD$ corresponds to the family of cylinders in \reff{socon}; 
see \cite{NY12} for further discussion (for the case of closed vesicles 
of spherical topology). 

Here we shall proceed similarly for the cylindrical vesicles; 
we start with well--posedness of (\ref{hflowu}a) for fixed $\Lam\in\R^2$, 
i.e., with {\em inactive} constraints (\ref{hflowu}b). Then 
we first switch on one constraint, wlog $q_1(u)=0$ and free $\lam_1=\lam_1(t)$, 
and subsequently consider the full constraints $q(u)=(q_1(u),q_2(u))=0$ with 
free $\Lam(t)=(\lam_1(t),\lam_2(t))$. 
\eex \erem

\subsection{Well-posedness} \label{well}
Given an initial condition $u_0\in \CX\subset L^2(\CC_L(r))$, 
 i.e., $u_0=u_0(x,\bphi)$ with periodic boundary conditions in $x$, 
where $\CX=H^m_\per(\CC_L(r))$ with the choice $m\in(3,4)$ explained 
below, we 
search for solutions $(t,x,\bphi))\mapsto (u(t,x,\bphi),\Lam(t))$ of \reff{hflowu} on $[0,T)\times \CC_r(L)$. Here 
$\CC_L(r)$ is topologically equivalent to a torus, 
and the periodic Sobolev spaces of order $m\ge 0$ are 
\huga{\label{hmdef}
H^m_\per=H^m_\per(\CC_L(r))=\{v\in L^2 (\CC_L(r)): 
\sum_{k\in\frac {2\pi}{L}\Z\times \Z} (1+|k|^2)^{m}|\hat v_k|^2<\infty\},
}
where   $\hat v_k=\spr{v,e^{ik\cdot (x,\bphi)}}$. 
Clearly, $H^{m_1}_\per{\subset}H^{m_2}_\per$ for $m_1{\geq}m_2$, and 
$[H^{m_1}_\per,H^{m_2}_\per]_\theta{=}H^{(1-\theta)m_1+\theta m_2}_\per$ 
for the associated real interpolation spaces, see \cite{LM72} and 
\cite[\S2.4.2]{tr78}, and $\nabla\in \CL(H^m_\per,H^{m-1}_\per)$. 

In the following, for $m>1$ we restrict the Sobolev spaces to open subsets  
\hueq{\label{resob} 
V^m_\eps:=H^m_\per(\CC_L(r))\cap \{|u|<r-\eps\} 
}
with some $0<\eps<r$. Then $\CC_L(r)-uN$ is an embedding of  $\CC_L(r)$, 
and we can easily state local Lipschitz properties of the operators 
involved in \reff{hflowu}. For convenience
we recall some Sobolev embedding results, which hold for compact manifolds 
with boundary as for compact domains in $\R^n$, see, e.g., \cite{Aubin82}.%
\footnote{The situation of non compact manifolds is more complicated, but if 
the manifold is complete, has a positive injection radius and bounded sectorial curvature, then the same results holds as well. In fact, in the following 
analysis near $\CC_L(r)$ we do not strictly need the manifold versions 
of Sobolev embeddings due to 
the global coordinates \reff{cylc} also underlying \reff{hmdef}, and hence 
we basically state the more general embedding results for possible 
generalization to other manifolds.} 
In 2D, $H^m(X)\hookrightarrow L^\infty(X)$ 
for $m>1$ embeds continuously, while $H^1(X)$ does not. Also, for $m=1+\al$, 
$H^m(X) \hookrightarrow C^{0,\al}(X)$, 
hence the point evaluation $u(x,\phi)$ and in particular 
the sets $\{|u(x,\bp)|<r-\eps\}$ are well defined for $m>1$. 

\subsubsection{Well-posedness of the PDE part}
Following Remark \ref{prrem} we first consider the PDE part (\ref{hflowu}a) 
with fixed $\Lam\in\R^2$, i.e., 
\huga{\label{qlf1}
\eta(u)\pa_t u=B(u)u+C(u), 
} 
where $B(u)$ collects the 4th order derivatives $b_{\al_1,\al_2}(u)\pa_x^{\al_1}\pa_\phi^{\al_2}$, $|\al|:=\al_1+\al_2=4$, namely 
\huga{\label{psB} 
\barr{l}
b_{4,0}(u)=g^{-2}((r{-}u)^2+u_\phi^2)^2, \quad 
b_{3,1}(u)=-g^{-2}(4u_xu_\phi((r+u)^2+u_\phi^2), \\
b_{2,2}(u)=g^{-2}(2((r{-}u)^2+u_\phi^2)(1+u_\phi^2)+4u_x^2u_\phi^2,\\
b_{1,3}(u)=-4g^{-2}(u_xu_\phi(1+u_x^2), \quad b_{0,4}=g^{-2}(1+u_x^2)^2, 
\earr
}
cf.\cite{LeSi16}, where our $B(u)$ appears as the principal part 
of the 4th order operator $P(u)$ in \cite{LeSi16}. 

\bthm \label{exthm}
Let $m\in (3,4)$ and fix $\Lam\in\R^2$. Then for all $u_0\in V^m_\eps$ 
there exists a $t_0=t_0(u_0)$ and a unique solution 
$u(\cdot,u_0)\in L^2((0,t_0),H^4(\Co))\cap H^1((0,t_0),L^2(\Co))\cap C([0,t_0],V_{\eps}^m))$ of \reff{qlf1}. 
\ethm 

To prove Theorem \ref{exthm} we apply \cite[Theorem 3.1]{MRS13}, 
which treats abstract 
quasilinear parabolic problems of the form 
\huga{\label{aa1}
\pa_t u=A(u)u+F(u),\quad t>0,\quad u(0)=u_0\in \CX, 
}
with Gelfand triple $X_0\hra \CX\hra X_1=D(A(u))$, i.e., phase space $\CX$. 
For this we could rewrite \reff{qlf1} as $\pa_t u{=}\wit B(u)u{+}\wit C(u)$ 
with $\wit B(u){=}\eta(u)^{-1}B(u)$ and $\wit C(u){=}\eta(u)^{-1}C(u)$, but 
since $\eta(u)=(1+u)g^{-1/2}$ contains at most first derivatives 
it will be clear from 
the following that we can as well keep $\eta(u)$ in \reff{qlf1} 
on the lhs and consider the rhs of \reff{qlf1} in the setting of 
\reff{aa1}. As already indicated above, we choose the 
base space $X_0=L^2$ and hence $X_1=D(B(u))=H^4$, 
and we need to find an open subset $V=V^m_\eps$ of the phase 
space $\CX$ with 
$X_0\subset \CX\subset X_1$ such that\\[-3mm]
\bci 
\item[(W1)] $B:V\to L(X_1,X_0)$ and $C:V\to X_0$ are locally Lipschitz; 
\item[(W2)] for all $u\in V$, $B(u)$ generates a $C_0$ (strongly continuous) 
analytic semigroup in $X_0$. 
\eci

By \cite[Theorem 3.1]{MRS13}, \reff{qlf1} is then 
locally well posed in $V$, i.e., there exists a 
$t_0>0$ and a unique solution $u\in C([0,t_0), V)\cap H^1((0,t_0),X_0)\cap 
L^2((0,t_0),H^4)$. 
To show (W1) we use Sobolev embeddings and the structures of $B$ and $C$, 
while (W2) essentially follows from $B(u)$ being uniformly elliptic 
(in the 4th order sense). In detail, 
the symbol $\sig[B(u)](x,\phi,\xi)$ satisfies 
\huga{\label{uell1}
\sig[B(u)](x,\phi,\xi)\ge g^{-2}((r-u)^2\xi_1^2+\xi_2^2)^2,
}
and hence for $\|u\|_\infty<r-\eps$, 
there exists constants $c_1,c_2$ such that 
$c_1\leq \sigma[B(u)](x,\phi,\xi)\leq c_2$
for all $(x,\phi,\xi)\in\CC_L\times \R^2$ with $\|\xi\|=1$.


\blem \label{regnlin}
For $m\in (3,4)$, $B(\cdot)\in C^\infty(V^m_\eps,\CL(H^4_\per,L^2))$ and $C\in C^\infty(V^m_\eps,L^2)$. In addition, the rhs of \reff{hflowu} is an analytic map from  $V^4_\eps\times\R^2$ 
into $L^2$. 
\elem
\begin{proof}
From the Sobolev embeddings, a function $F(u)=u^j$ is well defined for 
any $j$  
on $H^m$ if $m>1$, and the map $u\mapsto g(u)$ is smooth from 
$H^m(X)$ to $H^{m-1}(X)$ for $m>2$, and on $V^m$ bounded away from zero. 
Hence $g^{-s}\in L^\infty(X)$ for $s>0$ and due to the embedding 
$L^\infty(X)\hookrightarrow L^2(X)$, these terms are bounded in $L^2$. 
For $u$ small enough $g^{-s}$ can be expanded into a convergent power series, 
which gives the analyticity. 

The coefficients \reff{psB} of $B(u)$ are well defined 
for $u\in V^m_\eps$ with $m>3$ since $u,\pa_x u,\pa_\phi u\in L^\infty$. 
Actually, $m>2$ is enough here, and the condition $m>3$ comes from 
$C(u)$, which has summands 
$c_\al(u,\nabla u,\nabla^2u)=\tilde c_\al\, u^{j_1}(\pa_xu)^{j_2}(\pa_\bphi u)^{j_3} (\pa_x^2 u)^{j_4} (\pa_\bphi^2 u)^{j_5} (\pa_x\pa_\bphi u)^{j_6} g^{-s}(u,\nabla u),$ 
(nonlinear terms up to second order), and  
$
d_\al(u,\nabla u)\pa_{x,\phi}^\al u, |\al|=3,$ 
where $d_\al$ is nonlinear but the third order terms 
(from $\Delta H(u)$ with $H(u)$ from \reff{meanu}) only appear linearly.   
The maps $b_\al(u,y)$ are analytic on $\{|u|<r-\eps\}\times \R^2$, 
and so are the $c_\al(u,y,z)$ on 
$\{|u|<r-\eps\}\times \R^2\times \R^4$, and 
the $d_\al$ on $\{|u|<r-\eps\}\times \R^2$. 
By introducing a perturbation $u+\del\, v$ 
for $v\in H^\al$, and computing the derivative wrt $\del$, we have, e.g., 
\hueqst{\pa_u 
c_\al(u,\nabla u,\nabla^2 u)=
( \pa_x c_\al(x,y,z) v +\pa_y c_\al(x,y,z)\cdot \nabla v +\pa_z 
c_\al(x,y,z)\cdot \nabla^2 v) |_{(x,y,z)
=(u,\nabla u,\nabla^2 u)}.  
}
The structure of $\pa_u c_\al$ 
is hence the same as that of $c_\al$, and similar for the $b_\al$ and 
$d_\al$ such that these "coefficients" respectively summands 
are well defined and finite in $L^\infty$. 
This can be iterated for higher derivatives, to find the same result, 
and, e.g., for $w\in H^4$ and $u\in V^m_\eps$,
\hueqst{ \|\pa_u B(u)w\|_{L^2} \leq C \|w\|_{H^4}.}
This gives the first two assertions, and setting $u=w$ 
gives the analyticity of $G(u,\Lam)$, 
where the analyticity with respect to $\Lam$ is clear since $\Lam$ appears 
linearly. 
\end{proof} 

To show (W2) we use standard theory for sectorial operators 
\cite{Henry81,luna95}. Let $u\in V^m_\eps$. 
Since $D(B(u))=H^4$ is dense in $X=L^2$, it is sufficient to 
show that $B(u)$ is sectorial in $L^2$, i.e., that the spectrum 
of $B(u)$ is contained in a sector $S_{\al,\om}$ 
in $\C$ with tip at some $\om\in\R$ 
of angle $\al<\pi$ opening to the left, and that there exists and $M>0$ 
such that the resolvent estimate 
\hueq{\label{reso1} 
\|(\lam -B(u))^{-1}\|\leq \frac M{|\lam-\om|}
}
holds for all $\lam\not\in S_{\al,\om}$. 
Essentially, this follows from the uniform ellipticity \reff{uell1}. 
In more detail, we show that $B(u)$  induces a coercive 
sesquilinear form on $H^2$, from which (W2) follows via the Lax--Milgram 
Lemma, see, e.g., \cite[Corollary 2.29]{SE24}. 

\blem \label{coerlem} 
For $u\in V^m_\eps$ and $\lam$ sufficiently large, 
the operator $(\lam-B(u))$ induces a coercive 
sesquilinear form on $H^2_\per(\CC_L(r))$, i.e., 
there exists $\om,C$ independent of 
$u\in V^m_\eps$ such that for $\lam>\om$ we have 
\hueq{\label{coer1} \int_{\CC_r(L)} ((\lam w-B(u)w)w\dd S\geq C 
\|w\|^2_{H^2_\per (\CC_r(L))}.}
\elem 
\begin{proof} 
The idea is to use \reff{uell1}, 
but since $B(u)$ is not in divergence form we need 
to apply the standard trick of writing, e.g., 
$$b_{4,0}(u)\pa_x^4 w=\pa_x^2(b_{4,0}(u)\pa_x^2 w)-\pa_x(\wit b_{4,0}(u)\pa_x^2 w),
$$ 
where $\wit b_{4,0}(u)=\pa_x b_{4,0}(u)$ is still bounded in $L^\infty$, 
i.e., $\|\wit b_{4,0}(u)\|_\infty \le C_{4,0}$. 
Therefore, using $\ds |w_{xx}w_x|\le \frac {\del_1} 2 w_{xx}^2+\frac 1 {2\del_1}
w_x^2$, $\del_1>0$, and similarly 
standard $L^2$--interpolation estimates such as 
$\ds \int_{\CC_L(r)} w_x^2\dd S
\le \int_{\CC_L(r)}\frac {\del_2} 2 w_{xx}^2+\frac 1 {2\del_2} w^2\dd S$, 
we obtain, e.g., 
\hualst{
\spr{(\lam-b_{4,0}(u))w,w}&=\lam\int w^2\dd S+\int b_{4,0}(u)w_{xx}^2+\wit b_{4,0}(u)w_{xx} w_x \dd S\\
&\ge \left( c_1-\frac 1 2(\del_1+\frac{\del_2}{\del_1})\right)
\int w_{xx}^2\dd S
+\left(\lam-\frac {C_{4,0}} {4\del_1\del_2}\right)\int w^2\dd S, 
}
with $c_1>0$ from after \reff{uell1}. Proceeding similarly for 
the other terms from \reff{psB} and choosing suitable $\del_{1,2}$ and 
then $\om$  sufficiently large yields \reff{coer1},   uniformly in 
$u\in V^m_\eps$. 
\end{proof}

In summary we have verified (W1) and (W2) with $V=V^m_\eps$ 
with $m\in(3,4)$ corresponding to $\vt\in(3/4,1)$ for the interpolation 
space $[L^2,H^4]_{\vt,2}=H^m$, and Theorem \ref{exthm} 
now folows from \cite[Theorem 3.1]{MRS13}. 

\subsubsection{The constrained flows} 
In order to take the constraints into account, we 
first consider only one active constraint namely $q_1(u)=0$, and fix $\lam_2$, 
while the same arguments also apply for the case of active 
$q_2$ with fixed $\lam_1$. Defining 
$\overline H =\frac{H}{\|H\|}_{L^2}$ there exists a $\sig_1$ such that 
$-(-\pa_uE+\lam_2\pa_u\CV +\lam_1 H)= -(-\pa_uE+\lam_2\pa_u\CV +\sig_1 \overline H)$, 
and we define the projection $\CP_1:L^2(\CC_r(L))\ra 
\mathrm{span}\{\pa_u \CA(u)\}^\bot$ by 
$\CP_1 f= f- \spr{\overline H,f} \overline H$, 
and consider 
\huga{\label{1cfl}
V_t=\CP_1 (\pa_uE - \lam_2\pa_u\CV). 
} 
Then clearly 
$\ddt \CA(u)= 2\int_{\CC_r(L)} H(u) V_t \dd S =0$ and for $u\in V^m_\eps$ 
we find 
\hualst{ \lam_1&= \frac {1} {\int_{\CC_r(L)}H^2\dd S} \spr{\pa_uE{-}\lam_2\pa_u\CV,2H}\\
&=\frac {1} {\int_{\CC_r(L)}H^2\dd S} \int_{\CC_r(L)} (\Delta H{+}2 H(H^2{-}K){+}c_0K{-}2c_0H{-}\lam_2\pa_u\CV) H\dd S\\
&=\frac {1} {\int_{\CC_r(L)}H^2\dd S} \int_{\CC_r(L)} -|\nabla H|^2+2 H^2(H^2{-}K)+c_0HK-2c_0H^2-\lam_2\pa_u\CV H\dd S.
}
Importantly, the term $|\nabla H|^2$ includes third order 
derivatives of $u$ but only quadratically with 
lower order ``coefficients'' and hence $\int |\nabla H(u)|^2\dd S$
is well defined and analytic for 
$u\in H^m_\per$ for $m>3$. This yields the analogous result as 
Theorem \ref{exthm}. 

\bthm \label{exthm1} For $u_0\in V^m_{\eps}$  
there exists a $t_0=t_0(u_0)$ and a unique solution 
\hueqst{(u(\cdot,u_0), \lam_1)\in 
(L^2((0,t_0),H^4(\Co))\cap H^1((0,t_0),L^2(\Co))\cap C([0,t_0),V_{\eps}^s)))\times C([0,t_0),\R^1).
} 
of $\eta(u)\pa_t u=G(u,(\lam_1,\lam_2))$ under the constraint 
$q_1(u)=\CA(\CC_L-uN)-\CA_0=0$. 
\ethm 

If $\det \CD(u) \neq 0$, then we can consider the full constraints 
in \reff{hflowu} by projecting $V_t=\pa_u E$ onto the 
area and volume preserving flow. 
We construct an orthonormal basis 
 of $\mathrm{span}\{\pa_u\CA,\pa_u\CV\}=\mathrm{span}\{H,(r+u)g^{-1/2}(u) \}$
by Gram-Schmidt as 
\hualst{ \overline {\pa_u\CV} &=\frac 1{\|\pa_u \CV\|_{L^2}} \pa_u\CV= \frac 1{\|(r+u)g^{-1/2}(u)\|_{L^2}} (r+u)g^{-1/2}(u) \\
\overline {\pa_u\CA }&= \frac{1}{ \| H-\spr{\overline{\pa_u\CV}, H}\|_{L^2}} (H-\spr{\overline{\pa_u\CV}, H}) = \frac 1 {\sqrt{\|H\|_{L^2}^2-\spr{\overline{\pa_u\CV},H}^2}}\left(H-\spr{\overline{\pa_u\CV}, H}\right).}
We define the projection 
\hueqst{\CP_{2} f= f-\spr{\overline {\pa_u\CV},f}\overline {\pa_u\CV}-\spr{\overline {\pa_u\CA},f}\overline {\pa_u\CA},} 
and obtain $\ddt \CA=\spr{\pa_u\CA,V_t}=0$ and $\ddt \CV=\spr{\pa_u\CV,V_t}=0$ 
for the flow $V_t=\CP_{2} \pa_u E$. 
Since $\det \CD(u)\neq0$, \reff{eqlam} is uniquely solvable also in a small neighborhood of $u$. So by inverting $\CD(u)$ the Lagrange multipliers are 
\hualst{\lam_1&=\frac {1} {\det \CD}\big( \spr{\pa_u\CV,\pa_u\CV}\spr{\pa_u E,\pa_u\CA} -\spr{\pa_u\CA, \pa_u\CV}\spr{\pa_uE, \pa_u\CV}\big)\\
\lam_2&= \frac {1} {\det \CD} \big(\spr{\pa_u\CV,\pa_u\CV}\spr{\pa_u E,\pa_u\CV} -\spr{\pa_u\CA,\pa_u\CV} \spr{\pa_u E,\pa_u\CA}\big)}
Recall that $\pa_u\CA=2H(u)\eta(u)$, that $\pa_u\CV$ is of first order in $u$, 
and that the leading term in $\pa_uE$ is $\Delta H(u)$. 
Using integration by parts we 
find that the only term involving third order derivatives is $| \nabla H|^2$, 
which was already discussed. 
Hence we obtain the following result. 
\bthm \label{exthm2} For $u_0\in V^m_{\eps}$  
there exists a $t_0=t_0(u_0)$ and a unique solution of \reff{hflowu}, 
\hueqst{(u(\cdot,u_0), \Lam(\cdot, u_0))\in (L^2((0,t_0),H^4(\Co))\cap H^1((0,t_0),L^2(\Co))\cap C([0,t_0),V_{\eps}^s)))\times C([0,t_0),\R^2).
} 
\ethm 
\section{Stability of $\CC_L$}\label{stcy}
We consider the linearization of \reff{hflow} around $\CC_L(r)$. 
This requires some lengthy computations, but since $\CC_L(r)$ has 
constant mean and Gaussian curvatures (with $K=0$), we can explicitly 
compute the 
dispersion relation (Lemma \ref{leig}) and subsequently the stability ranges 
of and bifurcation points from $\CC_L(r)$ (Corollary \ref{destabcor}).
\subsection{Linearization and eigenvalues}\label{lineig}
Since $G(\cdot,\Lam):V^4_\eps\to L^2$ is an analytic map, 
we use \reff{meanu}, \reff{gaussu} and \reff{LBu} to calculate the 
Frech\'et derivative via the Gateaux differential in direction $v$, i.e., 
$\pa_u G(0,\Lam)v=\frac{\dd}{\dd \eps} G(\eps v,\Lam)|_{\eps=0}$; recall that  the Laplace-Beltrami operator on $\CC_L((r)$ is $\Delta_0=\pa_x^2+\frac 1 {r^2} \pa_\bphi^2$ in our coordinate system.
\blem \label{dcomp} The linearizations of 
$H(u), K(u)$, and $\Delta_uH(u)$ around $\CC_L(r)$, i.e., around $u=0$, 
are given by 
\hual{\label{dmean}\pa_u H(0)v&=-\frac 1 2\left(\Delta_0v +\frac 1 {r^2} v\right),\\
\label{dgauss}\pa_uK(0) v&=-r^{-1} \pa_x^2 v, \text{ and }\\[2mm]
\label{dlead}\pa_u \Delta_0H(0)v &= \Delta \pa_u H(0)v.
}
\elem

\begin{proof}
Using \reff{meanu} with $\eps v$ and sorting 
by $\eps$ powers gives
\hueqst{ 2H(\eps v)=-\frac{\eps(r^3 v_{xx}+rv_{\bphi\bphi}) +\CO(\eps^2)}{(r^2+\eps 2rv+\CO(\eps^2))^{3/2}}+\frac{1}{(r^2+\eps 2rv+\CO(\eps^2))^{1/2}}.}
Differentiating w.r.t. $\eps$ this gives 
\hueqst{\frac{\dd}{\dd \eps} 2 H(\eps v)=-\frac{(r^3 v_{xx}+rv_{\bphi\bphi})(r^2+\CO(\eps))^{3/2}+\CO(\eps)}{(r^2+\eps 2rv+\CO(\eps^2))^{3}}-\frac{r v+\CO(\eps)}{(r^2+\CO(\eps))^{3/2}}.}
The denominators are bounded away from zero, hence evaluating at $\eps=0$ gives \reff{dmean}.

For $K$ this works the same way, with 
\hueqst{ \frac{\dd}{\dd \eps} K(\eps v)= \frac{\dd}{\dd \eps}  \frac{-\eps r^3 v_{xx} +\CO(\eps^2)}{(r^2+2\eps 2rv+\CO(\eps^2))^2}=\frac{-r^3 v_{xx} ((r^2+\eps 2rv+O(\eps^2))^{2} +\CO(\eps)}{(r^2+\eps 2rv+\CO(\eps^2))^{4}}.}
Again evaluating at $\eps=0$ gives \reff{dgauss}. 

Finally for \reff{dlead} it is sufficient to notice that 
\hueqst{\frac{\dd}{\dd \eps} \Delta_{\eps v}H(\eps v)|_{\eps=0} 
= \left(\frac{\dd}{\dd \eps}\Delta_{\eps v} H(0)\right)\bigg|_{\eps=0} +\Delta_0  \left(\frac{\dd}{\dd \eps}H(\eps v)\bigg|_{\eps=0}\right)
=\Delta_0  \left(\frac{\dd}{\dd \eps}H(\eps v)\bigg|_{\eps=0}\right),}
since $H(0)$ is a constant and $\Delta_{\eps v}$ is non degenerate for small perturbations. 
\end{proof}

Inserting the expressions from Theorem \ref{dcomp} we obtain 
\blem $G(\cdot,\Lam)$ is Frechet differentiable at $u=0$, with differential, 
using \reff{socon}, 
\hual{\label{Gjac} \pa_uG(0,\Lam) v= -\frac 1 2\Delta^2_0 v-\left(r\lam_2+\frac{1}{r^2}\right) \Delta_0 v- \left(\frac{\lam_2}{r} +\frac{1}{2r^4}\right) v +\left(\frac 1 r-2c_0\right) r^{-1} \pa_x^2 v.
}
\elem
Over the unbounded $\CC(r)$, $\pa_uG(0,\Lam))$ has a continuous spectrum, 
but by restricting to the compact cylinder $\CC_L(r)$ of length $L$, 
 $\pa_uG(0,\Lam))$ has a discrete spectrum. 
\blem \label{leig}
The eigenvalue problem 
\hueq{ \label{eigeq}\pa_uG(0,\Lam)\psi=\mu\psi}
has a countable set of solutions $(\mu_{mn},\psi_{mn})\in 
\R\times H^4_{\per}(\CC_L(r))$,  $(m,n)\in\Z^2$, 
\bsub\label{eigs}
\hual{
&\text{$\psi_{mn}(x,\bphi)=e^{\ri (k_mx+n\bphi)}$, 
$k_m=2m\pi/L$, and $\mu_{mn}=\mu(k_m,n)$ with}\\
&\mu(k,n,\Lam)=-\frac 1 2  \Big[k^2{+}\frac {n^2} {r^2} \Big]^2 
+\left(r\lam_2{+}\frac{1}{r^2}\right)\Big[k^2{+}\frac {n^2} {r^2} \Big]
-\left(\frac{\lam_2}{r}{+}\frac{1}{2r^4}\right) 
-\left(\frac 1 r{-}2c_0\right) \frac{k^2}{r}.
}
\esub 
For all $\Lam\in\R^2$ the eigenfunctions $\psi_{0\pm1} \in\ker(\pa_uG(0,\Lam))$.
\elem
\begin{proof}
Recall that $H^4_{\per}(\CC_L(r))$ has a compact embedding in $L^2(\CC_L(r))$, 
hence the resolvent of \reff{Gjac} is compact. Using the $L^2$ basis of 
eigenfunctions of $\Delta_0$, 
\hueqst{ \psi_{mn}(x,\bphi)=e^{\ri(k_mx+n\bphi)}, \quad 
\Delta_0\psi_{mn}= -\Big[k_m^2+\frac 1 {r^2} n^2\Big]\psi_{mn}(x,\bphi), 
}
we obtain (\ref{eigs}b). 
Clearly $\mu(0,\pm1,\Lam)=0$, hence 
$\sin(\bphi),\cos(\bphi)\in \ker(\pa_uG(0,\Lam))$. 
\end{proof}

Using Lemma \ref{aderlem}, we write the linearization of \reff{hflowu} as 
\hual{ \label{genev0}\pa_U \sG(0,\Lam) U
& =\bpm \pa_u G(0,\Lam ) & \pa_{\lam_1} G(0,\Lam ) &\pa_{\lam_2} G(0,\Lam ) \\ 
		\pa_u \CA(0 ) & \pa_{\lam_1} \CA(0 ) &\pa_{\lam_2} \CA(0 )\\
		\pa_u \CV(0 ) & \pa_{\lam_1} \CV(0 ) &\pa_{\lam_2} \CV(0 )
		 \epm  \bpm u\\ \sig_1\\\sig_2\epm\nonumber\\
	&=\bpm \pa_u G(0,\Lam )u & -2H(0) \sig_1&-\sig_2 \\ 
		-\int_{\CC_L(r)} 2H(0) u\dd S &0 & 0\\
		\int_{\CC_L(r)} u \dd S& 0 & 0
\epm.}
The domain of $\pa_U\CG(0,\Lam)$ is $H^4_\per(\CC_L(r))\times \R^2$,  
equipped with the norm $\|U\|_{\pa_{(u,\Lam)}\sG}=\|u\|_{H^4}+\|\Sigma\|_{\R^2}$, and altogether we have the generalized eigenvalue problem 
\hual{ \label{genev}\pa_U \sG(0,\Lam) \Psi
= \mu \bpm 1 & 0& 0\\0 & 0& 0\\0 & 0& 0 \epm \Psi,\quad \Psi=(\psi,\sig_1,\sig_2),
}
with eigenvalues $\mu \in \{ \mu(k_m,n), \, (m,n)\in\Z^2\backslash \{(0,0)\}\cup \{0\}$ and eigenfunctions 
\hueqst{\Psi_{mn}= \bpm \psi_{mn}\\ 1\\ -2\epm, \, (m,n)\in\Z^2\backslash \{(0,0)\}, \text{ and }\Psi_{00}=\bpm 0\\ 1\\ -2\epm.
} 
\brem (a) The eigenvalues $\mu(0,\pm 1,\lam)=0$ correspond to 
translations of $\CC_L(r)$ in $y$ and $z$ directions. 
 Each eigenvalue of the form $\mu(k_m,0,\Lam)$ and $\mu(0,n,\Lam)$ has 
multiplicity two and every ``double'' eigenvalue $\mu(k_m,n,\Lam)$ 
with $m\neq0\ne n$ has multiplicity four with eigenfunctions $\psi_{\pm m\pm n}$.

(b) The eigenvalues $\mu(k_m,n,\Lam)$ do not explicitly depend on $\lam_1$, but implicitly. 
Recall that $\lam_1=\frac 1 2 \left( \frac 1 {2r^2} -2r\lam_2-2c_0^2\right)$, hence the parameters are not independent on the cylinder. 

(c) In the following we also use the notation $\Psi_{m,n}$ for the first 
two components of $\Psi_{m,n}$, which are the eigenvectors for the case 
of just the area constraint being active, cf.~Theorem \ref{exthm1} and 
\reff{areahelf} in \S\ref{cmr}. 

(d) The crucial difference of the spectrum of \reff{genev} and just the PDE part 
\reff{eigeq} is that $(\psi_{00},1,-2)$ with $\psi_{00}\equiv 1$ 
is not an eigenvector of $\pa_u\sG (0,\Lam)$. 
\eex\erem

\subsection{Center manifold construction}\label{cmr}
Based on the local existence Theorem \ref{exthm1} and the spectral 
computations from Lemma \ref{leig} we construct the center manifold 
for the flow near $\CC_L(r)$ which will be the basis for the subsequent 
stability and bifurcation results. 
At $\CC_L(r)$, $\det \CD=0$, and we cannot perturb $\CC_L(r)$ 
while preserving the area and the volume. Therefore we start with 
\hueq{\label{areahelf} \bpm \pa_t u\\ 0\epm = 
\bpm \eta^{-1}(u)G(u,\Lam)\\q_1(u,\Lam)\epm=:\CG_\CA(u,\Lam),}
where $\Lam=(\lam_1,\lam_2)$ with $\lam_2$ a {\em fixed} external 
pressure difference. 
The $\R^3$ has the Euclidean symmetry group $E(3)$ of rigid body motions, 
which can be constructed using two smooth curves 
$\tau,\brho:\R\ra\R^3$ with $\tau(0)=\brho(0)=0$, such that the translated and rotated surface is given by 
\hueqst{  
\tilde X=\{y+\tau(\eta)+y\times \brho(\eta): y\in X\}.
}
However, since we only consider normal variations of $X$, 
admissible rigid body motions are 
\hueqst{\CM X=\{X+u\, N: u =\spr{ \tau(\eta)+y\times \brho(\eta),N(y)}\}.} 
The dimension of $\CM X$ is given by the dimension of the tangent space 
$T_e \CM X$ at the identity, 
hence spanned by $\spr{\tau'(0),N(y)}$ and $\spr{y\times \brho'(0),N(y)}$. 
Since $\tau'(0)$ and $\brho'(0)$ are arbitrary vectors, 
\hueq{\label{ridmo}T_e\CM X=\mathrm{span} \{\spr{e_i,N},\spr{y\times e_i,N},\, i=1,2,3\}.} 
But for $\CC_L(r)$ we have 
\hueq{ \label{tcmfc} T_e \CM\Cr=\mathrm{span}\{ \spr{e_2,N}, \spr{e_3,N}\}=\mathrm{span}\{ \psi_{01}, \psi_{0-1}\}
}
where a straightforward calculation shows that translation in and 
rotation around $e_1$ generate tangential motions to $\CC_L(r)$, while 
rotations around $e_2$ and $e_3$ violate the periodicity in $x$. 

In \S\ref{stabrange} we will see that there are parameter values 
$\lam_2$ such that the spectrum of the linearization of \reff{areahelf} 
is contained in the negative half plane, with one zero eigenvalue 
with multiplicity according to \reff{tcmfc}. The spectrum on the finite 
cylinder is discrete, with a spectral gap between the zero and first non zero 
eigenvalue. Hence the spectrum decomposes into the center 
spectrum $\mathrm{spec}_c$ and the stable spectrum $\mathrm{spec}_s$. 
Using the well-posedness Theorem \ref{exthm1}  we may apply 
\cite[Theorem 4.1]{Si95} on the existence of a center manifold under 
the above conditions. Since \reff{tcmfc} shows that the center eigenspace 
coincides with the tangent space 
of the manifold of translated cylinders, both manifolds also coincide locally. 
This gives the following results, 
where again $y$ and $z$ are 
the euclidian coordinates orthogonal to the cylinder axis $x$. 
\blem\label{cmf}
At the trivial solution $(0,\lam_2)$ with $\lam_1$ and $\lam_2$ related by \reff{socon}, there exists a smooth center manifold $\mathcal M^c$ for \reff{areahelf} 
of dimension at least two. If the kernel of $\pa_{(u,\lam_2)}\CG_\CA$ 
is spanned by $\Psi_{0\pm1}$ then 
\hueq{\CM^c = \text{Shifts of $\CC_L$ in $y$ and $z$ direction}.}
\elem

\bthm\label{stabthm}
If $\lam_2\in\R$ such that the spectrum of $\pa_{(u,\lam_2)} \CG_\CA$ 
is contained in $(-\infty,0]$ and $\ker(\pa_{(u,\lam_2)} \CG_\CA)$ is spanned 
by $\Psi_{0\pm1}$, then there exists $\del,\om,M>0$ such that for all 
$u_0\in V^m_\eps$ with $\|u_0\|_{H^m}<\del$  
 the solution $U=(u,\lam_1)$ of \reff{areahelf} exists globally in time and 
there exist a unique $U_\infty\in \CM^c$ such that 
\hueqst{ \| U(t,U_0)-U_\infty\|_{H^m\times\R}\le Me^{-\om t}\|U_0\|_{H^m\times\R}.
} 
\ethm

\subsection{Stability and instability} \label{stabrange}
Following Theorem \ref{stabthm}, to check for stability of $\CC_L(r)$ 
we seek $\Lam\in\R^2$  such that all 
eigenvalues $\mu(k_m,n,\Lam)$ are in $(-\infty, 0]$, and the only 
zero eigenvalues belong to $\Psi_{0\pm1}$ (translations in $y$ and $z$).  
{\em From here on, to simplify notation we fix $r=1$}, and recall 
that the results for general $r$ follow from rescaling $r \mapsto 1$, 
$L\mapsto rL$ and $c_0\mapsto c_0/r$. 
Practically,  we vary $\lam_2$, with $\lam_1$ given by \reff{socon} 
along the trivial branch, and monitor the signs of eigenvalues. 
Depending on $c_0,L$, the main questions are:
\bci
\item Are there intervals for $\lam_2$ such that $\CC_L$ is linearly stable?
\item Which modes destabilize the cylinder at what $\lam_2$?
\eci 

\begin{figure}[h]
\bce
\btab{l}{{\sm (a)}\hs{49mm}{\sm (b)}\hs{49mm}{\sm (c)}\\
\hs{-2mm}\ig[width=0.35\tew]{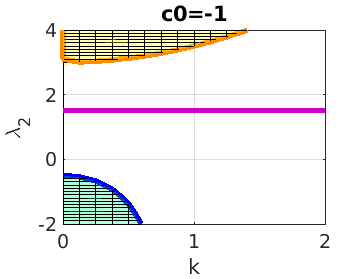}
 \hs{-2mm}\ig[width=0.35\tew]{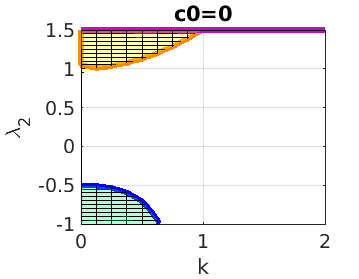}
 \hs{-2mm}\ig[width=0.35\tew]{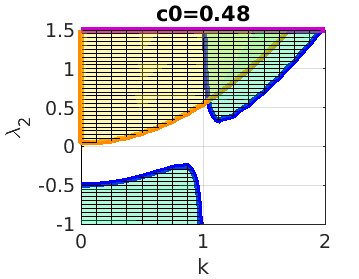}\\
{\sm (d)}\hs{65mm}{\sm (e)}\\
\ig[width=0.22\tew]{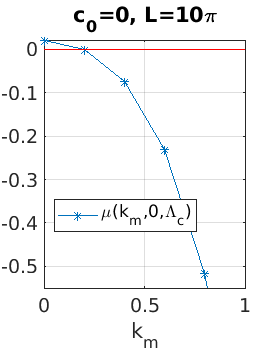}
\hs{0mm}\ig[width=0.22\tew]{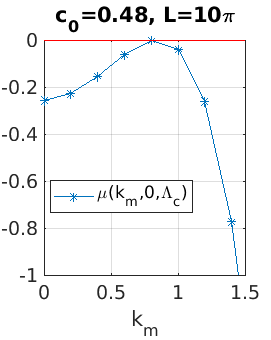}
\hs{-2mm}\ig[width=0.35\tew]{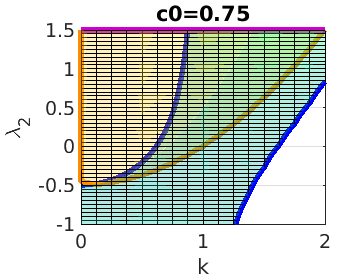}\hs{-4mm}
\hs{4mm}\ig[width=0.22\tew]{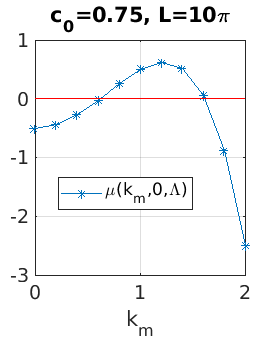}
}
\ece
\vs{-2mm}
\caption{\small{(a--c) Instability curves in the $k$--$\lam_2$ plane 
($\lam_1$ according to \reff{socon}) 
for different values of $c_0$, delimiting the $\lam_2$ parameter range 
for stability of $\CC_L$. 
The blue curve corresponds to $(k,0)$ (pearling), 
and the orange to $(k,1)$ (coiling/buckling, for $k\ne 0$). 
The respective unstable sets $(k,\lam_2)$ in like colors, 
where the grid lines in $k$ correspond to $k_m$ for $L=16\pi$.  
The violet line $\lam_2=3/2$ illustrates the wrinkling 
(mode $(0,2)$). For given $(c_0,\lam_2^*)$, $\CC_L$ is stable 
for all $L$ if the line $\lam_2=\lam_2^*$ is below $\lam_2=3/2$ and 
entirely in the white set. The lower and upper blue curves 
from (c) meet at $k=1$ for $c_0=1/2$ beyond which no stable $\lam_2$ range for 
$\CC_L$ remains. The maximum  near $k=1$ in the lower blue curve in 
(c) gives the only finite wave number instability in these plots. 
(d) Selected eigenvalue curves at criticalities. 
(e) Fully unstable regime for $\CC_L$, and 
$k_m\mapsto \mu(k_m,0,\Lam)$ with $\lam_2=0$. 
 }}  
\label{neigs}

\end{figure}

In Lemma \ref{leig}, we fixed $L{>}0$ yielding the discrete allowed 
wave numbers $k_m=2\pi m/L$, $m\in\Z$. Analytically it is more convenient  
to return to continuous $k\in\R$. Thus we consider 
\huga{\label{mumu}
\mu(k,n,\lam_2)
=-\frac 1 2 (k^2+n^2)^2+(\lam_2+1)(k^2+n^2)-(\lam_2+1)-(1-2c_0)k^2, 
}
where due to the reflection (in $x$) 
symmetry of $\pa_u G(0)$ in \reff{Gjac} only even powers of $k$ appear, and 
hence we may restrict to $k\ge 0$, and similarly to $n\ge 0$. 
Subsequently we reintroduce 
$k_m$ to formulate  Corollary \ref{destabcor}. 
See also Fig.\ref{neigs} for illustration of the results. 

For $k=0$ we have 
\huga{\label{azieig}
\mu(0,n,\lam_2)=-\frac{n^4}{2} +\left(\lam_2+1\right)n^2-\left(\lam_2+\frac1 2\right),}
with $\mu(0,1,\lam_2)=0$ for all $\lam_2$ (translational eigenvalues), and 
$\pa_{\lam_2} \mu(0,n,\lam_2)=n^2-1$
is positive for $n>1$. Furthermore $\mu$ is zero at 
$\lam_2^{(0,n)}=(1-n^2)^{-1}(-\frac{n^4}{2} +n^2-\frac 1 2)$, 
which is monotonically increasing in $n$. 
Hence, {\em under axially constant perturbations $\CC_L$ is stable for 
$\lam_2<\frac 3 2$}, 
indicated by the violet line in Fig.\ref{neigs}.

Next we consider $k>0$ and $n\ne 0$, in particular  $n=1$ which 
turns out to be the most dangerous $n$ for $k\ne 0$. We have 
\hualst{\mu(k,1,\lam_2)&=-\frac 1 2 ((k^2+1)^2+(\lam_2+1)(k^2+1) 
-\left(\lam_2+\frac{1}{2}\right)-(1-2c_0)k^2\nonumber \\
&= -\frac 1 2 (k_1^2+1)^2 +(\lam_2+2c_0)k^2+\frac 1 2=-k^2(\frac 1 2 k^2+1-\lam_2-2c_0).
}
One sees directly that $\mu(k,1,\Lam)$ is strictly increasing 
in $\lam_2$ and has a zero at $k>0$ for 
\huga{\label{inmixed} 
\lam_2= \frac 1 2 k^2+1-2c_0,}
indicated by the yellow line in Fig.\ref{neigs}. In particular,  
{\em if $c_0>-1/4$, then there exists an explicitly computable 
$L$ large enough so that the $(k_1,1)$ coiling/buckling 
instability,  $k_1=2\pi/L$, 
in $\lam_2$ is below the wrinkling instability of $\CC_L$.} 

For $k>0$ and $n=0$ (pearling) 
we have, with $\kap:=k^2$, i.e., 
$\muti(\kap,\lam_2):=\mu(k^2,0,\lam_2)$, 
\hual{\label{l2p1} 
&\muti(\kap,\lam_2)=-\frac 1 2 \kap^2+(\lam_2+2c_0)\kap-\lam_2-\frac 1 2, 
}
and solving $\mu(\kap,\lam_2)=0$ for $\lam_2=\lam_2(\kap)$ yields 
\huga{\label{l2hu}
\lam_2(\kap)=\frac 1 2 \frac{\kap^2-4c_0\kap+1}{\kap-1}=
\frac 1 2 (\kap-1)+\frac{(1-2c_0)\kap}{\kap-1}.
}
For $c_0\ne \frac 1 2$, this has a singularity at $\kap=1$ with 
\hugast{
\left\{\barr{ll} \lam_2(\kap)\to -\infty&\text{ if $\kap\nearrow 1$ and 
$c_0<1/2$, or $\kap\searrow 1$ and $c_0>1/2$},\\
\lam_2(\kap)\to \infty&\text{ if $\kap\nearrow 1$ and 
$c_0>1/2$, or $\kap\searrow 1$ and $c_0<1/2$}.
\earr\right.
}
For $c_0=\frac 1 2$ we have $\lam_2(\kap)=-1/2+\kap/2$. In any case, 
$\pa_{\lam_2}\muti(\kap,\lam_2)=\kap-1$ such that the unstable 
region is below $\lam_2(\kap)$ for $0<\kap<1$ and above 
$\lam_2(\kap)$ for $\kap>1$. See the blue curves/sets in Fig.\ref{neigs}. 
Finally, searching for extrema of $\lam_2(\kap)$ 
via $\pa_\kap \lam_2(\kap)=0$ yields 
\huga{\label{kape}
\kap_\pm=1\pm\sqrt{2-4c_0}.
} 
Hence, $\kap_\pm$ exist for $c_0<1/2$, and $\kap_-\in(0,1)$ for 
$c_0\in(1/4,1/2)$, and (only) in this case we {\em have a finite 
wave number $k=\sqrt{\kap}$ instability of $\CC_L$, for sufficiently 
large $L$}. 

Again, see Fig.\ref{neigs} (a)--(d) for illustration of these instability 
curves and selected dispersion relations, \huc{while in (e) we show 
the fully unstable case, which will be of interest for the outlook 
in Fig.\ref{c1f1}.}
In the following Corollary we summarize the 
above, with maximal wavelength refering to the minimal allowed 
$k_1=2\pi/L$. 

\bcor \label{destabcor}
The cylinder $\CC_L$ is volume and area preserving stable for every 
$L>0$ if:\\[-3mm]
\bcena
\item $c_0\in (1/4,1/2)$ and 
$\lam_2\in(\lam_2^{\text{pearl}},\lam_2^{\text{cb}})$ where 
$\lam_2^{\text{pearl}}=\lam_2(\kap_-)$ with $\kap_-=1-\sqrt{2-4c_0}$ 
and $\lam_2$ from \reff{l2hu},  and where $\lam_2^{\text{cb}}=1-2c_0$. 
In particular, $\CC_L$ with large $L$ destabilizes near 
$\lam_2^{\text{pearl}}$ to a finite wavelength 
pearling, and near $\lam_2^{\text{cb}}$ to 
maximal wavelength buckling and coiling. 
\item $c_0\in [-1/4,1/4]$ and $\lam_2\in(\lam_2^{\text{pearl}},\lam_2^{cb})$, 
where $\lam_2^{\text{pearl}}=-1/2$ and $\lam_2^{\text{cb}}=1-2c_0$. 
In particular, $\CC_L$ with large $L$ destabilizes near 
$\lam_2=-1/2$ to maximal wavelength pearling, and near 
$\lam_2^{\text{cb}}$ to maximal wavelength coiling and buckling.  
\item $c_0\in (-\infty,-1/4]$ and $\lam_2\in(\lam_2^{\text{pearl}},
\lam_2^{\text{wri}})$, where $\lam_2^{\text{pearl}}=-1/2$ and 
$\lam_2^{\text{wri}})=3/2$. In particular, 
$\CC_L$ with large $L$ destabilizes near $\lam_2=-1/2$ 
to maximal wavelength pearling, and at $\lam_2=3/2$ to 
wrinkling (homogeneous in $x$). 
\ecen
\vs{2mm}
If $c_0>1/2$, then for all $\lam_2\in\R$ there exists an $L_0>0$ such 
that $\CC_L$ with $L>L_0$ is spectrally unstable. 
\ecor

\huc{
\brem\label{carem}
The values $\lam_2^{\text{cb}}=1-2c_0$ and $\lam_2^\text{pearl}=-1/2$ 
are the asymptotics for $k_1\to 0$, i.e., for the infinitely 
long cylinder. In the following section we return to some finite 
$L>0$, i.e., some finite $k_1=2\pi/L$, and in the numerics in \S\ref{numex} 
we use, e.g., $L=10$, hence $k_1=\pi/5\approx 0.628$, or 
$L=30$, hence $k_1=\pi/15\approx 0.209$. The corresponding values 
for  $\lam_2^{\text{cb}}$ and $\lam_2^\text{pearl}$ (maximum wavelength 
pearling) can then be read off from the orange and blue lines 
in Fig.\ref{neigs}. For, e.g., $c_0=0.48$ and $L=30$, 
the finite wavelength pearling 
corresponds to $k_m=2\pi m/30\approx 0.838$ for $m=4$. 
\eex\erem 
}

\section{Bifurcations from $\CC_L$}\label{aesec}
We now compute the branches bifurcating at the instabilities given 
in Corollary \ref{destabcor}. In \S\ref{ampeq} we give the 
general bifurcation results, based on the $O(2)\times O(2)$ 
symmetry of the problem at $\CC_L$, and in \S\ref{redeq} we derive 
the associated amplitude equations (AEs) and moreover discuss the 
stability of the bifurcating branches.%
\huc{\footnote{We exclude (leave for future work) the exceptional 
(co--dimension two) cases of combinations $(L,c_0)$ such that at 
some $\lam_2$ two different modes become unstable simultaneously. 
This can for instance happen for large $L$ and $c_0\stackrel>\approx 
-0.25$, namely mode crossing of wrinkling and coiling/buckling.}}

\subsection{General bifurcation results} \label{ampeq}
As a bifurcation problem, the desired primary bifurcation parameter is 
either $\CV$ or $\CA$. However, these 
parameters do not change along the trivial $u=0$ branch, and therefore 
we first choose $\lam_2$ as the primary bifurcation parameter and only 
consider \reff{areahelf}, which we write as 
\bsub\label{bifprob} 
\hual{
\eta(u) \pa_t u&= G(u,\Lam),\\ 
0&=\CA(u)-\CA_0.
}
\esub 
The linearization of \reff{bifprob} is the first $2\times 2$ block of 
$\pa_U \CG$, hence with eigenvalues \reff{eigs} with 
$(m,n)\in\Z^2\backslash\{(0,0), (0,\pm 1)\}$, and eigenvectors 
\hueq{
\Psi_{mn}(x,\bphi)=\bpm A_{mn}e^{\ri(k_mx+n\bphi)} + \cc\\ 0 \epm. 
}
Here $\mu_{00}$ is not present since it corresponds to 
changes of radius along the cylinder branch. 
Solving \reff{neigs} gives (possible) bifurcation points 
(for completeness with $r$--dependence) 
\hual{\label{bifp}\lam_2=\Bigg[\frac {(1-n^2)}{r^2}-k_m^2\Bigg]^{-1} \bigg(-\frac 1 2  \Big[k_m^2+\frac  {n^2}{r^2}\Big]^2 +\frac{1}{r^2}\Big[k_m^2+\frac{ n^2}{r^2}\Big]
-\left(\frac 1 r-2c_0\right) \frac{k_m^2}{r}-\frac{1}{2r^4} \bigg).} 
Due to symmetries the bifurcation points have 
multiplicities at least two, and hence the equivariant branching lemma \cite{CL00,GoS2002,Hoyle} gives the needed framework for the bifurcation analysis. 

\bthm \label{bifthm} At \reff{bifp} in the system \reff{bifprob} we have 
equivariant pitchfork bifurcations  
 of steady state solutions $X(\eps)=\CC_L-\eps u(\eps)N$, $\Lam(\eps)=\Lam_0+\CO(\eps^2)$, for $0<|\eps|<\eps_0$ sufficiently small. 
The maps $\eps\mapsto(u(\eps),\Lam(\eps))\in H^4\times\R^2$, or, equivalently 
$\eps\mapsto (X(\eps),\Lam(\eps))\in H^4_\per(\CC_L,\R^3)\times\R^2$ 
are analytic, and we characterize the bifurcations as follows: 
\bcena \item$n=0$ (pearling): 
\hueq{\label{bifaxi}
	u=\left(A_{m}e^{\ri k_mx} + \cc\right)+\CO(\eps) }
	with amplitude $A_m\in\C$; 
the isotropy class can be represented by the isotropy subgroup 
$\Sigma:=\{\sig\in\Gamma: \sig X=X\}=\{
(0,m_x,\tau_\phi,m_\phi)\}\leqslant O(2)\times O(2)$, 
where $m_x$ stands for reflections in $x$ at $x_0$ with maximal or minimal 
$r(x_0)$, $\tau_\phi$ stands for translations  in $\phi$ (rotations), 
and $m_\phi$ for reflections in $\phi$. The 
fixed point subspace ${\rm Fix}\ \Sigma$ consists of surfaces of revolution. 
\item$m=0$ (wrinkling): 
\hueq{ \label{bifazi}
	u= \left(A_{n}e^{\ri n\bphi} + \cc\right)+\CO(\eps) }
	with amplitude $A_n\in\C$; an isotropy subgroup is 
$\Sigma=\{(\tau_x,m_x,0,m_\phi)\}$ (translations and reflections 
in $x$, and reflections in $\phi$). ${\rm Fix}\ \Sigma$ consists of $x$-independent 
surfaces. 
\item$m\neq0\neq n$ (coiling and buckling): Two types of branches bifurcate. 
Coiling, 
 \hueq{\label{bifnaxi} 
	u=\left(A_{mn}e^{\ri (k_mx+n\bphi)} + \cc\right) +\CO(\eps), 
}
with $\Sigma=Z_2$ generated by {\em joint} reflections in $x$ {\em and} 
$\bphi$  (at suitable $x_0, \bphi_0)$, and continuous symmetry 
$\tau_{\phi,-k_m\xi/n}\circ \tau_{x,\xi}$, i.e., spatial translation 
followed by angular rotation. In particular, spatial translation and 
angular rotation yield the same non--trivial group orbits. 

The other branch corresponds to buckling as equal amplitude 
superpositions of coiling, 
 \hueq{\label{bifnaxispi} 
	u=\left(A_{mn}e^{\ri (k_mx+n\bphi)} +A_{m-n}e^{\ri (k_mx-n\bphi)} + \cc\right)+\CO(\eps), } 
with $\Sig=Z_2\times Z_2$ (independent reflections in $x$ or $\bphi$ at 
suitable $x_0, \bphi_0$), and no continous symmetries left.  
\ecen
\ethm

\begin{proof}
The proof proceeds via the equivariant branching lemma; see, for example, 
\cite[Theorem 2.3.2]{CL00}. 
By \reff{uell1}, $B(u)$ is elliptic; hence the 
linearization of $G(u,\Lambda)$ is elliptic as well, 
and thus the upper left $2\times2$ block of $\partial_U \CG$ is a 
Fredholm operator. 
At every $\lambda_2$ from \reff{bifp}, the zero eigenvalue is 
of finite multiplicity, namely multiplicity 2 in (a,b), and multiplicity 4 
in (c), with the symmetries acting on 
the amplitudes $(A_{mn},A_{m-n})\in\C^2$. 
For (c), the joint 
reflection gives the condition 
that either $A_{mn}$ or $A_{m-n}$ must be zero (coiling), while 
enforcing a $Z_2\times Z_2$ symmetry yields $A_{mn}=A_{m-n}$ (buckling). 
In all cases, the isotropy subgroups acting on the amplitudes 
yield one-dimensional ${\rm Fix}(\Sig)$. 
\end{proof}

\brem 
The equivariant group actions correspond to the rigid body motions stated in \reff{ridmo}.
As the solution $X$ on the bifurcating branch has the same symmetry group as 
$\CC -u N$, with $u$ one of the functions from Theorem \ref{bifnaxispi}, we shall compute that, for example, in (a), the translations in $x$ are nontrivial rigid body motions of 
$X - (\eps A_m e^{ik_mx} + \mathrm{c.c.}+\CO(\eps^2))N$. 

Assume that $A_{-m} = A_m$ with $A_m \in \R$, so that $u = \cos(k_mx)$. 
The (not normalized) normal vector is 
\hualst{
\nu = \partial_x X(\varepsilon) \times \partial_\bphi X(\varepsilon)
&= \bpm 1\\ -\eps k_m\sin(k_mx) \cos(\bphi)\\ -\eps k_m\sin(k_mx)\sin(\bphi)\epm 
\times 
\bpm 0\\ -(1+\eps \cos(k_mx)) \sin(\bphi)\\ (1+\eps \cos(k_mx))\cos(\bphi)\epm +\CO(\eps^2)\\
&= -(1+\eps \cos(k_mx)) 
\bpm \eps k_m\sin(k_mx)\\  \cos(\bphi)\\ \sin(\bphi)\epm+\CO(\eps^2).
}
Hence $\langle \nu, e_j \rangle \neq 0$ for $j=1,2,3$; thus, ${\rm span}\{\langle \nu, e_j \rangle, j=1,2,3\}$ is a subspace of the tangent space of the center manifold of $X$. However, as this coincides with the number of equivariant group actions close to the bifurcation point, the tangent space coincides with the 
three translations. In contrast, a rotation around the $x$-axis yields
\hualst{ 
X(\varepsilon) \times e_1 = 
\bpm 0\\ -(1+\eps \cos(k_m x))\sin(\bphi)+\CO(\eps^2)
\\ (1+\eps \cos(k_m x))\cos(\bphi)+\CO(\eps^2)\epm, 
} 
so $\langle X(\varepsilon) \times e_1, \nu \rangle = 0$ and the rigid body motion is trivial.
A similar argument applies for (b). 

In (c), the 4D kernel generates two different types of branches. 
For the buckling, 
a combination of the above arguments yields that translations 
in $x$ and the rotations in $\bphi$ are nontrivial rigid body motions.
For the coiling we have, choosing $u=\cos(k_m x+n\bphi)$, 
\hualst{ 
\partial_x X(\varepsilon) &= \bpm 1\\ 0\\ 0 \epm - \eps k_m\sin(k_mx+n\bphi) \bpm 0\\ \cos(\bphi)\\\sin(\bphi)\epm +\CO(\eps^2),\\
\partial_\bphi X(\varepsilon) &= (1+\eps \cos( k_mx+n\bphi)) \bpm 0\\ -\sin(\bphi)\\ \cos(\bphi)\epm 
- \eps n\sin(k_m x+n\bphi)\bpm 0 \\ \cos(\bphi)\\ \sin(\bphi)\epm +\CO(\eps ^2), 
}
so the normal vector is 
\hualst{
\nu = &(1+\eps \cos( k_mx+n\bphi)) 
\bpm 0\\ -\cos(\bphi)\\ -\sin(\bphi)\epm 
- \eps n\sin(k_m x+n\bphi) \bpm 0\\ -\sin(\bphi)\\ \cos(\bphi)\epm \\
&- \eps k_m\sin(k_mx+n\bphi)(1+\eps \cos(k_mx+n\bphi)) \bpm 1\\0\\0\epm+\CO(\eps^2).
}
This results in a nontrivial translational symmetry in $x$. 
Similarly, 
\hualst{ 
X(\varepsilon) \times e_1 = 
\bpm 0\\ -(1+\cos(k_mx+n\bphi))\sin(\bphi)\\ (1+\cos(k_mx+n\bphi))\cos(\bphi)\epm+\CO(\eps^2),
}
such that $\langle X \times e_1, \nu \rangle = \eps n\sin(k_mx+n\bphi)(1+\eps \cos(k_mx+n\bphi))+\CO(\eps^2)$ 
yields the same action as the translation. 
These computations are also important for the numerics, 
as they give the formulas and numbers of phase conditions needed for the 
numerical continuation in \S\ref{numex}. 
\hfill\mbox{$\rfloor$} 
\erem 

\subsection{Reduced equation} \label{redeq}
To further describe the bifurcations from \autoref{bifthm}, 
we derive the reduced equations (or amplitude equations AEs) at $\CC_L$. 
These are straight forward computations, but somewhat 
cumbersome in our geometric setting and thus we partly use 
{\tt Maple}. Here we explain the basic steps, and 
the computations for the pearling and wrinkling instability, but for 
the coiling vs buckling we relegate some details to Appendix \ref{aeapp}. 
Additionally we give the numerical values of the coefficients of 
the AEs for some of examples, which will later be compared to numerics 
for the full problem. 

According to \autoref{bifthm}, at ``simple'' (strictly speaking of multiplicity 
two, but from here we factor out the translation or rotation) bifurcation 
points \reff{bifp}, i.e., $(m,n){=}(m,0)$ or $(m,n){=}(0,n)$, the 
 bifurcating 
branch is tangential to solutions of 
 $D\CG(0,\Lam)[(\psi,0)] = 0$ where $\CG$ is the rhs of \reff{bifprob}. 
We thus introduce the small parameter 
\hual{ 
&\eps=\sqrt{\left|\frac{\lam_2-\lam_2(m,n)}{\beta_2}\right|},
}
where $\beta_2=\pm 1$ is introduced for later convenience to deal 
with sub--and supercritical pitchforks, and altogether we use the ansatz 
\hual{
&\eps u=\eps A\psi_{m,n} +\eps^2\sum_{(i,j)\in \CB}  
A_{i,j} \psi_{i,j} +\cc +\CO(\eps^3),\\
&\lam_2=\lam_2(m,n)+\beta_2 \eps^2+\CO(\eps^3), \quad 
\lam_1=\lam_1(m,n)+\al_2\eps^2+\CO(\eps^3), \quad 
}
where $\CB$ is the index set of the modes generated from the critical 
modes in quadratic interactions, $A=A(T)$, $A_{i,j}=A_{i,j}(T)$, 
and $T=\eps^2 t$.  
The $O(2)\times O(2)$ symmetry of the system (translations and reflections 
in $x$, and in $\bphi$, where for the simple bifurcations only one is 
nontrivial, translation for pearling, and rotation for wrinkles) is 
inherited by the reduced equation which hence has the form 
$
\ddT A= A g(A,\overline A,(\al_j)_{j=1,\dots,N-1},(\beta_j)_{j=1,\dots,N-1}),
$
 where $g$ is an even function in all variables, i.e., in lowest order 
\hueq{ \label{normform}
\ddT A= A\left(a\beta_2 + b|A|^2\right).
}
At the ``double'' BPs $m\ne 0\ne n$ two types of branches bifurcate, and we 
give the ansatz below. 

\brem\label{rerem1}
Near a bifurcation point from a stable branch, 
the critical mode determines the evolution, while all others are 
exponentially damped. Hence, in a standard AE setting, the signs of 
$a$, $\beta_2$ and $b$ in \reff{normform} determine the stability on 
the bifurcating branch. Here, from \reff{normform} 
we get stability information for steady states under area constraint, 
with $\lam_2$ as extrinsic fixed parameter, but this can be extended 
to stability under area and volume constraints as follows. 
Close to the bifurcation point $\lam_2(m,n)$ the solution reads 
$X(\eps)=\CC_L-\eps A \psi N$. 
Then the solution is linearly stable iff $\pa_U \CG(X(\eps),\Lam)$ 
(which includes the $q_1$ area constraint) has only negative eigenvalues. 

For admissible (area and volume preserving) 
perturbations $\phi:\CC_L\ra \R$ of $X(\eps)$, we have the 
necessary and sufficient conditions 
\huga{ 
\int_{X(\eps)} H(\eps A \psi)\phi\dd S=0\quad\text{ and }\quad
\int_{X(\eps)} \eps A\psi g^{-1/2}(\eps A\psi) \phi \dd S=0.}
For the mean curvature we have 
$H(\eps A \psi)= \frac 1 2 - \eps A (\mu(m,n) +1) \psi +\CO(\eps^2)$, 
and hence admissible perturbations must be orthogonal to $\psi$. 
In order to construct an admissible perturbation $\tilde\phi$, we 
take a Fourier mode $\psi_{kl}$ as a leading order term. Then we find that if 
$\pa_u\CA(\eps)\psi_{kl}=\CO(\eps^2)$ and $\pa_u\CV(\eps)\psi_{kl}=\CO(\eps^2)$, 
we can construct a $\phi=\psi_{kl}+g(\eps)$, with $g=\CO(\eps)$. Hence, for a 
branch described by $X(\eps)=\CC_L-\eps A \psi_{mn} N+\CO(\eps^2)$, 
the only possible 
perturbations have a leading term with $(m,n)\neq (k,l)$. 
Thus, near bifurcation points associated with loss of stability 
of $\CC_L$ we have the following results: 
\bci 
\item Branches bifurcating from a simple BP (pearling or wrinkling) 
are stable wrt \reff{hflow}. 
\item Solutions on branches from a double BP 
(buckling and coiling) can only be destabilized by the other 
kernel vectors. We will see examples of this in \S\ref{numex}. 
\hfill\mbox{$\rfloor$}
\eci 
\erem

\subsubsection{Pearling instability}\label{pea-stab}
For the pearling instability (with for notational simplicity $m=1$) we have 
\hueqst{
\eps\, u= \eps A\,\psi_{10} + \eps^2 \left( \frac 1 2 A_{0}\psi_{00} +A_{2}\, \psi_{20}\right) +\cc +\CO(\eps^3),\quad 
\eps=\sqrt{\left|\frac{\lam_2(1,0)-\lam_2}{\beta_2}\right|}, 
}
with $A=A(T),\ T=\eps^2 t$. Plugging into \reff{bifprob} and 
sorting by powers of $\eps$, with the 
abbreviations $\lam_2=\lam_2(1,0)$ and $k=k_1$, we obtain the following: \\[-3mm]
 \bci
 \item[$\eps \psi_{10}:$] 
$\ds 0=-\frac{1}{2}\left(\left(k^{4}+\left(-4 c_{0}-2 \lambda_{2}\right) k^{2}+2 \lambda_{2}+1\right) \right)A.$
 \item[$\eps^2\psi_{00}:$] $\ds 0=\left(-\lambda_{2}-\frac{1}{2}\right) A_{0}+\frac{1}{2}\left(\left(k^{4}+\left(-8 c_{0}-2 \lambda_{2}-1\right) k^{2}+4 \lambda_{2}+5\right) {| A |}^{2}\right)-\beta_{2}-\alpha_{2}.$
 \item[$\eps^2\psi_{20}:$] $\ds 0=\left(-\frac{7 k^{4}}{4}{+}\frac{\left(-3-8 c_{0}+2 \lambda_{2}\right) k^{2}}{4}{+}\lambda_{2}{+}\frac{5}{4}\right) A^{2}+\left(-8 k^{4}{+}\frac{\left(32 c_{0}{+}16 \lambda_{2}\right) k^{2}}{4}{-}\lambda_{2}{-}\frac{1}{2}\right) A_{2}$. 
 
 \item[$\eps^3\psi_{10}:$] 
\hualst{\hs{-4mm}\ddT A= &\left(\frac{5 k^{6}}{2}{+}\frac{ \left(1{-}8 c_{0}{-}3 \lambda_{2}\right) k^{4}}{2}{+}\frac{\left(\frac{7}{2}{+}12 c_{0}{+}\lambda_{2}\right) k^{2}}{2}{-}\frac{27}{4}{-}3 \lambda_{2}\right) A {| A |}^{2}\\
&+\left(\frac{\left({-}4 c_{0}{-}1\right) k^{2}}{2}{+}2 \lambda_{2}{+}\frac{5}{2}\right) A A_{0}\\
 &{+}\left(4 k^{4}{+}4 \left({-}\frac{3}{8}{-}\frac{5 c_{0}}{2}{-}\frac{\lambda_{2}}{2}\right)  k^{2}{+}4 \left(\frac{5}{8}{+}\frac{\lambda_{2}}{2}\right) \right) \overline{A}A_{2}+\left(-k^{2}+1\right) \alpha_{2} A. }
 \eci
 From the area constraint we get the condition $0=A_{0} +k^2|A|^2$ for the 
$\psi_{00}$ mode, 
which gives 
 \hueqst{ \al_2=\frac{1} 2\left(k^{4}-8 k^{2} c_{0}+4 \lambda_{2}+5\right) {| A |}^{2}-\beta_{2}.}
We used directly the $\eps$ scale of $\lam_2$, hence we only get back the zero eigenvalue at $\eps\psi_1$. The above equations at $\eps^2 \psi_{00}$ and 
$\eps^2\psi_{20}$ give 
\hual{\label{a0a} 
A_0&=\frac{1}{2 \lambda_{2}+1}\left(\left(k^{4}+\left(-8 c_{0}-2 \lambda_{2}-1\right) k^{2}+4 \lambda_{2}+5\right) {| A |}^{2}+-2 \alpha_{2}-2 \beta_{2}\right), 
\text{ and} \\
A_2&=-\frac{A^{2} \left(7 k^{4}+\left(3+8 c_{0}-2 \lambda_{2}\right) k^{2}-4 \lambda_{2}-5\right)}{2 \left(16 k^{4}+\left(-16 c_{0}-8 \lambda_{2}\right) k^{2}+2 \lambda_{2}+1\right).}
}
The dependence on $\beta_2$ here comes from the dependence of 
$\lam_1$ on $\lam_2$. 
Combining everything leads at $\eps^3\psi_{10}$ 
to the solvability condition \reff{normform} with 
\hualst{a=&k^2-1, \\
b=& \frac{5}{\left(k^{2}{-}1\right) k^{2} \left(4 k^{4}{-}5 k^{2}{+}4 c_{0}{-}1\right)}  \Bigg(\frac{3}{20}k^{14}{+}\left(4 c_{0}{-}\frac{73}{20}\right) k^{12}{+}\left(\frac{8 c_{0}}{5}{-}\frac{13}{10}\right) k^{10}\\
&{+}\left(\frac{32}{5} c_{0}^{2}{-}\frac{52}{5} c_{0}{+}\frac{73}{20}\right) 
k^{8}{+}\left(\frac{16}{5} c_{0}^{2}{-}8 c_{0}{+}\frac{31}{5}\right) k^{6}{+}\left(\frac{4 c_{0}}{5}{-}\frac{43}{20}\right) k^{4}{+}\left(\frac{12 c_{0}}{5}{-}\frac{3}{2}\right) k^{2}\Bigg). }
Recall that here $k=k_m=2m\pi/L$, hence $a<0$ for $L>2m\pi$, and similarly 
$b<0$ for non small $L$. For instance, $a\approx -0.6$ and $b\approx -3$ 
in Fig.\ref{ampc00f}, where we compare the AE predictions for $m=1$ pearling 
with $c_0=0$ and $L=10$ with numerics. In any case for $ab<0$ we set 
$\beta_2=-1$  (subcritical case), and 
otherwise $\beta_2=1$ (supercritical case), 
giving the amplitude $A=\sqrt{\left|\frac a b\right|}$, with also 
$\al_2$ uniquely determined.  
The change of the volume can be approximated by inserting the amplitude Ansatz, 
giving  $0=A_0+A^2+\beta_2$ at $\eps^2\psi_{00}$, and comparison with 
\reff{a0a} shows that preserving the area and the volume is not possible 
along the pearling branch. 

\subsubsection{Wrinkling instability}
For the primary wrinkling instability at $\lam_2(0,2)=3/2$ we have the ansatz 
\hueqst{ \eps\, u(\bphi)= \eps Ae^{2\ri\bphi} + \eps^2\left(\frac 1 2A_0+A_4e^{4\ri\bphi}\right) +\cc +\CO(\eps^3),
}
again with $\eps=\sqrt{ \frac{|\lam_2-\lam_2^{(0,2)}|}{\beta_2}}$ and $A=A(T)$, 
$T=\eps^2 t$, and obtain the following. 
\bci
\item[$\eps \psi_{02}$:] $\ds 0=3 A  \left(\lambda_{2}-\frac{3}{2}\right)$. 
\item[$\eps^2 \psi_{00}$:] $\ds 0=\frac{1}{4}\left(-40 \lambda_{2}+82\right) {| A |}^{2}-\left(\lambda_{2} +\frac{1}{2}\right)A_{0}-\beta_{2}-\alpha_{2}$. 

\item[$\eps^2\psi_{04}$:] $\ds 0=\left(-9 \lambda_{2}+\frac{369}{4}\right) A^{2}+\left(15 \lambda_{2}-\frac{225}{2}\right) A_{4}$.
\item[$\eps^3\psi_{02}$:] $\ds \ddT A=-6 \left(\lambda_{2}-\frac{11}{4}\right) A A_{0}+\frac{1}{4}\left(-120 \lambda_{2}+1170\right) \overline{A} A_{4}-3 A \alpha_{2}-3 \left(-\lambda_{2}+\frac{221}{4}\right) {| A |}^{2} A$. 
\eci
From the area constraint at $\eps^2\psi_{00}$ we have 
$0= 8|A|^2+2A_0$, 
leading to $\al_2= \frac{3 \left(15-4 \lambda_{2}\right) {| A |}^{2}}{2}-\beta_{2}$. 
Thus $\lam^{(0,2)}_2=\frac 3 2$, $A_0=\frac {11} 4 |A|^2 -\frac{\beta_2} 2-\frac{\al_2} 2$, $A_4=\frac 7 8 A^2$ and finally $\al_2=\frac{27} 2 |A|^2-\beta_2$.
This ultimately leads to the solvability condition 
\hueq{\label{wria1} \ddT A = (3 \beta_2 -15.19 |A|^2) A,
}
which is notably independent of $c_0$ and $L$, and yields supercritical 
wrinkling bifurcations for $\beta_2=1$. 

\subsubsection{Coiling/buckling instability}
For the instability with critical wave vector 
$(m,n)=(1,1)$ and $(m,n)=(1,-1)$ we choose the ansatz 
\hueq{\label{mma1} \eps\, u(x,\bphi)= \eps\left( A\psi_{11}{+}B \psi_{1-1}\right) 
+ \eps^2\left(\frac 1 2 A_{00}{+}A_{22}\psi_{22}{+}B_{22}\psi_{2-2}{+}A_{20}\psi_{20}{+}A_{02} \psi_{02} \right) +\cc +\CO(\eps^3),}
with $\eps=\sqrt{ \frac{|\lam_2-\lam_2^{mm}|}{\beta_2}}$ and $T=\eps^2 t$. Using the amplitude formalism, this gives a cubic {\em system} of reduced 
equations for $A,B$, which by $O(2)\times O(2)$ symmetry must read 
\bsub\label{ampsys}
\hual{ 
\ddT A &= A (a\beta_2+b_1 |A|^2+b_{2} |B|^2),\\
\ddT B &=  B (a\beta_2+b_{2} |A|^2+b_1 |B|^2).
}
\esub 
 From that we can directly recover the two types of solutions: 
Coiling with $\ds A=\sqrt {\left|\frac{a}{b_1}\right|}$ and $B=0$ or 
the other way round; 
buckling with $\ds A=B =\pm\sqrt{\left|\frac {a}{b_1+b_2}\right|}$. 
Here again $\beta_2=\pm 1$, depending on the signs of $a$, $b_1$, and $b_1+b_2$. 

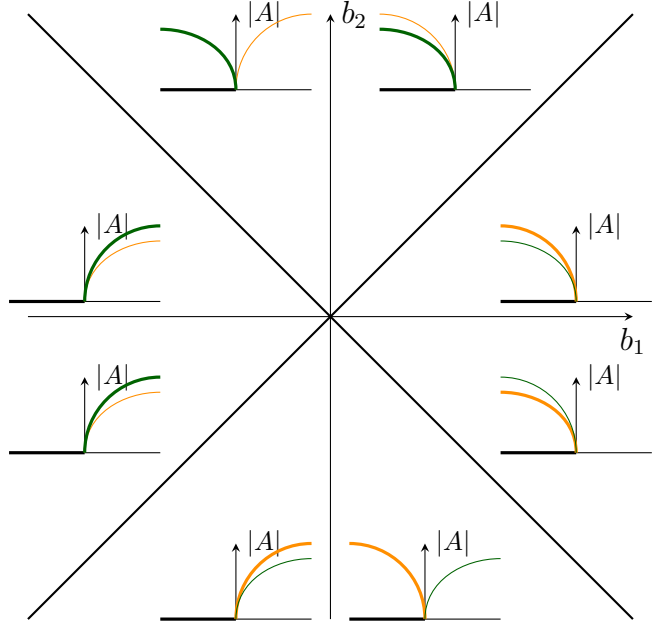
\begin{figure}
\bce
\begin{tikzpicture}
\draw[ -stealth] (-4,0) -- (4,0) node[below ] {$b_1$};
\draw[ -stealth] (0,-4) -- (0,4) node[right ] {$b_2$};
\draw[ thick] (0,0) -- (-4,-4)  ;\draw[ thick] (0,0) -- (4,-4);\draw[ thick] (0,0) -- (4,4);
\draw[ thick] (0,0) -- (-4,4);
\definecolor{o1}{rgb}{1,0.56,0};
\definecolor{g1}{rgb}{0,0.39,0};
\coordinate (ref) at (-1.25,-4) ; 
\draw[very thick] ($(ref)+(-1,0)$)-- ($(ref)$);
\draw($(ref)$)-- ($(ref)+(1,0)$);
\draw[ -stealth] ($(ref)$)-- ($(ref)+(0,1)$) node[right, font=\footnotesize ] {$|A|$};
\draw[ very thick, color=o1] (ref) to [out=90,in=180 ] ($(ref)+(1,1)$);
\draw[ color=g1] (ref) to [out=90,in=180 ] ($(ref)+(1,.8)$);
 
 \coordinate (ref) at (1.25,-4) ; 
\draw[very thick] ($(ref)+(-1,0)$)-- ($(ref)$);
\draw($(ref)$)-- ($(ref)+(1,0)$);
\draw[ -stealth] ($(ref)$)-- ($(ref)+(0,1)$) node[right, font=\footnotesize ] {$|A|$};
\draw[ very thick, color=o1] (ref) to [out=90,in=0 ] ($(ref)+(-1,1)$);
\draw[ color=g1] (ref) to [out=90,in=180 ] ($(ref)+(1,.8)$);

 \coordinate (ref) at (3.25,-1.8) ; 
\draw[very thick] ($(ref)+(-1,0)$)-- ($(ref)$);
\draw($(ref)$)-- ($(ref)+(1,0)$);
\draw[ -stealth] ($(ref)$)-- ($(ref)+(0,1)$) node[right, font=\footnotesize ] {$|A|$};
\draw[ very thick, color=o1] (ref) to [out=90,in=0 ] ($(ref)+(-1,.8)$);
\draw[ color=g1] (ref) to [out=90,in=0 ] ($(ref)+(-1,1)$);

 \coordinate (ref) at (3.25,.2) ; 
\draw[very thick] ($(ref)+(-1,0)$)-- ($(ref)$);
\draw($(ref)$)-- ($(ref)+(1,0)$);
\draw[ -stealth] ($(ref)$)-- ($(ref)+(0,1)$) node[right, font=\footnotesize ] {$|A|$};
\draw[ very thick, color=o1] (ref) to [out=90,in=0 ] ($(ref)+(-1,1)$);
\draw[ color=g1] (ref) to [out=90,in=0 ] ($(ref)+(-1,.8)$);

 \coordinate (ref) at (1.65,3) ; 
\draw[very thick] ($(ref)+(-1,0)$)-- ($(ref)$);
\draw($(ref)$)-- ($(ref)+(1,0)$);
\draw[ -stealth] ($(ref)$)-- ($(ref)+(0,1)$) node[right, font=\footnotesize ] {$|A|$};
\draw[  color=o1] (ref) to [out=90,in=0 ] ($(ref)+(-1,1)$);
\draw[very thick, color=g1] (ref) to [out=90,in=0 ] ($(ref)+(-1,.8)$);

 \coordinate (ref) at (-1.25,3) ; 
\draw[very thick] ($(ref)+(-1,0)$)-- ($(ref)$);
\draw($(ref)$)-- ($(ref)+(1,0)$);
\draw[ -stealth] ($(ref)$)-- ($(ref)+(0,1)$) node[right, font=\footnotesize ] {$|A|$};
\draw[  color=o1] (ref) to [out=90,in=180 ] ($(ref)+(1,1)$);
\draw[ very thick,color=g1] (ref) to [out=90,in=0 ] ($(ref)+(-1,.8)$);

 \coordinate (ref) at (-3.25,.2) ; 
\draw[very thick] ($(ref)+(-1,0)$)-- ($(ref)$);
\draw($(ref)$)-- ($(ref)+(1,0)$);
\draw[ -stealth] ($(ref)$)-- ($(ref)+(0,1)$) node[right, font=\footnotesize ] {$|A|$};
\draw[color=o1] (ref) to [out=90,in=180 ] ($(ref)+(1,.8)$);
\draw[ very thick,color=g1] (ref) to [out=90,in=180 ] ($(ref)+(1,1)$);

 \coordinate (ref) at (-3.25,-1.8) ; 
\draw[very thick] ($(ref)+(-1,0)$)-- ($(ref)$);
\draw($(ref)$)-- ($(ref)+(1,0)$);
\draw[ -stealth] ($(ref)$)-- ($(ref)+(0,1)$) node[right, font=\footnotesize ] {$|A|$};
\draw[color=o1] (ref) to [out=90,in=180 ] ($(ref)+(1,.8)$);
\draw[ very thick,color=g1] (ref) to [out=90,in=180 ] ($(ref)+(1,1)$);
\end{tikzpicture}

\vs{4mm}
\caption{ {\small Schematic linearized stability diagram of the amplitude equation depending on $b_1$ and $b_2$. Orange stands for the coiling and green for 
the buckling branch. Thick lines indicate stable solutions. }}
\label{stabfig}
\ece
\end{figure}

We have $a=k_m^2$, but the computation of the coefficients $b_1$ and $b_2$ 
naturally is more 
cumbersome than for the scalar case (pearling and wrinkling), and therefore 
relegated to Appendix \ref{aeapp}. Here we instead briefly comment 
on the bifurcation directions, and stability of the bifurcating branches under  
volume and area preserving flows, see also Fig.\ref{stabfig}. 
The AEs are the same as in the classical $D_4$ symmetry problem 
\cite[\S4.3.1]{Hoyle} and lead to the same direction of bifurcations. 
However the stability information gained is completely different, due to the constraints, cf.~Remark \ref{rerem1}. Namely, 
the Jacobian of the right hand side of \reff{ampsys} at a steady state $(A,B)$ is 
 \hualst{ \bpm 
a+3b_1 A^2+b_2B^2 & 2b_2 AB\\ 2b_2 AB &a+b_2A^2+3b_1 B^2\epm.
}
For $a>0$, on the coiling solution branch the eigenpairs are 
\huga{\text{$\sig_1= (a\beta_2+3\sign(b_1)|a\beta_2|)$, $v_1=(1,0)^T$,  and 
$\sig_2=( a\beta_2+\frac{b_2}{|b_1|} |a\beta_2|)$, $v_2=(0,1)^T$.}}
The first direction $v_1$ is excluded by the constraints as it 
points along the coiling branch, and hence the coiling is stable for 
$\frac {b_2}{|b_1|} < -\sign(\beta_2)$. 
For the buckling the Jacobian has the eigenpairs  
  \hueq{
\text{$\sig_1=2 \frac{b_2+b_1}{|b_1+b_2|}, v_1=(1,1)^T$, and 
$\sig_2=\frac{b_2-b_1}{|b_1+b_2|} , v_2=(1,-1)^T$. } 
}
From the constraints again $v_1$ is excluded and hence 
the relevant eigenvalue is $\sig_2$, which is negative for $b_2<b_1$. 
In \S\ref{numex} we give two comparisons of the predictions from 
\reff{ampsys} with numerics, namely $c_0=0$ and $L=10$ yielding 
$(a,b_1,b_2)\approx (0.39, -0.67, -0.82)$ in Fig.\ref{c00f}, 
and $c_0=0.48$ and $L=10$ yielding 
$(a,b_1,b_2)=(0.39,-0.083, 0.56)$ in Fig.\ref{c005f2}. 
In both cases we get the correct branching and stability predictions, 
and also very good quantitative agreement.  

\section{Numerical examples} \label{numex}
We now use numerical continuation and bifurcation methods 
to corroborate the local results from 
\S\ref{aesec}, and to extend them to a more global picture. 
This includes some secondary bifurcations with interesting stable 
shapes away from the straight cylinder. The numerical setup follows 
the stability analysis in \autoref{stabrange} and \S\ref{redeq}. 
Namely, we fix a period $L$ (e.g., $L=10$, $L=15$ or $L=30$, see 
also Remark \ref{carem}),  and the area $\CA=\CA(X)$, and 
initially use $\lam_2$ as the primary bifurcation parameter 
and $\lam_1$ and $\CV(X)$ as secondary (dependent) parameters. 
The linearized stability of the solutions, however, is then computed for 
fixed $\CA(X)$ and $\CV(X)$ while $\lam_1$ and $\lam_2$ are free parameters. 
Due to the radius-length scaling invariance of the Helfrich energy \reff{helfen} and the area and volume constraint in the dynamical 
problem, it is useful to introduce the reduced volume  
\hueq{\label{redvol} v=\CV(X)/\CV_0,} 
 where $\CV_0$ is the volume of the straight cylinder with the same $\CA(X)$, 
compare \cite{JSL93}, and we normalize the energy 
\hueq{\label{enorm}
\tilde E(X)=\frac 1 \CA \int_X (H-c_0)^2\dd S.} 
Also for our numerics we focus on the three exemplary cases 
\hugast{
\text{(a) $c_0=0$,\quad (b) $c_0=-1$, \quad and (c) $c_0=0.48$}
} 
from Fig.\ref{neigs}, and mostly 
we use $L=10$. For (a) and (b) 
this is no significant restriction as in these cases the primary 
bifurcations (to pearling, buckling and coiling in (a), and to pearling and 
wrinkling in (b)) are of long wave type for any (sufficiently large) $L$, 
i.e., a different $L$ gives the same types of primary bifurcations with 
different axial (maximal) periods $L$ (and the wrinkling is independent of $L$). 
However, for (c) and sufficiently large $L$ (preferably an integer 
multiple of the critical period $\ell\approx 7.5$) we have a 
finite wave number pearling instability. Thus, besides the case $L{=}10$ 
(for comparison with (a) and (b)), for (c) we shall consider different 
multiples of 7.5 for $L$, namely $L{=}15, L{=}30$, and $L{=}45$, allowing 
2, 4, and 6, pearls on the primary branch. The influence of $L$ 
will then be the behavior of this branch away from bifurcation, 
which loses stability in a bifurcation to a single pearl branch, 
and this BP moves closer to the primary BP with increasing $L$. 
This is further discussed in \S\ref{dsec}. 
\huc{Additionally, in \S\ref{cpsec} we give an outlook on case 
(d) $c_0=0.75$, 
from Fig.\ref{neigs}, where $\CC_L$ is unstable for all $\Lam$.}

In the numerics we (as always) have to find a compromise between 
speed and accuracy. For the default $L=10$ we start with an initial 
discretization of the straight cylinder $\CC_L$ by $n_p\approx 3000$ 
mesh points; the continuation of the nontrivial branches then often 
requires mesh refinement (and coarsening), typically leading 
to meshes of 4000 to 5000 points. For the longer cylinders for $c_0=0.48$ 
we increase this to up to 12000 mesh points. A simple measure of numerical 
resolution is the comparison of the $\lam_2$ values of the 
numerical BPs from $\CC_L$ with their analytical values from \S\ref{lineig}, 
and here we note that in all cases we have agreement of at least 4 digits. 
Moreover, we also compare 
some numerical nontrivial branches with their descriptions by the AEs, 
and find excellent agreement. 

\brem\label{vrem}{\rm 
Some cylindrical solutions have similarities with 
axisymmetric closed vesicle solutions of \reff{helf0}. 
We recall from \cite{JSL93} or \cite{geompap}, that from the sphere 
the first bifurcation is to  prolates (rods) and oblates 
(wafers). 
For mildly negative or positive $c_0$  the prolates are stable near 
bifurcation and the oblates unstable, and continuing the prolates 
they become dumbbells (two balls connected by a  neck), 
from which pea shaped vesicles bifurcate which become stable as one of the 
ball shrinks. For sufficiently 
negative $c_0$, the oblates are stable near bifurcation, and continuing 
the branch they become discocytes of biconcave red blood cell shape. 
Here we see some reminiscent behavior: For sufficiently negative $c_0$ 
we find secondary bifurcations of stable wrinkles with 
``embedded biconcave shapes'', see, e.g., Fig.\ref{cmf1}E. 
For sufficiently positive $c_0$, the primary pearling bifurcation 
leads to ``strings of dumbbells'' rather than strings of (round) pearls, 
see, e.g., samples C in Figures \ref{c05f1} and \ref{c05f1l}, or the 
samples A$_1$--A$_3$ in Fig.\ref{c1f1}. 
}\eex\erem

We complement our results by some (area and volume conserving) 
numerical flows close to steady states, and for these 
perturb steady states $X_0$ in two different ways: one is 
\huga{\label{flp1}
\tilde{X}_0=X_0+\del u_0N_0,\quad u_0=\frac 1{1+\xi x^2}\cos(\bphi), 
}
where $x$ is the coordinate along 
the cylinder axis, and $\bphi$ the polar angle orthogonal to it, 
and $\del\in\R$ and $\xi>0$ are parameters, with typically $\xi=1$ and 
$\del=-0.125$ or $\del=-0.25$ (recall that $N$ is the inner normal).  
The other way is  
\huga{\label{flp2}
\tilde{X}_0=X_0+\del u_0N_0,\quad 
\quad \text{with $u_0=\psi_{\text{crit}}$}, 
}
where $\psi_{\text{crit}}$ is 
(the $u$ component of) the eigenvector to the largest 
eigenvalue of the Jacobian 
at $X_0$, and typically $\del=-0.5$. 
Due to the negative $\del$ combined with the 
inner normal $N$, perturbations of type \reff{flp1} 
in general induce a small increase of the initial volume 
$\CV_0=\CV+\CV_\del$ and area $\CA_0=\CA+\CA_\del$, which are conserved 
under the flow. 
On the other hand, perturbations of type \reff{flp2} are 
area and volume preserving to linear order.

Our numerical flows 
$t\mapsto (X(t),\Lam(t))=\Psinum(t;h)X_0$ are
{\em ad hoc} approximations  with a time stepsize $h$ 
and tolerance {\tt tol} (see Appendix \ref{num}), and  
we summarize the results from our flow experiments as follows:
\bcen 
\item Linearly stable $X_0$ (at $(\CV,\CA)$) are also stable under the flow 
at the perturbed $(\CV_0,\CA_0)$, and the solution flows back to the respective 
type (pearling, coiling, ...) at $(\CV_0,\CA_0)$, with a small 
change in $(\lam_1,\lam_2)$. See, e.g., Fig.\ref{c0f2}(a,b). 
\item If $X_0$ at $(\CV,\CA)$ is linearly unstable, then we generally 
expect $X(t)$ to flow to the nearest stable steady state 
(of in general a different class), see Fig.\ref{c0f2}(c$_2$). 
However, if we have a significant 
perturbation of $(\CA_0,\CV_0)$, then $X(t)$ may still flow to a similar 
steady state close to $X_0$ by adapting 
$(\lam_1,\lam_2)$, see, e.g., Fig.\ref{c0f2}($c_1$). 
This underlines the fact 
that $\lam_{1,2}$ are dynamical variables. 
\item Naturally, we want $E$ to decrease monotonously in $t$, cf.\reff{helf1b}, 
which if violated gives an easy indicator of poor numerics. 
The monotonous decrease of $E$ often requires rather small stepsizes 
and tolerances, and a result is that in many examples 
a single numerical flow to steady state is as expensive as 
the computation of the full BD of steady states, 
which we take as a strong argument 
for our steady state continuation and bifurcation approach. 
\item \huc{For $c_0=0.75$, where $\CC_L$ is always unstable and hence 
no stable primary bifurcations occur, we can numerically flow 
to stable finite wavelength pearling, which we subsequently continue in 
parameters, see Fig.\ref{c1f1}. However, this is only an outlook, 
and we leave a detailed and systematic study of this for future work.}
\ecen

\subsection{$c_0=0$}
\begin{figure}[ht]
  \begin{center}
\btab{l}{
\btab{ll}{{\small (a$_1$)}&{\sm (a$_2$)}\\  
\hs{-0mm}\ig[height=0.33\tew]{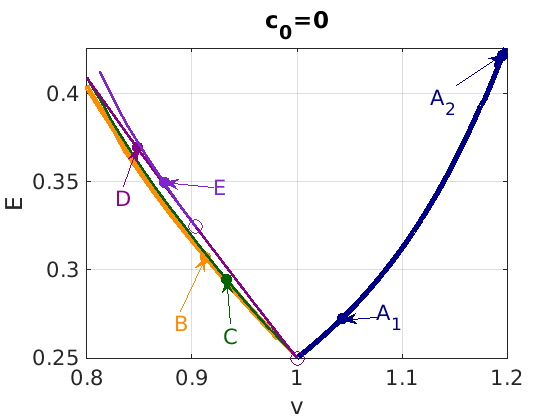} &
\hs{-4mm}\ig[height=0.33\tew]{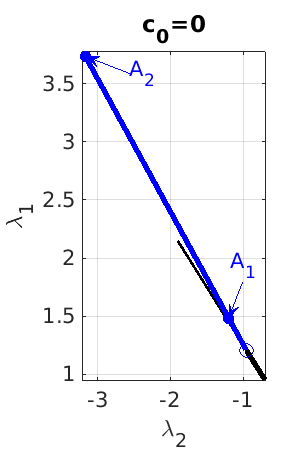} 
\hs{-2mm}\ig[height=0.33\tew]{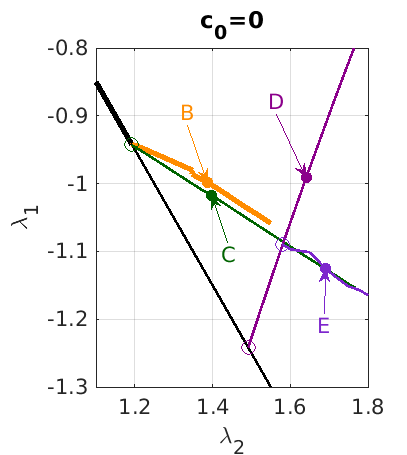}}\\
{\small (b)}  \\[-0mm]
\ig[width=0.16\tew]{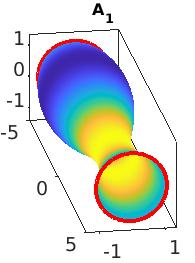}\hs{-0mm}\ig[width=0.16\tew]{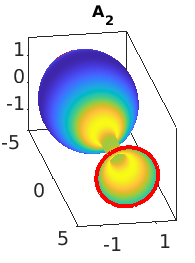}
\ig[width=0.16\tew]{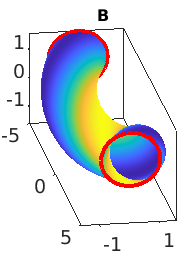}\hs{-0mm}\ig[width=0.16\tew]{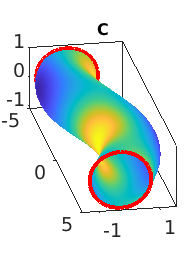}
\hs{-0mm}\ig[width=0.16\tew]{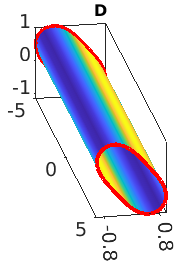}\hs{-0mm}\ig[width=0.16\tew]{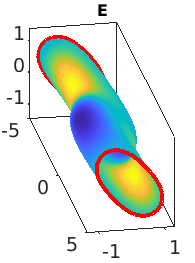}\\
} 
\vs{-0mm}
  \end{center}
  \caption{{\small Basic BD for $c_0=0$, $L=10$. (a$_1$) Energy $E$ over reduced volume 
$v$. The family of cylinders $\Co_{10}$ parameterized by $\lam_1$  
 corresponds to the point where all primary branches bifurcate; 
A=pearls, B=coils, C=buckling, D=wrinkling; E is a secondary bifurcation 
from D, yielding a mix of buckling and wrinkling. 
(a$_2$) Relations between $\lam_1$ and $\lam_2$. The cylinders are the 
black branch, stable for $\lam_2\in(-0.95,1.19)$. Bifurcation of A at
      $\lam_2\approx -0.95$, bifurcation of B,C at $\lam_2\approx 1.19$, 
and bifurcation of D at $\lam_2=3/2$. Samples in (b); here the colors 
indicate $u$ from the current continuation step, with blue for $u<0$ 
(yellow for $u>0$) 
indicating that $X$ moves outward (inward) along the branch. 
The periodic boundaries are shown in red for better visibility. 
\label{c00f}}}
\end{figure}

For $L=10$ and $c_0=0$, \S\ref{stabrange} gives that the cylinder 
destabilizes at $\lam_2\approx -0.95$ to pearling $(k_1,0)$, 
and at $\lam_2\approx 1.19$ to coiling/buckling $(k_1,1)$. 
Figure \ref{c00f} shows the BDs and samples starting in 
(a$_1$) with the energy $E$ over the reduced volume $v$.  

Here all branches 
bifurcate from $v=1$, but (a$_2$) shows that they are well separated in 
the $(\lam_1,\lam_2)$ plane. 
As stated in \autoref{pea-stab}, the blue pearling branch starts stable 
close to the bifurcation point (sample A$_1$), and stays stable up 
to rather strong deformations of the original cylinder, sample A$_2$. 
The area constraint prevents the branch from growing to arbitrary volumes, 
but due to neck development the continuation becomes more difficult for 
larger $v$. 
From the mode $(m,n)=(1,1)$ at $\lam_2\approx 1.19$ we find 
buckling and coiling. The coefficients in \reff{ampsys} 
are $(a,b_1,b_2)=(0.39,-0.4,-0.67)$, 
 and as predicted (see Remark \ref{arem1}) 
the coiling (sample B) is stable close to bifurcation, 
while the buckling (sample C) is unstable. This 
behavior continues to small $v$, with no secondary bifurcations. 

Additionally we show one wrinkling 
branch (sample D) as a  later bifurcation from the cylinder, and 
a secondary bifurcation to a wrinkling--buckling mix (sample E). 
However, these solutions are all linearly unstable in the $\lam_{1,2}$ range 
shown, at the given $\CV$ and $\CA$.

\begin{figure}[ht]
\bce
\btab{lll}{{\small (a)}\\
\ig[width=0.29\tew]{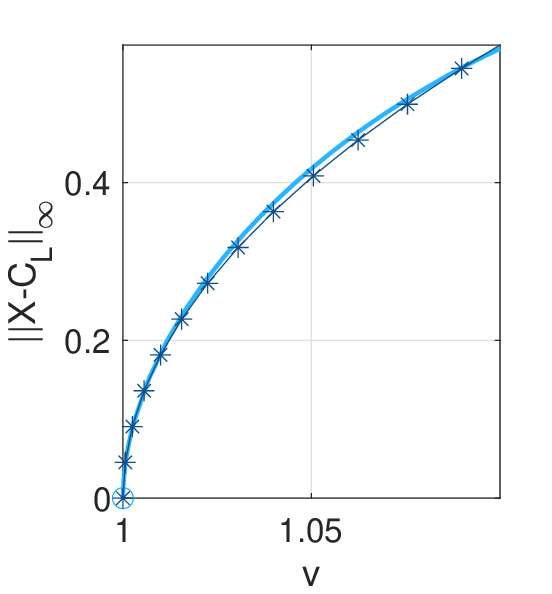} &
\ig[width=0.29\tew]{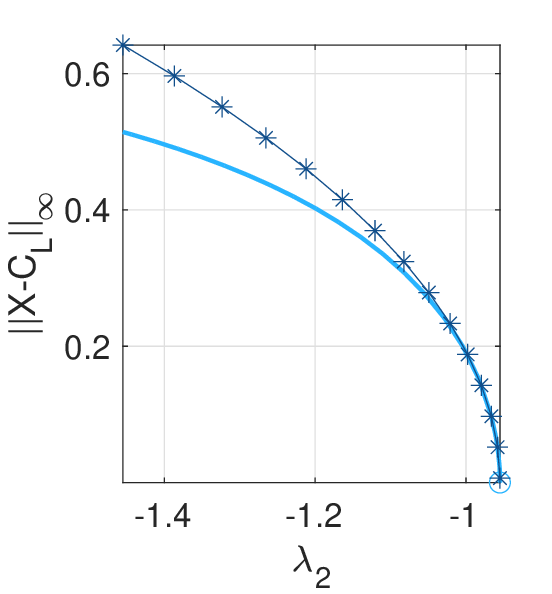} &
\ig[width=0.29\tew]{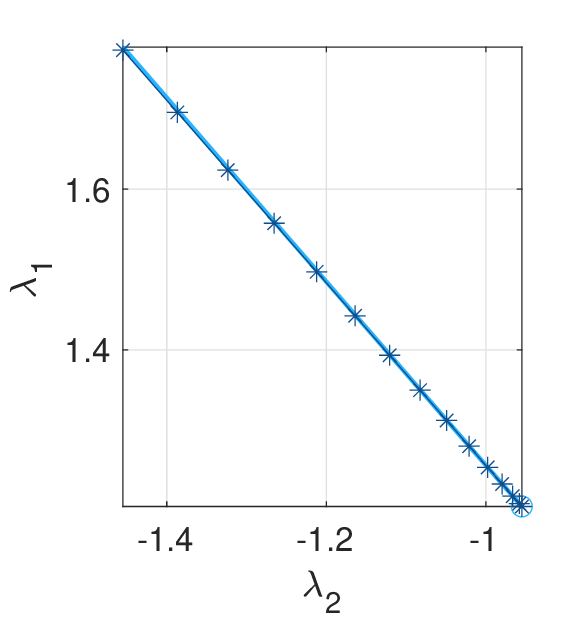} \\[-4mm]
{\small (b)}  \\
\ig[width=0.29\tew]{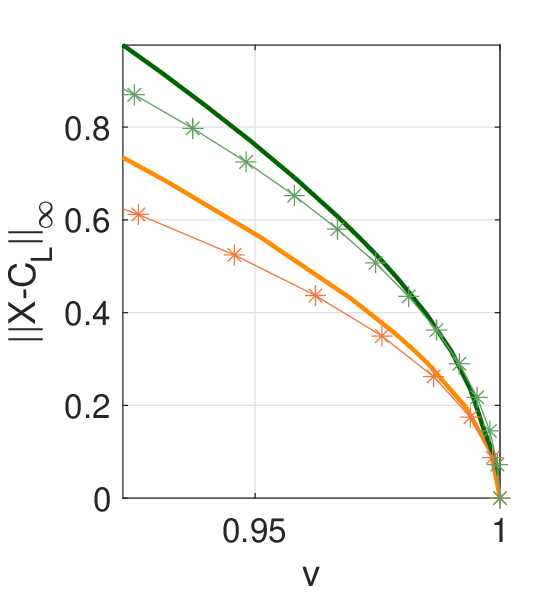} &
\ig[width=0.29\tew]{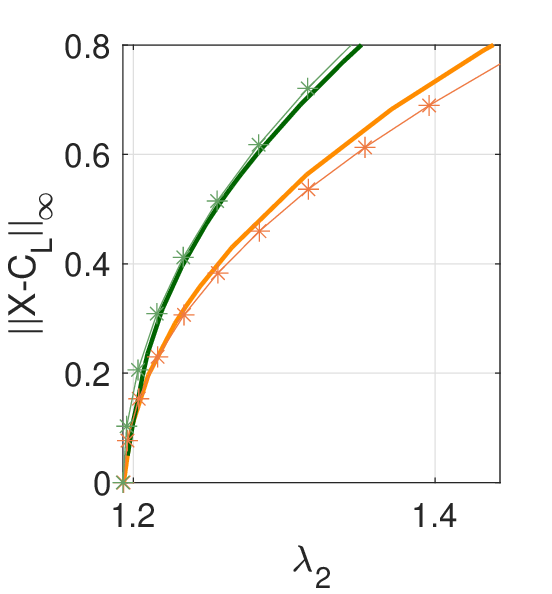} &
\ig[width=0.29\tew]{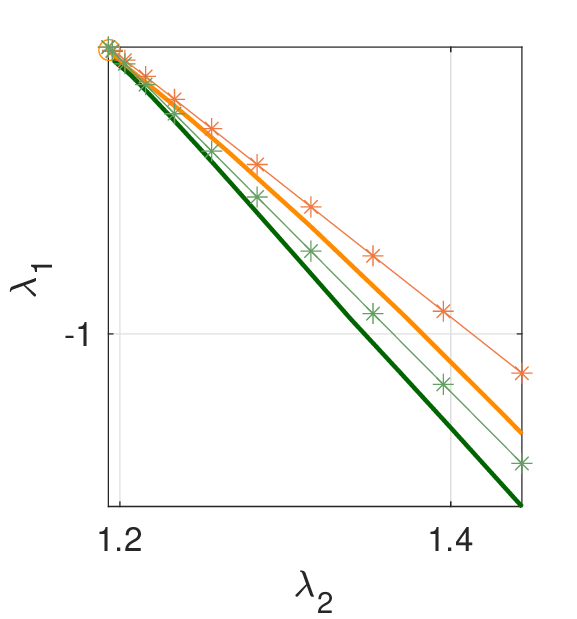} \\
}
\ece
\vs{-0mm}
\caption{{\small $c_0=0$, comparison of numerics and AE predictions 
(lines with *), (a) pearling, (b) coils (orange) and buckling (green). 
Left and center: Amplitude over reduced volume $v$, and over $\lam_2$. 
Right: $\lam_1$ over $\lam_2$.}}\label{ampc00f}
	\end{figure}

\brem\label{arem1} 
The amplitude formalism in \autoref{redeq} is only expected 
to be valid in a small neighborhood of the bifurcation points.  
In Fig.~\ref{ampc00f} we show comparisons between the numerical 
computations and the AE predictions. 
Here the amplitude formalism 
for the pearling
\hugast{ 
\text{$u=\eps A\psi_{10}+\CO(\eps^2)$,\quad $\lam_2=\lam_2(1,0)-\eps^2$, 
\quad $\lam_1=\lam_1(1,0)+\al_2\eps^2+\CO(\eps^4)$,}
}  
yields $\ds \ddT A =A(0.61-2.96A^2)$ and $\al_2= 0.66 A^2-1$, 
and that for the wrinkling $u=\eps A\psi_{01}+\CO(\eps^2)$, 
$\lam_2=\lam_2(0,1)+\eps^2$, $\lam_1=\lam_1(0,1){+}\al_2\eps^2{+}\CO(\eps^4)$
yields the universal AE \reff{wria1}, i.e., 
$\ds \ddT A=A(3{-}15.19A^2)$ and $\al_2=\frac {27} 2 A^2{-}1$. 
For both, pearling and wrinkling,  we find an excellent match with the numerics. 
\eex\erem

\begin{figure}[H]
  \begin{center}
\btab{l}{
\btab{ll}{
{\small (a),  from Fig.\ref{c00f}A$_2$}\\[-0mm]
 \hs{-4mm}\ig[height=0.24\tew]{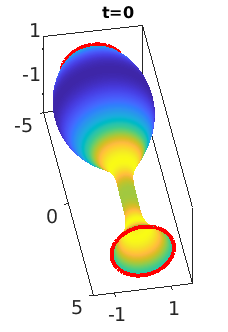} 
\hs{1mm}\ig[height=0.24\tew]{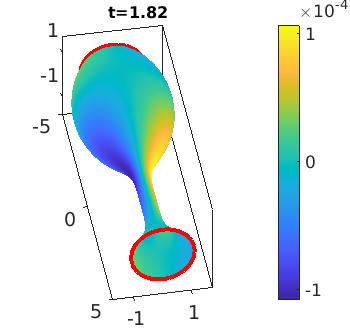} 
&\hs{-2mm}\ig[height=0.24\tew]{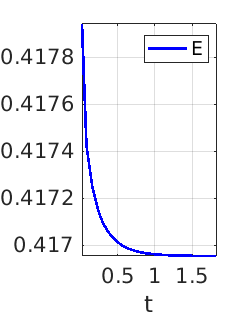}
\hs{-2mm}\ig[height=0.24\tew]{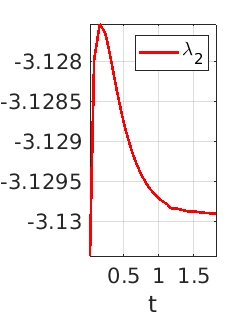}
\hs{-0mm}\ig[height=0.24\tew]{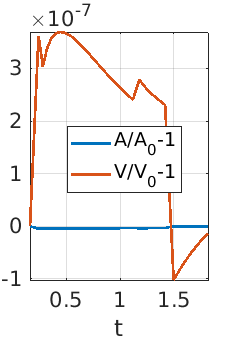}\\[-0mm]
{\small (b),  from Fig.\ref{c00f}B}\\[-0mm]
\hs{-4mm}\ig[height=0.24\tew]{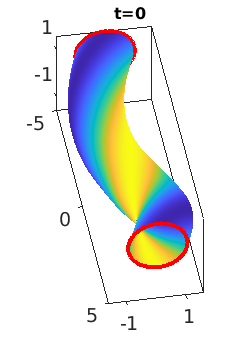} 
\hs{-1mm}\ig[height=0.24\tew]{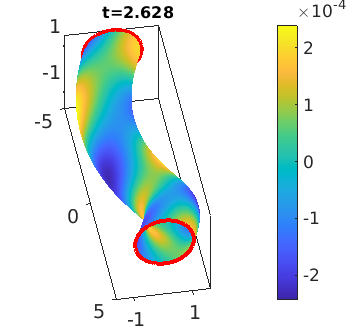} 
&\hs{-2mm}\ig[height=0.24\tew]{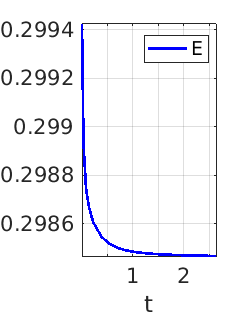}
\hs{-2mm}\ig[height=0.24\tew]{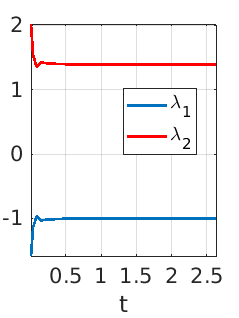}
\hs{-0mm}\ig[height=0.24\tew]{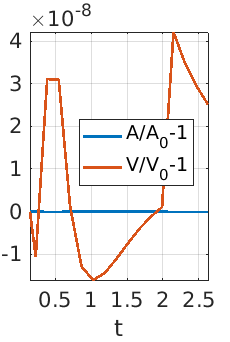}
\\[-0mm]
{\small (c$_1$), from Fig.\ref{c00f}C,  perturbation by \reff{flp1}}\\
\hs{-4mm}\ig[height=0.24\tew]{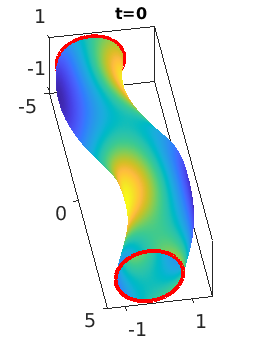} 
\hs{-1mm}\ig[height=0.24\tew]{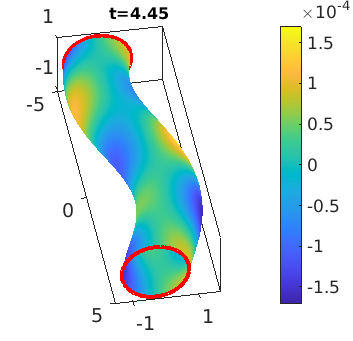} 
&\hs{-2mm}\ig[height=0.24\tew]{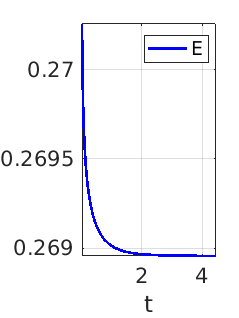}
\hs{-2mm}\ig[height=0.24\tew]{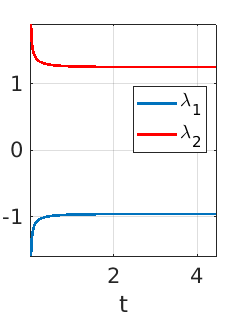}
\hs{-0mm}\ig[height=0.24\tew]{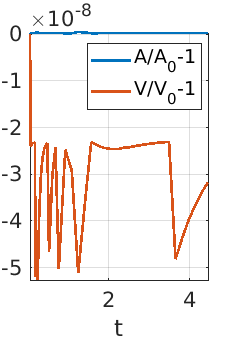}
\\[-0mm]
{\small (c$_2$), from Fig.\ref{c00f}C, perturbation by \reff{flp2}}\\
\hs{-4mm}\ig[height=0.24\tew]{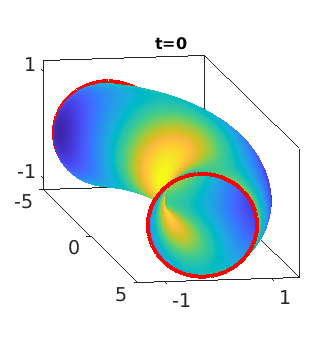} 
\hs{-1mm}\ig[height=0.24\tew]{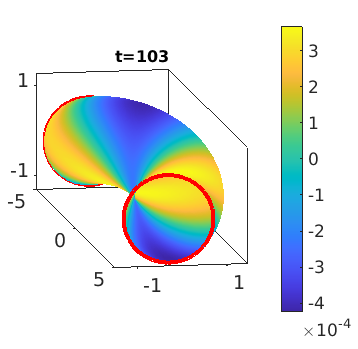} 
&\hs{-2mm}\ig[height=0.24\tew]{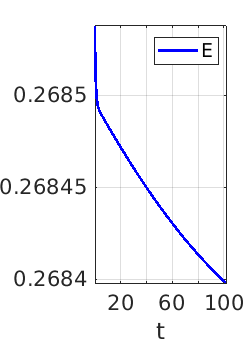}
\hs{-0mm}\ig[height=0.24\tew]{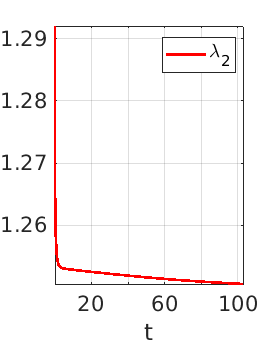}
\hs{-0mm}\ig[height=0.24\tew]{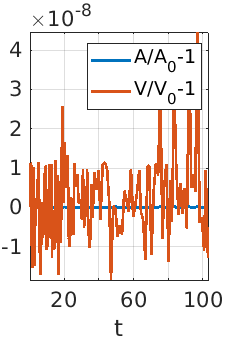}
\\[-0mm]
 }
}
\vs{-0mm}
  \end{center}
  \caption{{\small Numerical flows from perturbations of states from 
Fig.\ref{c00f}; initial and 'final' state, and evolution of $E$, $\lam_{1,2}$ 
(only $\lam_2$ in (a) and (c$_2$)), and deviations from initial area and volume; 
very slow evolution in (c$_2$). \label{c0f2}}}
\end{figure}

In Fig.\ref{c0f2} we show numerical flows from the steady states A--C 
from Fig.\ref{c00f}, perturbed according to \reff{flp1}, with parameters 
given in Table \ref{ptab1}. In all cases we integrate until 
$\|\pa_t u\|_\infty<0.001$. The plots on the left show the initially perturbed 
$X_0$, and $X$ at the end of the DNS colored by the last update $\del u$, 
such that $u_t\approx \del u/h$, where $h$ is the current stepsize. 
The plots on the right show the evolution of $E(t)$ and $\lam_{1,2}$, 
and the relative errors in $\CA$ and $\CV$ (with the error in $\CV$ 
always two orders of magnitude larger than that in $\CA$). 

The first two flows start near stable steady 
states $(X_0,\lam_{1,0},\lam_{2,0})$, and as expected flow back to 
near $X_0$ rather quickly, with partly significant 
changes of $(\lam_1,\lam_2)$. 
However, the third flow (c$_1$) is somewhat surprising. 
State C from Fig.\ref{c00f} 
is linearly unstable, but the perturbed $X_0$ flows to another 
buckling by adjusting $\lam_{1,2}$ at the increased area and 
volume, where this buckling is stable. On the other hand, in (c$_2$) 
for the 
(linearly) area and volume preserving initial perturbation into the 
most unstable direction (which point towards coiling), 
the behavior is different: The solution now converges to coiling, 
with a rather small flow throughout; 
$\lam_1$ hardly changes at all, $\lam_2$ only very slightly, 
and the convergence is very slow (since the leading eigenvalue 
of the target coiling is negative but very close to zero). For all flows, we can use 
the final state as an initial guess for steady state 
continuation (e.g., in $\lam_2$), 
and the obtained numerical steady states then are (volume preserving) 
linearly stable.

\taskip
\begin{table}[H]\caption{{\sm Spectral data ($n_{\text{unst}}$=number of ($V$-preserving) 
unstable directions, $\mu_{\text{crit}}$=least stable/most unstable eigenvalue), perturbation parameter 
$\del$ ($\xi{=}1$ throughout), and results for DNS in Fig.\ref{c0f2}; 
$A_0{=}20\pi, V_0{=}10\pi$. Initial perturbations for A--C of type 
(\ref{flp1}), and for C$^*$ of type \reff{flp2}. $(\lam_1,\lam_2)_0$ 
are the values at the (perturbed) steady state. 
 \label{ptab1}}}
\centering 
{\small
\begin{tabular}{l|llll}
IC type&A$_2$ (pearling)&B (coiling)&C (buckling)&C (buckling)$^*$\\
\hline
$n_{\text{unst}}, \mu_{\text{crit}}$&0, -2.22&0, -0.018&
1, 0.012&1, 0.012\\
$\del$&-0.125&-0.25&-0.25&0.5\\
$(A-A_0)/A_0$&0.0012&0.0003&0.0125&0.0024\\
$(V-V_0)/V_0$&0.0011&0.0004&0.02&0.004\\
$(\lam_1,\lam_2)_0$&(3.734, -3.165)&(-0.999, 1.387)&
(-0.972, 1.29)&(-0.972, 1.29)\\
$(\lam_1,\lam_2)_{\text{end}}$&(3.701, -3.133)&(-0.999, 1.389)&
(-1.003, 1.367)&(-0.965, 1.25)
\end{tabular}
}
\end{table}
\teskip

\subsection{$c_0=-1$}
\begin{figure}[ht]
  \begin{center}
\btab{l}{
\btab{ll}{
\hs{-4mm}\ig[height=0.33\tew]{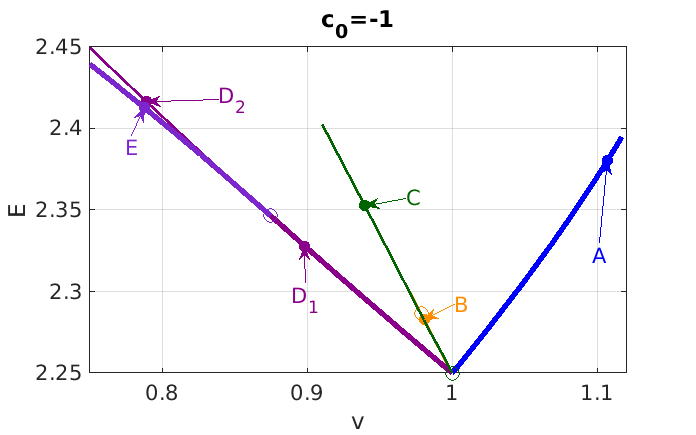} 
& \hs{-4mm}\ig[height=0.33\tew]{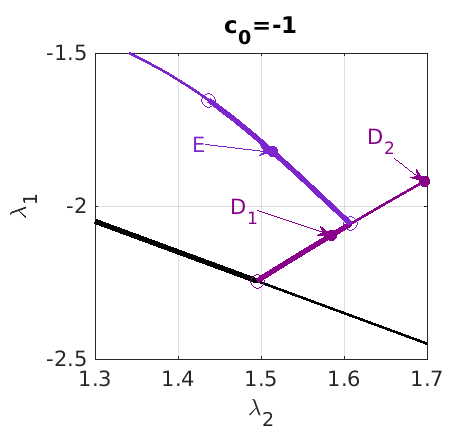} 
}\\
\ig[width=0.16\tew]{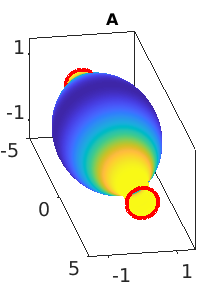}\hs{-0mm}\ig[width=0.16\tew]{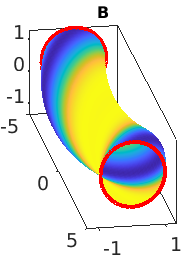}
\ig[width=0.16\tew]{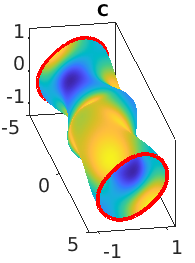}\hs{-0mm}\ig[width=0.16\tew]{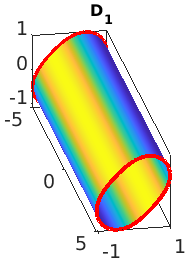}
\hs{-0mm}\ig[width=0.16\tew]{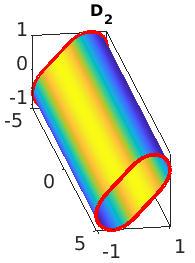}\hs{-0mm}\ig[width=0.16\tew]{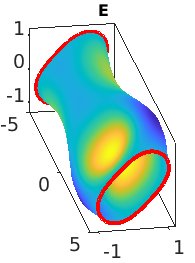}
} 
\vs{-0mm}
  \end{center}
  \caption{{\small Basic BD for $c_0=-1$, $L=10$, and samples.  
\label{cmf1}}}
\end{figure}

For (strongly) negative $c_0$, the wrinkling branches and their secondary  
bifurcations
become relevant. Fig.\ref{cmf1} shows basic results 
for $c_0{=}-1$ and $L{=}10$. 
As before, at negative $\lam_2\approx -2.254$ 
and $\lam_1\approx 1.5$ (cf.Fig.\ref{neigs}), $\Cr$ loses stability 
to pearling, but now at positive $\lam_2=3/2$ and $\lam_1=-9/4$ 
it loses stability to the wrinkling branch (samples D$_{1,2}$), while the 
bifurcation to buckling and coiling is at a later $(\lam_1,\lam_2)\approx 
(-3.77,3.021)$ and yields initially unstable branches. 
The wrinkling is initially stable, but 
at $\lam_2\approx 1.6$ loses stability to ``oblates embedded in wrinkling'', 
see sample E, and Remark \ref{vrem}. 

\begin{figure}[ht]
  \begin{center}
\btab{l}{
\btab{ll}{
{\small (a), from Fig.\ref{cmf1}D$_2$}\\
 \hs{-4mm}\ig[height=0.24\tew]{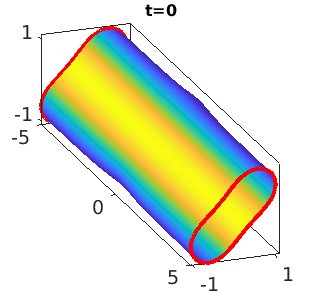} 
\hs{-1mm}\ig[height=0.24\tew]{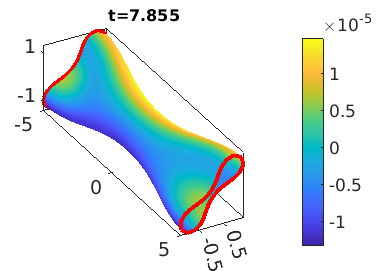} 
&\hs{-2mm}\ig[height=0.24\tew]{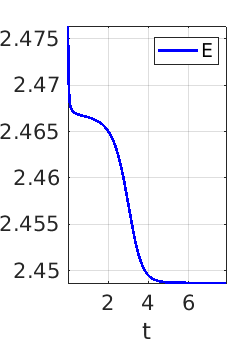}
\hs{-2mm}\ig[height=0.24\tew]{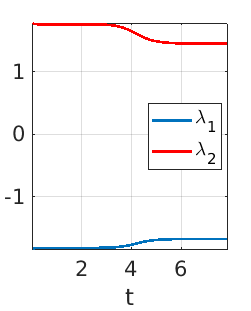}\\
{\small (b),  from Fig.\ref{cmf1}B}\\
 \hs{-4mm}\ig[height=0.24\tew]{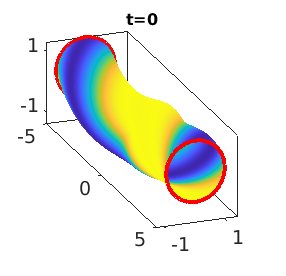} 
\hs{1mm}\ig[height=0.24\tew]{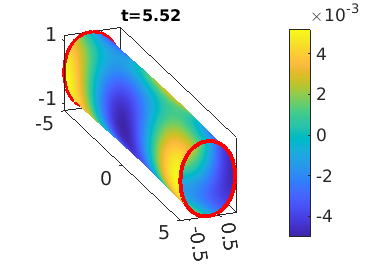} 
&\hs{-2mm}\ig[height=0.24\tew]{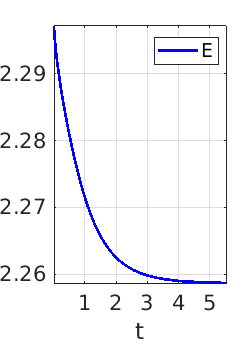}
\hs{-2mm}\ig[height=0.24\tew]{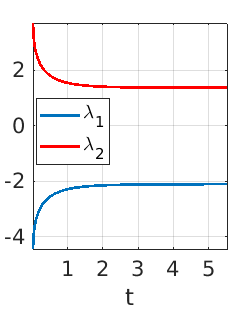}\\
{\small (c),  from Fig.\ref{cmf1}C}\\
 \hs{-4mm}\ig[height=0.24\tew]{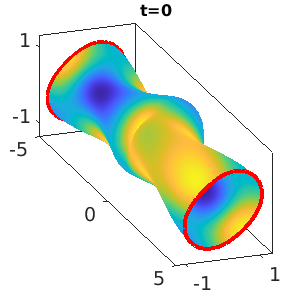} 
\hs{1mm}\ig[height=0.24\tew]{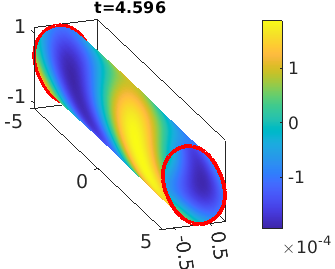} 
&\hs{-2mm}\ig[height=0.24\tew,width=27mm]{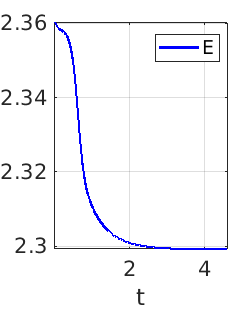}
\hs{-2mm}\ig[height=0.24\tew]{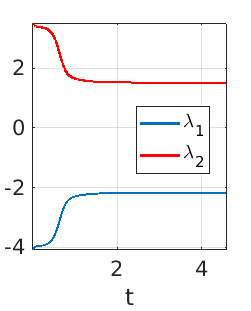}
 }
}
\vs{-0mm}
  \end{center}
  \caption{{\small Numerical flows near states from 
Fig.\ref{cmf1}; initial perturbations of type \reff{flp1} with $\del=-0.25$. 
Flow to ``embedded oblate'', here at the boundary, in (a), 
and to wrinkling in (b,c). \label{cmf2}}}
\end{figure}

\begin{figure}[ht]
  \begin{center}
\btab{l}{
\ig[width=0.25\tew]{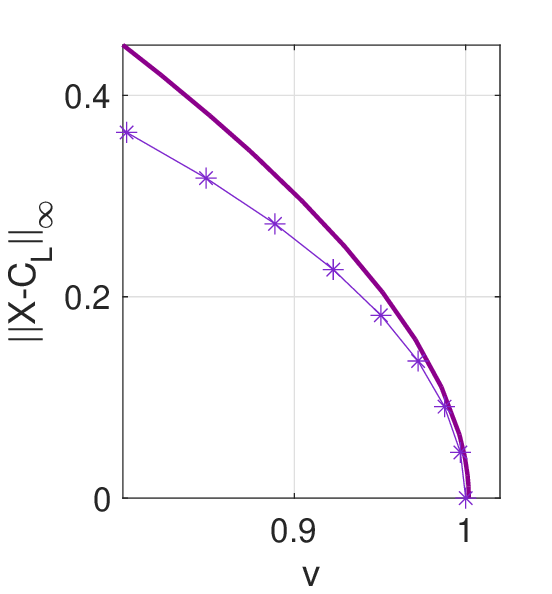} \ig[width=0.25\tew]{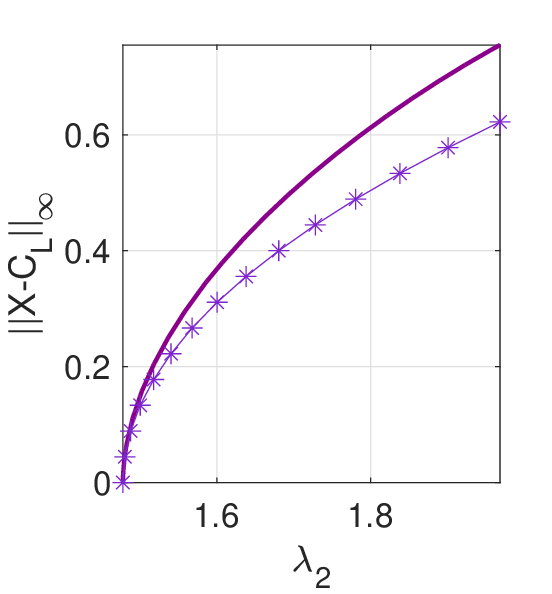} 
\ig[width=0.25\tew]{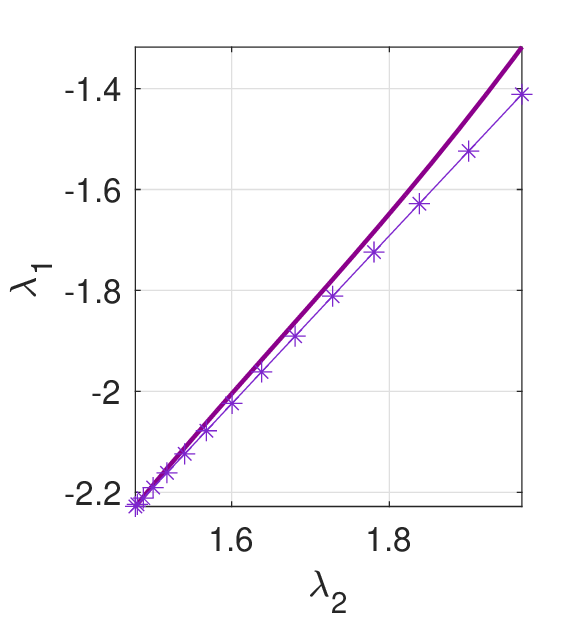} 
} 
\vs{-0mm}
  \end{center}
  \caption{{\small 
Comparison with AE prediction for primary wrinkling.  \label{cmf1b}}}
\end{figure}

Figure \ref{cmf2} shows some numerical flows from states from Fig.\ref{cmf1}, 
focusing on wrinkling perturbed according to \reff{flp1} 
(pearling is again stable, also under rather large perturbations). 
D$_2$ at its given $(\CA_0,\CV_0)$ is linearly (weakly) unstable, 
and the perturbation yields a slow evolution to an oblate wrinkling, 
with a very small change in $E$. In (b) and (c) we perturb a coiling and 
a buckling, and in both cases flow to wrinkling, and the bottom 
line of Fig.\ref{cmf2} and similar 
experiments is that (stable) wrinkling and oblate wrinkling 
have rather large domains of attraction. 
For completeness, Fig.\ref{cmf1b} shows a comparison with 
the AE predictions from \reff{wria1} for wrinkling, again with very good 
agreement, cf.~Remark \ref{arem1}.

\subsection{$c_0=0.48$}
\huc{We now consider 
$c_0\in(1/4,1/2)$}, specifically $c_0=0.48$, cf.~Fig.\ref{neigs}(c). 
By Corollary \ref{destabcor}, 
the cylinder destabilizes in a {\em finite wavelength} pearling instability, 
and in a long wave instability to coiling and 
buckling. The finite pearling period is $\ell:=\ell_p\approx 7.5$, and thus 
we will be mostly interested in length $L$ as multiples of $7.5$. 

\begin{figure}[ht]
\btab{l}{ 
\btab{ll}{ 
{\small (a)}&{\small (b)}\\
\hs{-5mm}\ig[height=0.29\tew]{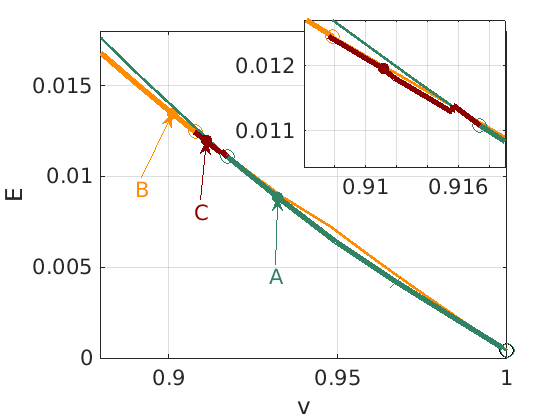}
\hs{-5mm}\ig[width=0.26\tew]{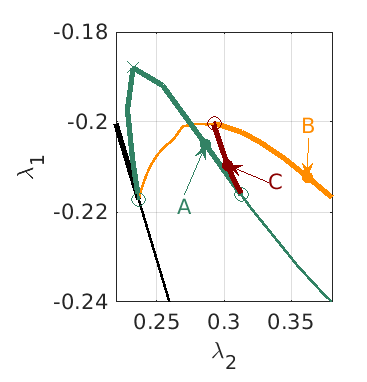}
&\rb{0mm}{\hs{-0mm}\ig[width=.21\tew]{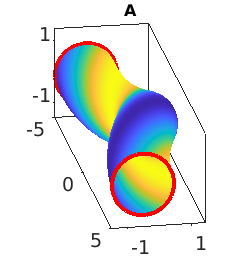}}
\rb{0mm}{\hs{-3mm}\ig[width=.18\tew]{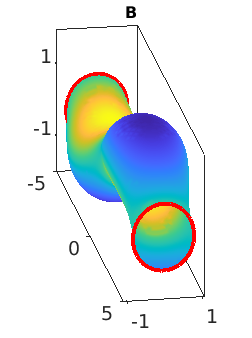}}
}\\[-2mm]
\btab{ll}{{\small (c)}&{\sm (d)}\\
\rb{0mm}{\hs{-0mm}\ig[width=.2\tew]{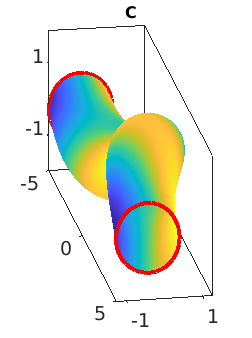}}
&\ig[width=0.25\tew]{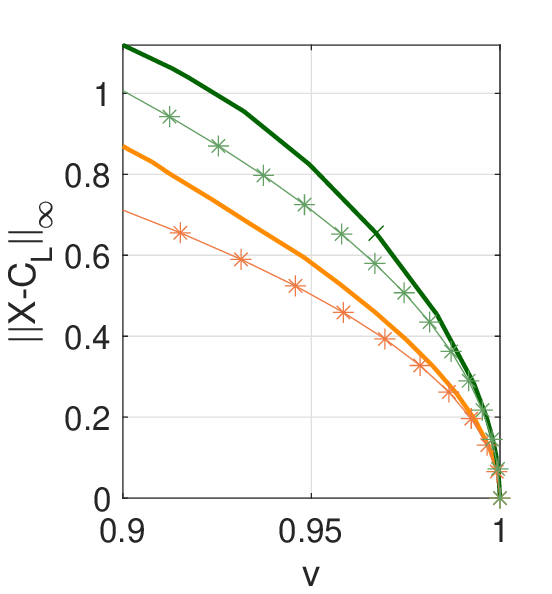} 
\ig[width=0.25\tew]{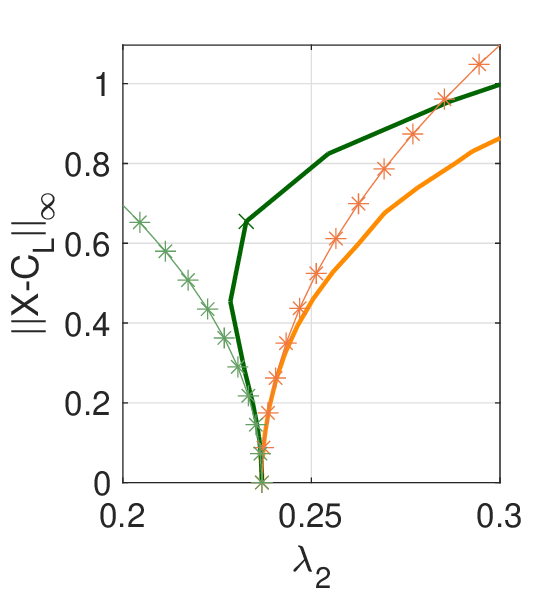} 
\ig[width=0.25\tew]{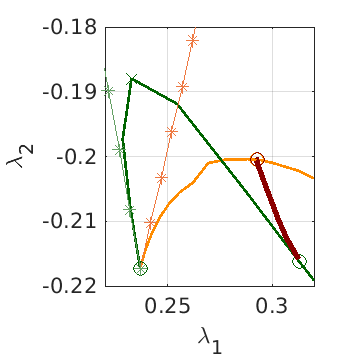} 
}
}

\vs{-0mm}

\caption{{\small  $c_0=0.48$, $L=10$. Numerical BDs and samples in (a,b,c), 
where C shows a stable connection between the stable regimes of coiling and 
buckling.
Comparison to AE predictions in (d). 
}}\label{c005f2}
\end{figure}

However, for comparison with the cases $c_0=0$ and $c_0=-1$ 
we first again choose $L=10$, and in Fig.\ref{c005f2} give the BD of buckling 
and coiling. The coefficients for the AE description \reff{ampsys} now 
are $(a,b_1,b_2)=(0.39,-0.083, 0.56)$, and in agreement with Fig.\ref{stabfig} 
the buckling is now stable at bifurcation and the coiling is unstable. 
Moreover, we find interesting secondary bifurcations on both branches, 
where the buckling loses (coiling gains) stability, and these two BPs 
are connected by a stable intermediate branch, see sample C. The comparison 
with the AEs in (c) again shows the agreement close to bifurcation, but 
naturally the (cubic order) AEs cannot capture the fold of the coiling  
in $\lam_2$. 

\begin{figure}[ht]
  \begin{center}
\btab{l}{
\btab{ll}{{\small (a$_1$)}&{\sm (a$_2$)} \\
\hs{-4mm}\ig[height=0.32\tew]{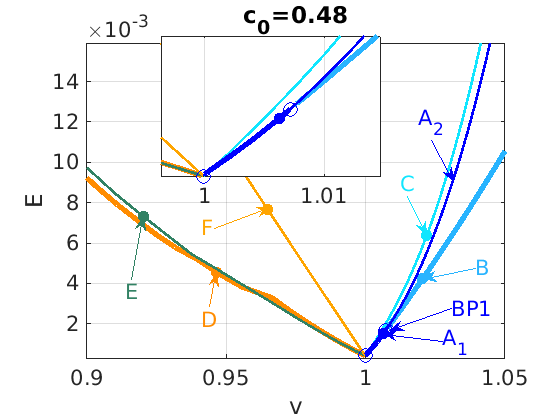} 
&\hs{-2mm}\ig[height=0.32\tew]{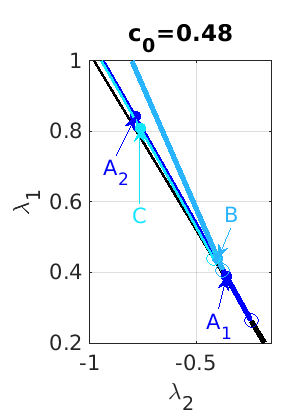}  
\hs{-1mm}\ig[height=0.32\tew]{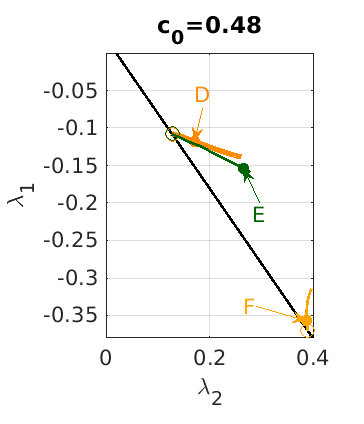} 
}\\
{\small (c)}\\
\hs{-3mm}\ig[height=0.32\tew]{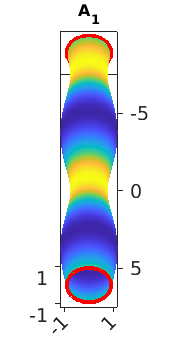}\hs{-3mm}\ig[height=0.32\tew]{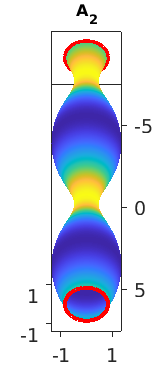}
\hs{-2mm}\ig[height=0.32\tew]{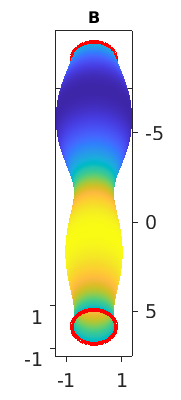}\hs{-2mm}\ig[height=0.32\tew]{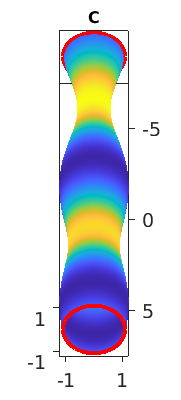}
\hs{-3mm}\ig[height=0.32\tew]{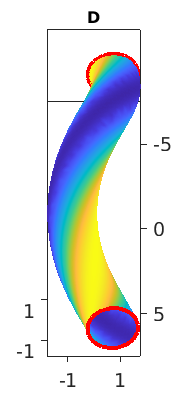}\hs{-0mm}\ig[height=0.32\tew]{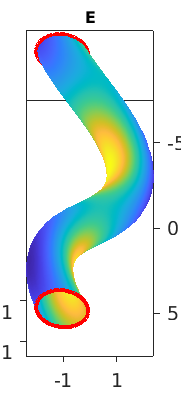}
\hs{-0mm}\ig[height=0.32\tew]{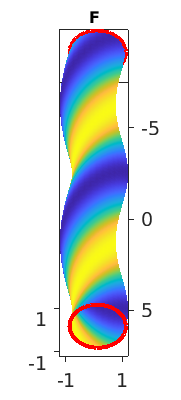}
}
\vs{-0mm}
  \end{center}
  \caption{{\small Basic BD for $c_0=0.48$, $L=15$. (a) Energy $E$ over 
reduced volume $v$, and $\lam_1$ over $\lam_2$. (b) Samples from 
primary pearling, A$_1$ stable, A$_2$ unstable. (c) Further samples: 
B is period doubling from A branch, stable; C is the primary $m=1$ 
pearling bifurcation, unstable; D is a (1,1)--coiling, stable; 
E is a (1,1)--buckling, unstable; F is a (2,1)--coiling. 
\label{c05f1}}}
\end{figure}

\begin{figure}[ht]
  \begin{center}
\btab{ll}{
\hs{-3mm}\ig[height=0.3\tew]{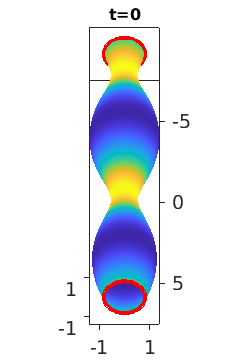}\hs{-2mm}\ig[height=0.3\tew]{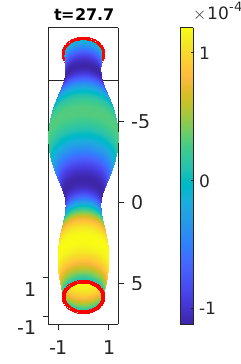}
\hs{-0mm}\ig[height=0.3\tew]{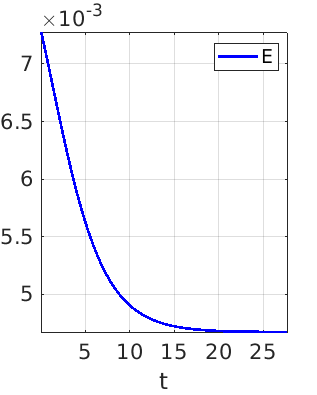}\ig[height=0.3\tew]{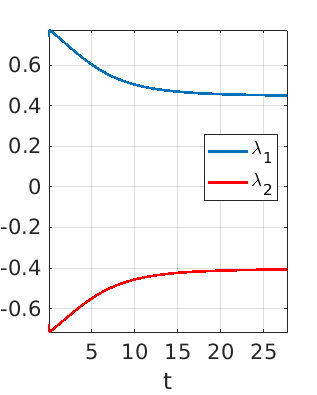}
}
\vs{-0mm}
  \end{center}
  \caption{{\small Sample flow from perturbation of A$_2$ of type 
(\ref{flp2}) with $\del=0.5$, going to B type solution.\label{c05f2}}}
\end{figure}

\begin{figure}[ht]
  \begin{center}
\btab{l}{
\btab{ll}{{\small (a)}&{\sm (b)} \\
\hs{-7mm}\ig[height=0.32\tew]{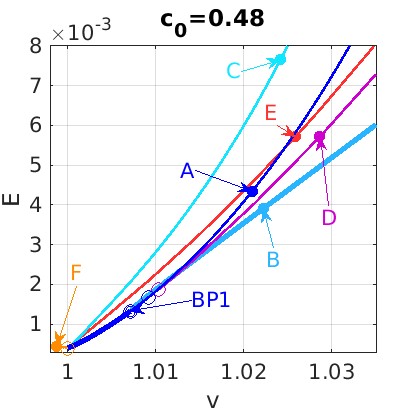} 
&\hs{-6mm}\ig[height=0.36\tew]{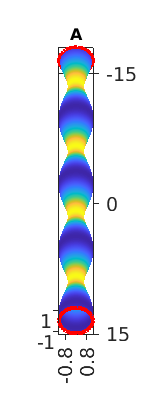}\hs{-2mm}\ig[height=0.36\tew]{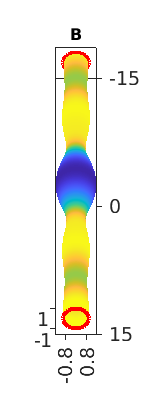}
\hs{-3mm}\rb{2mm}{\ig[height=0.34\tew]{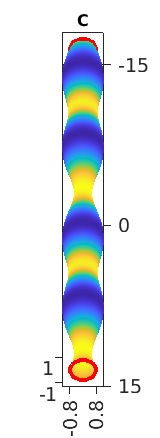}}\hs{-2mm}\ig[height=0.36\tew]{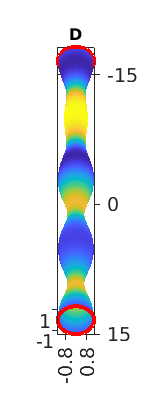}
\hs{-3mm}\ig[height=0.36\tew]{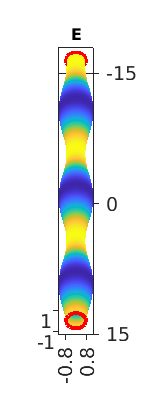}\hs{-2mm}\ig[height=0.36\tew]{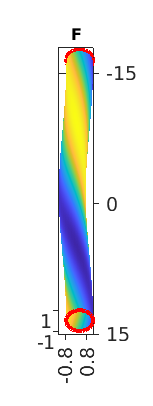}
}
}
\vs{-0mm}
  \end{center}
  \caption{{\small $c_0=0.48$, $L=30$. The $m=4$ (period $\ell=7.5$) 
pearling A now loses  
stability at BP1 much closer to its primary bifurcation, to 
a single ``localized'' pearling branch, sample B. Further branches: 
C: $m=2$ ($\ell=15$) primary pearling; E: $m=3$ ($\ell=10)$ primary pearling;
D: 2ndary pearling from 2nd BP on $m=1$ branch.  
F: primary $(m,n)=(1,1)$ coiling.   
\label{c05f1l}}}
\end{figure}

Figure \ref{c05f1} shows the BD and samples for $L{=}15$. 
To the right of $v{=}1$ bifurcates the $k_2{\approx}0.838$ pearling branch 
(at $\lam_2{\approx}-0.243$), and 
is initially stable (sample A$_1$). However, shortly after 
($\lam_2{\approx}-0.3756$) there 
is a secondary bifurcation to a stable branch with double the period, 
sample B. 
On the other hand, the primary $m{=}1$ pearling branch, sample C, 
(bifurcating at $\lam_2{\approx}-0.421$) is unstable throughout.%
\footnote{\huc{At first sight, C may look similar to the $m=2$ pearling A$_1$ and 
A$_2$, but the difference is that C has one weak and one stronger constriction, 
which intuitely makes C one period in a string of dumbbells; 
see also C in Fig.~\ref{c05f1l}, and Remark \ref{vrem}.}}
The $(m,n){=}(1,1)$ coiling (sample D) and buckling (sample E) 
bifurcate at $\lam_2{\approx}0.138$, 
and again the coiling is stable, and stays stable to rather small $v$, 
while the buckling is unstable. Wrinkling does not occur in the 
parameter regime considered; however, the 
(compared to $L=10$ before) slightly larger $L=15$  
lets the pearling and coiling/buckling BPs move closer together, which is 
why we chose 
to also present the $(m,n)=(2,1)$ coiling F (unstable). 
Figure \ref{c05f2} shows the flow from a perturbation of A$_2$ to a B 
type solution but rather close to the BP. 
Therefore, the convergence to steady state is again very slow, as 
near bifurcation the stable B branch has a leading eigenvalue close to $0$.

The main point of Fig.\ref{c05f1} is that the choice of $c_0=0.48$ 
together with $L=15$ yields the destabilization 
of the cylinder to the $k_2\approx 0.838$ pearling, but these solutions are only 
stable in a somewhat narrow range (in $v$, say). To further assess this, 
in Fig.\ref{c05f1l} we consider the primary pearling with $\ell=7.5$ 
on twice the domain $L=30$. 
The result shows that the range  of stable $\ell=7.5$ pearling 
shrinks with increasing $L$, and this continues for further increase of 
$L$, 
yielding the following dependence of the reduced volume $v_s$ of single 
pearling BPs from primary $\ell=7.5$ pearling on $L$: 
\huga{\label{ptt1}
\begin{tabular}{c|cccc}
$L$&15&39&45&60\\
$v_s$&1.0081&1.0063&1.0046&1.0035
\end{tabular}
}
We also ran some numerical flows from 
perturbations of A in Fig.\ref{c05f1l}, and similarly at larger $L$, 
with results as in Fig.\ref{c05f2}, i.e., convergence to the 
single pearling solutions. However, these DNS are very expensive, 
partly due to the larger domain and hence larger scale discretization, 
but more importantly due to the very slow convergence of the flows.  

\subsection{Outlook: $c_0=0.75$}\label{cpsec}
Lastly, following Fig.\ref{neigs}(e) we consider 
$c_0>1/2$, specifically $c_0=0.75$ for which 
$\CC_L$ is never stable if $L\ge 2\pi $ ($k_1\le 1$, say), and    
branches as before still bifurcate from $\CC_L$, 
but not at loss of stability of $\CC_L$. 
Hence they can also still be described by AEs as before, but they cannot 
be stable close to bifurcation, and we can only expect stable ``patterns'' 
at $O(1)$ amplitude. This is more difficult than 
the analysis and numerics so far, and hence here we only give 
an outlook on some interesting (numerical) results we obtain, 
again for a rather small length $L=10$ of the cylinder, and 
leave more systematic and detailed studies to later work.

\begin{figure}[ht]
  \begin{center}
\btab{l}{
\btab{ll}{{\small (a)}&{\sm (b)} \\
\hs{-4mm}\ig[height=0.32\tew]{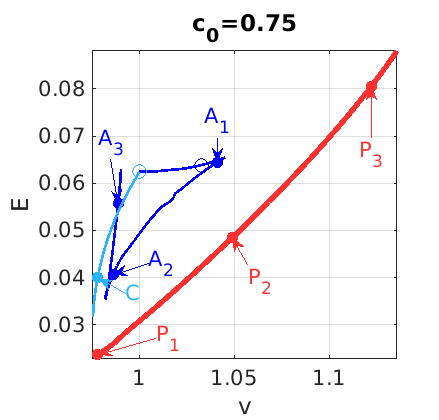} 
\hs{-2mm}\ig[height=0.32\tew]{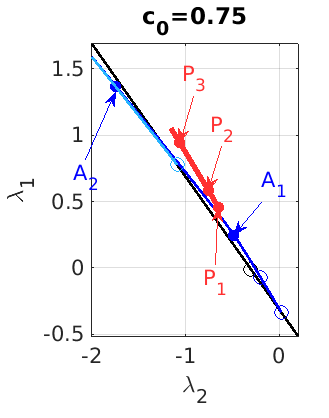}  
&\hs{0mm}\ig[height=0.3\tew]{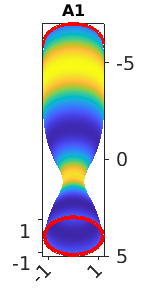} \hs{-1mm}\ig[height=0.3\tew]{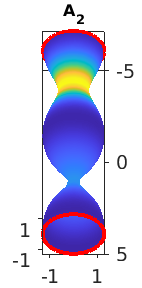} 
 \hs{-1mm}\ig[height=0.3\tew]{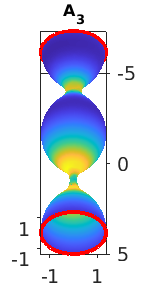} 
}\\
{\small (c)}\\
 \hs{1mm}\ig[height=0.22\tew]{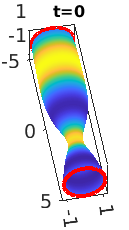} \hs{-1mm}\ig[height=0.22\tew]{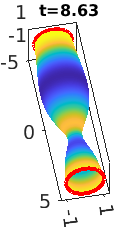} 
 \hs{-1mm}\ig[height=0.22\tew]{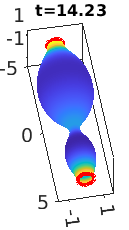} \hs{-0mm}\ig[height=0.22\tew]{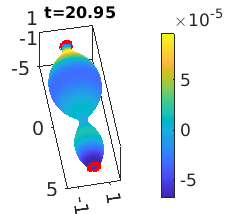} 
\rb{-3mm}{
\hs{-3mm}\ig[height=0.27\tew]{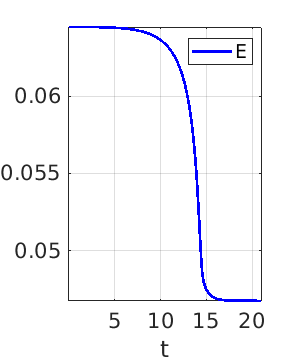} \hs{-1mm}\ig[height=0.27\tew]{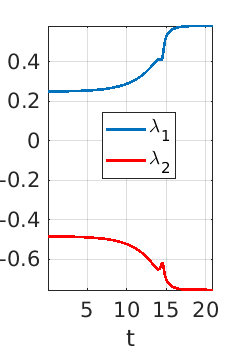}} \\
{\small (d)}\\
\hs{-3mm}\ig[height=0.2\tew]{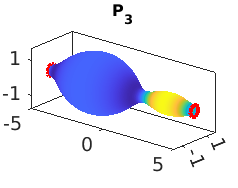} 
 \hs{-1mm}\ig[height=0.2\tew]{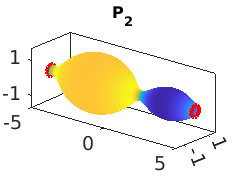} \hs{-1mm}\ig[height=0.2\tew]{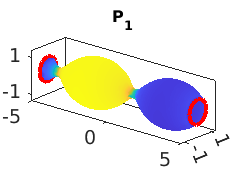} 
\hs{-1mm}\ig[height=0.2\tew]{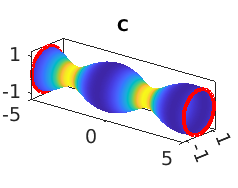} 
}
\vs{-2mm}
  \end{center}
  \caption{{\small Selected branches for $c_0=0.75$, $L=10$. 
(a) Energy $E$ over 
reduced volume $v$, and $\lam_1$ over $\lam_2$. (b) Samples from 
primary pearling (blue branch in (a)) all unstable. (c) DNS 
from near $A_2$, converging to a stable pearling $P_2$ with unequal pearl sizes. 
(d) Samples from the stable P branch, and sample from (unstable) 
primary $m=2$ pearling branch (light blue in (a)). 
\label{c1f1}}}
\end{figure}

Figure \ref{c1f1} shows the BD of $k_1\approx 0.628$ pearling 
(A, blue) and $k_2\approx 1.257$ pearling (C, light blue), 
which are both unstable at bifurcation as they must, and 
moreover are unstable throughout, 
and addtionally an $O(1)$ amplitude pearling branch (P, red), 
where P$_2$ is obtained from DNS  starting 
near A$_1$ (panel (c)). In summary, we see a number of new effects which 
strongly suggest further study, but here we restrict to the following comments.

$\bullet$ We stop the continuations shortly after A$_3$, shortly 
after C, and shortly after P$_1$ as further continuations there 
require excessive mesh adaptations. We can continue significantly 
beyond P$_3$, but our main region of interest is near $v=1$. 
There are bifurcations on the blue branches which break the 
rotational symmetry, but we restrict to presenting axisymmetric branches. 

$\bullet$ The $k_1$ pearling branch A shows two folds in $v$, 
the first at $v_{f,1}\approx 1.045$, the second at 
$v_{f,2}\approx 0.98$. In particular, between A$_1$ and A$_2$ 
it goes through $v=1$ again, while keeping the area $\CA$ 
constant throughout. This is not a contradiction to the fact 
that normal variations of $\CC_L$ that fix $\CA$ and $\CV$ 
locally do not exist; the passage through $v=1$ here 
corresponds to a non--small normal variation. 

$\bullet$ The folds on the 
A branch naturally call for fold continuation in $c_0$, in 
particular for decreasing $c_0$ as for $c_0<1/2$ we do not see 
folds on the analogous A branches. We believe that they vanish 
at some $c_0>1/2$ as the blue swallowtail shrinks and the folds 
collide, but this remains to be studied. 

$\bullet$ The flow in (c) starts at A$_1$ with a small pertubation 
in direction of the most unstable eigenvalue; 
A$_1$ is below the fold, and there is exactly one 
positive eigenvalue $\mu\approx 0.2$. In any case, 
the results of flow experiments seem rather robust wrt 
the choice of the initial conditions (choice of steady state 
and perturbation), i.e., we always ``land'' on the P branch, 
though naturally at different $v$ and $\lam_2$, or analoga 
of the P branch at different area $\CA$. In the particular 
case shown, the perturbation of A$_1$ yields a relative 
change of $\CV$ of order $10^{-3}$, and a relative change of 
$\CA$ of order $10^{-4}$. (The initial perturbation 
{\em linearly} conserves $\CA$ and $\CV$, i.e., 
is tangential to the constraint manifold, but 
for any finite perturbation we must expect a (nonlinear) 
change of $\CA$ and $\CV$. However, the change in $\CA$ is 
so small that we ignore it and keep the normalization of 
$E$ in \reff{enorm} by the original $\CA$.) 

$\bullet$ The last plots in (c) suggest almost steady states for $t>20$, 
and we use the solution at $t=20.95$ as an initial guess for 
Newton's method for the steady problem. This indeed takes us to the 
(spectrally) stable steady state P$_2$, and continuing 
from here to larger/smaller $v$ yields the red branch. 

$\bullet$ It is remarkable that the P branch with two pearls of different 
sizes  is stable, that along the branch the pearls change their relative 
sizes (different from the primary pearling branches, where only 
the global amplitude changes), 
and that in fact P$_2$ attracts (perturbations of) A$_1$.  
However, some finite size effects are involved here, because some stability 
of the P branch is lost when 
repeating the experiment from Fig.\ref{c1f1} on larger domains, 
e.g., $L=20$. In this case, initial conditions for numerical flows 
must be chosen more carefully, namely closer to $v=1$ to reach the analogue 
of the P branch (4 pearls on $L=20$), and continuation of the P branch 
then shows a loss of stability (wrt to still axially symmetric solutions) 
away from $v=1$. 
We plan to study in particular these effects 
in a systematic way in future work, including their relation 
to the transition from A to B in Fig.~\ref{c05f1l}. 

$\bullet$ Shape wise (compare P$_1$ and C in (d)) it seems possible 
that the P branch connects to the C branch in a period--doubling 
bifurcation, but samples  P$_1$ and C are at rather different values of 
$\lam_{1,2}$. In our preliminary experiments so far we could not 
identify to what  other branches the P branch connects, and again we 
leave this for future work.


\section{Discussion}\label{dsec}
Motivated by experimental results as in Fig.\ref{f0}(a,b), 
we studied the bifurcations from straight cylinders $\CC_L$ of lengths $L$ 
(and wlog radius $r=1$) 
in the Helfrich model \reff{hflow}. After showing well--posedness 
we discussed the linearization at $\CC_L$ and derived the amplitude 
equations for the bifurcations of pearling, wrinkling, and coiling and buckling. 
Subsequently we used numerical continuation and bifurcation methods 
to extend the local bifurcations to a more global picture. 

Our results initially suggest that \reff{hflow} besides being 
very interesting mathematically is a possible 
model for the experiments in Fig.\ref{f0}(a,b). However, 
\huc{when restricting to branches bifurcating at loss 
of stability of $\CC_L$, 
there also are  
strong discrepancies. One is that coiling (and buckling, which as far as 
we know is not clearly documented in experiments) always bifurcate 
with the maximum wavelength 
$L$ allowed by the domain. For $c_0>1/4$, pearling can occur with a finite 
wavelength (Figs.\ref{c05f1}--\ref{c05f1l}), but the stability 
range of finite wavelength pearling shrinks with increasing $L$, 
probably to zero for $L\to\infty$. This means that the experimental 
results in Fig.\ref{f0}(a,b) are not described by the primary 
 branches bifurcating at loss of stability of $\CC_L$. 
On the other hand, when further increasing $c_0$ beyond $c_0{=}1/2$ we find 
stable finite wavelength $O(1)$ amplitude pearling, but we only gave 
an outlook on this in \S\ref{cpsec}, and this clearly needs to be studied 
more systematically and in more detail, and for larger $L$.} 

Mathematically, the bifurcations of 
localized pearling from finite wave number pearling, with the pertinent 
BP moving closer to the primary bifurcation for larger domain size $L$, 
is similar to the bifurcation of 
localized patterns from periodic patterns in various 
dissipative systems, e.g., reaction diffusion (RD) systems or 
Swift--Hohenberg (SH) type equations. See, e.g., 
\cite{bukno2007,hokno2009,BKLS09,kno2015}, and 
\cite[Chapters 8 and 9]{p2pbook} for further references and 
numerics using \pdep. 
In these systems, the localized patterns are usually attributed to 
subcritical primary bifurcations of periodic patterns, 
which are hence unstable, but stabilize after a fold 
at large amplitude. The localized patterns 
then show some ``snaking'' behavior with several folds in which 
the pattern grows, and with alternating stable and unstable segments 
of this snaking branch. Moreover, in the limit $L\to\infty$, 
the BP of the snaking 
branch converges to the primary BP of the periodic patterns. 

The data in \reff{ptt1} for BP1 moving closer  to $v=1$ for larger $L$ 
suggest a similar effect, but there are clear differences 
to the snaking scenario, and many open questions. In our numerics: 
\bci 
\item After the loss of stability of the primary pearling A--branch, 
it never regains stability; in particular, we cannot find a 
fold (in $v$, say) for the A--branch. 
\item Similarly, we cannot find a fold or other loss of stability of 
the localized pearling B--branch. 
\eci 

The numerics here are harder for large $L$ than for 
periodic or localized 1D (or 2D) patterns in, e.g., reaction diffusion 
systems, 
and for instance for $L=60$ we only continue the A--branch and the 
B--branch 
to $v=1.01$, such that these results are in some sense still local. 
\huc{As mentioned in \S\ref{cpsec}, one way to proceed should be 
fold continuation in $c_0$ of the folds that we find on the A branch 
at larger $c_0$.}   

All our numerics and associated analysis was for finite (and mostly 
rather small) $L$. 
For linearizations of (classical) dissipative systems on intervals $(0,L)$ 
at some steady state, the spectral gap 
between the critical eigenvalues and the rest of the spectrum shrinks 
as $L\to\infty$, and over $\R$ the 
amplitude equations (ODEs) for finite $L$ must be replaced by 
modulation equations (PDEs) of Ginzburg--Landau type, incorporating 
the continuous curve of critical eigenvalues. See, 
e.g., \cite[Part IV]{SU17}, and the references therein. For the Helfrich 
cylinders it is not clear how to proceed for $L\to\infty$. 
For long wave instabilities 
(buckling, coiling, wrinkling, and pearling except in a small range of positive 
$c_0$ as in Fig.\ref{c05f1}--\ref{c05f1l}), without the volume constraint 
the critical curve $\R\ni k\mapsto \mu(k,n)$ is maximal at wave number 
$k=0$, and we believe that instead of \reff{hflow} 
a better way to incorporate the volume constraint is to consider 
the $H^{-1}$ gradient flow for \reff{helfen} with area constraint. 
This automatically conserves $\CV$, and the critical 
curve is changed to $k\mapsto \tilde\mu(k,n)=k^2\mu(k,n)$,  
such that we expect modulation equations with an (additional) conservation 
law; see, e.g., \cite{MC00,ZS16,Hi20}. 
This, however, is also beyond the scope of the present paper.

\appendix 
\section{Buckling and coiling amplitude system computations}\label{aeapp}
To derive the coefficients $a,b_1,b_2$ for the amplitude system \reff{ampsys} 
for the mixed mode instability we as usual plug the extended ansatz 
\reff{mma1} into \reff{bifprob} and sort wrt powers of $\eps$. 
 \bci
 \item[$\eps \psi_{11}$:] $\ds -\frac{1}{2}k^{2} \left(k^{2}-4 c_{0}-2 \lam_{2}+2\right)A=0.$
  \item[$\eps \psi_{1-1}$:] $\ds -\frac{1}{2} k^{2} \left(k^{2}-4 c_{0}-2 \lam_{2}+2\right) B=0$. 
  \item[$\eps^2\psi_{00}$:]  $\ds 0=\left(-\frac{1}{2}+\frac{k^{4}}{2}+\frac{\left(5-8 c_{0}-2 \lam_{2}\right) k^{2}}{2}-\lam_{2}\right) ({| A |}^{2}+ {| B |}^{2})
  -(\lam_{2}+1/2) A_{00} -\beta_{2}-\alpha_{2}.$

  
  \item[$\eps^2\psi_{22}$:] $0=\left(-\frac{7 k^{4}}{4}{+}\frac{\left(-8 c_{0}{+}2 \lam_{2}{+}19\right) k^{2}}{4}{-}\frac{3 \lam_{2}}{2}{+}\frac{9}{4}\right) A^{2}+\left(-8 k^{4}{+}\frac{\left(16 \lam_{2}{+}32 c_{0}{-}64\right) k^{2}}{4}{+}3 \lam_{2}{-}\frac{9}{2}\right) A_{22}.$
  \item[$\eps^2\psi_{2-2}$:] $0=\left(-\frac{7 k^{4}}{4}+\frac{\left(-8 c_{0}{+}2 \lam_{2}{+}19\right) k^{2}}{4}-\frac{3 \lam_{2}}{2}{+}\frac{9}{4}\right) B^{2}+\left(-8 k^{4}{+}\frac{\left(16 \lam_{2}{+}32 c_{0}{-}64\right) k^{2}}{4}+3 \lam_{2}{-}\frac{9}{2}\right) B_{22}$. 
  
   \item[$\eps^2\psi_{20}$:] $0=\left(-\frac{7 k^{4}}{2}{+}\frac{\left(2 \lam_{2}+8 c_{0}+3\right) k^{2}}{2}{-}\lam_{2}{-}\frac{1}{2}\right) B A +\left(-8 k^{4}+\frac{\left(16 c_{0}{+}8 \lam_{2}\right) k^{2}}{2}{-}\lam_{2}{-}\frac{1}{2}\right) A_{20} $. 
   
      \item[$\eps^2\psi_{02}$:] $0=4 A \left(\frac{k^{4}}{8}+\left(c_{0}-\frac{\lam_{2}}{4}+\frac{5}{8}\right) k^{2}-\frac{3 \lam_{2}}{4}+\frac{9}{8}\right) \overline{B}+\left(3 \lam_{2}-\frac{9}{2}\right) A_{02}$. 
      
       \item[$\eps^3\psi_{11}$:] 
{\sm \hualst{\hs{-10mm}\ddT A =& \bigg[(-k^{2} \alpha_{2}-4 \left(-\frac{5 k^{6}}{8}+\left(c_{0}+\frac{3 \lam_{2}}{8}-2\right) k^{4}+\left(-\frac{c_{0}}{2}+\frac{5 \lam_{2}}{8}+\frac{23}{16}\right) k^{2}-\frac{3 \lam_{2}}{4}+\frac{9}{8}\right) {| A |}^{2}\\
       &{-}8 \left({-}\frac{5 k^{6}}{8}{+}\left(c_{0}{+}\frac{3 \lam_{2}}{8}\right) k^{4}{+}\left(\frac{3 c_{0}}{2}{-}\frac{3 \lam_{2}}{8}{+}\frac{15}{16}\right) k^{2}{-}\frac{3 \lam_{2}}{4}{+}\frac{9}{8}\right) {| B |}^{2}{-}4 k^{2} \left(\frac{1}{2} c_{0}{-}\frac{3}{8}  \right) A_{00}\bigg] A \\
       &{+}6\left( \left(c_{0}{+}\frac{1}{4}\right) k^{2}{-} \left(\lam_{2}{-}\frac{3}{2}\right)  \right) A_{02}B {-}10  \left(-\frac{2 k^{4}}{5}{+}\left(c_{0}{+}\frac{\lam_{2}}{5}{-}\frac{49}{20}\right) k^{2}{+}\frac{3 \lam_{2}}{5}{-}\frac{9}{10}\right) A_{22}\overline{A}\\
       &-2  k^{2} \left(-2 k^{2}+c_{0}+\lam_{2}-\frac{1}{4}\right) A_{02}\overline{B}.}
}
    \item[$\eps^3\psi_{1-1}$:] 
{\sm \hualst{\hs{-10mm}\ddT B=&\bigg[{-}k^{2} \alpha_{2}{-}4 \left({-}\frac{5 k^{6}}{8}{+}\left(c_{0}{+}\frac{3 \lam_{2}}{8}{-}2\right) k^{4}{+}\left(-\frac{c_{0}}{2}{+}\frac{5 \lam_{2}}{8}{+}\frac{23}{16}\right) k^{2}{-}\frac{3 \lam_{2}}{4}{+}\frac{9}{8}\right) {| B |}^{2}\\
     &{-}8 \Big({-}\frac{5 k^{6}}{8}{+}\left(c_{0}{+}\frac{3 \lam_{2}}{8}\right) k^{4}{+}\left(\frac{3 c_{0}}{2}{-}\frac{3 \lam_{2}}{8}{+}\frac{15}{16}\right) k^{2}{-}\frac{3 \lam_{2}}{4}{+}\frac{9}{8}\Big) {| A |}^{2}
   {-}2 k^{2} \left(c_{0}{-}\frac{3}{4} \right)A_{00}\bigg] B \\
   &+6  \left(\left(c_{0}{+}\frac{1}{4}\right) k^{2}{-}\lam_{2}{+}\frac{3}{2}\right) A\overline{A_{02}}{-}10 \bigg({-}\frac{2 k^{4}}{5}{+}\left(c_{0}{+}\frac{\lam_{2}}{5}{-}\frac{49}{20}\right) k^{2}
     +\frac{3 \lam_{2}}{5}-\frac{9}{10}\bigg) B_{22} \overline{B}\\
     &{-}2  k^{2} \left({-}2 k^{2}{+}c_{0}{+}\lam_{2}{-}\frac{1}{4}\right) A_{20} \overline{A}.}}
 \eci
From the area constraint we get, at $\eps^2\psi_{00}$, 
$0=\left(k^{2}+1\right) (| A |^{2}+| B |^{2})+A_{00}$, 
and using $\lam_2^{(1,1)}=\frac 1 2 k^2-2c_0+1$ we solve the algebraic equations for 
the $\eps^2$ modes: 
 \hualst{&A_{00}=-\frac{\left(4 k^{2} c_{0}-2 k^{2}-4 c_{0}+3\right) }{k^{2}-4 c_{0}+3}(| A |^{2}+|B|^2)-\frac{2 \alpha_{2}+2 \beta_{2}}{k^{2}-4 c_{0}+3},\\
 &A_{22}=-\frac{ \left(2 k^{4}+4 k^{2} c_{0}-6 k^{2}-4 c_{0}-1\right)}{2 \left(4 k^{4}+7 k^{2}+4 c_{0}+1\right)}A^{2},\quad 
B_{22}=-\frac{ \left(2 k^{4}+4 k^{2} c_{0}-6 k^{2}-4 c_{0}-1\right)}{2 \left(4 k^{4}+7 k^{2}+4 c_{0}+1\right)}B^{2},\\
&A_{20}=-\frac{ \left(6 k^{4}-4 k^{2} c_{0}-4 k^{2}-4 c_{0}+3\right)}{12 k^{4}-7 k^{2}-4 c_{0}+3}AB, \quad 
A_{02} =-\frac{ \left(4 k^{2} c_{0}+4 c_{0}+1\right) }{k^{2}-4 c_{0}-1}A\overline{B}. 
}
For the Lagrange multiplier expansion we get 
 \hualst{\al_2= \frac{k^{2} \left(k^{2}-8 c_{0}+6\right) }{2}{| A |}^{2}+\frac{k^{2} \left(k^{2}-8 c_{0}+6\right) }{2}{| B |}^{2}-\beta_{2}, }
and altogether this gives the solvability condition \reff{ampsys} for $A$ and 
$B$ 
at $\eps^3$ with 
\hualst{a=&k^2, \\
b_1=&\frac{5}{4} k^{6}+5 k^{4} c_{0}+\frac{3}{4} k^{4}+9 k^{2} c_{0}-\frac{33}{4} k^{2}-6 c_{0}-\frac{3}{2}\\ 
&+\frac{3 \left(2 k^{4}+4 k^{2} c_{0}-6 k^{2}-4 c_{0}-1\right) \left(-\frac{1}{2} k^{4}+k^{2} c_{0}-\frac{13}{4} k^{2}-2 c_{0}-\frac{1}{2}\right)}{4 k^{4}+7 k^{2}+4 c_{0}+1},\\
b_2=&3 k^{6}+4 k^{4} c_{0}-6 k^{4}-16 k^{2} c_{0}-3 k^{2}-12 c_{0}-3 \\
&+\frac{6 \left(4 k^{2} c_{0}+4 c_{0}+1\right) \left(-\left(c_{0}+\frac{1}{4}\right) k^{2}+\frac{k^{2}}{2}-\frac{1}{2}-2 c_{0}\right) }{k^{2}-4 c_{0}-1}\\
&+\frac{2 \left(6 k^{4}-4 k^{2} c_{0}-4 k^{2}-4 c_{0}+3\right) k^{2} \left(-\frac{3 k^{2}}{2}-c_{0}+\frac{3}{4}\right) }{12 k^{4}-7 k^{2}-4 c_{0}+3}. 
}
\section{Numerical algorithms}\label{num}
The basic numerical setup for continuation in the differential geometric 
setting is explained in \cite{geompap}, see also \cite{p2phome} for 
software download and demos. 
Here we recall some basic ideas and comment on some extensions. 
The numerics are based on finite element 
discrete differential geometry operators, as 
for instance implemented in \gptool~ \cite{gpgit}.  
Let ({\tt X}, {\tt tri}) be a triangulation of the surface $X$ with 
node coordinates {\tt X}$_j\in\R^{1\times3}$, 
$j\in\{1,\ldots,n_p\}$ and triangle list {\tt tri}$\in\R^{n_t\times 3}$. 
Using linear hat functions $\phi_i:X\ra\R$ with $\phi_i(${\tt X}$_j)=\del_{ij}$ with $i\in\{1,\ldots,n_p\}$, the operator (cotangent Laplacian) 
\hueq{
L_{kl}=\int_X\nabla\phi_k\nabla\phi_l\dd S=-\frac 1 2(\cot \al_{kl}+\cot \beta_{kl}).
}
approximates the (positive definite) Laplace--Beltrami operator, while for 
the FEM mass matrix $M$ we use the so--called Voronoi mass matrix. 
Then, e.g., $MH(X)=\frac 1 2 \spr{LX,N}$, where the vertex normals $N$ 
are obtained from interpolating the normals of the adjacent triangles. 
For the Gaussian curvature $K$ we use a discrete version of Gauss--Bonnet, 
and refer to \cite{BGN19} for more sophisticated schemes. 
See \cite[\S2.2]{geompap} for comments on and illustrations of 
convergence of the used FEM. 

Technically, for the 4th order problem 
\reff{hflow} we use a mixed formulation as two 2nd order problems, 
first analyzed in \cite{Rusu05} (in a more complicated scheme); 
see \cite{SGJW22} for a convergence analysis in the frame of the 
biharmonic equation. Additionally to the PDE $\eta(u)\pa_t u=G(u,\Lam)$ 
and the constraints $q(u)=0$, the continuation needs phase conditions (PCs), 
for which we introduce additional Lagrange multipliers $s_x,s_y,s_z$ and 
$s_{rx}$ pertaining to translations in $x,y$ and $z$, respectively, 
and rotations around the $x$--axis. The $y$--and $z$--translation 
PCs are always on, the $x$--translations is on for pearling and coiling, 
the $x$--rotation is on for the wrinkling, and $x$--translation {\em and} 
rotation are on for the buckling. In all continuations, the associated Lagrange 
multipliers stay $\CO(10^{-6})$. For the discrete systems, an effect of 
all the constraints (also of $q(u)=0$) is that the linear systems 
(to be solved in Newton loops and for eigenvalue computations) are 
bordered, with border widths up to six rows/columns, and we usually 
use bordered elimination \cite[\S4.5.1]{p2pbook} for their solution.

\paragraph{Mesh handling} 
The initial cylinder in discretized by a standard mesh from \pdep, mapped 
via a local map on the cylinder. The number of nodes in $x$ and $\bphi$ 
direction is chosen to minimize the mesh--quality quotient 
\hueq{\label{mq}q=\max_{T\in\CT}\frac{e_M(e_1+e_2+e_3)}{2|T|},} 
with $e_j$ the edge lengths, $e_M$ the maximal edge length, and $|T|$ 
the triangle areas.  
The triangulation is oriented inwards, matching with our 
analytical convention. 
On the bifurcating branches triangles distort, due to the change of {\tt X}. 
We use mesh adaptation by refinement and coarsening of {\tt X} 
and {\tt tri} based on \reff{mq} and triangle areas. Namely,  
we refine if the area is above a threshold, typically 1.3 times the size 
of the original cylinder triangles. For the conservation of the periodic 
boundary, we enforce the same refinement at ``opposite'' boundary triangles. 
For the coarsening \reff{mq} is the control parameter, activated if $q>8$. 
We select triangles smaller than a threshold $A_t$, 
typically $A_t=\max_{T\in\CT}|T|/3$ 
and collapse acute triangles along the shortest edge, 
while obtuse triangles are handled with edge flipping. 

\paragraph{Periodicity in $x$} 
During the continuation of bifurcating branches, it is not necessary that the boundary of $X$ stays in a plane. 
However, for mesh refinement along the branch the legacy 
{\tt box2per} setting looks for periodic boundary nodes in planes.  
For that reason we only allow an $(y,z)$ displacement at the boundary, 
by setting $N_x=0$.

\paragraph{Extracting the amplitude} For the comparison between the 
numerics and the AE predictions, we evaluate the 
deviation at each mesh point. In each continuation step we obtain 
$X+uN$ as the new solution, but multiple steps result in 
tangential movement of mesh points relative to the cylinder. 
Therefore we project each mesh point back onto the cylinder 
orthogonally to 
$Cp = \left( p_1, \frac{p_2}{\|(p_2, p_3)\|}, \frac{p_3}{\|(p_2, p_3)\|} \right)^T.$  
From this projection, we can measure the normal deviation from the cylinder and compare it with the amplitude approximation. 

\begin{figure}[ht]
\bce
\ig[height=0.38\tew,width=0.55\tew]{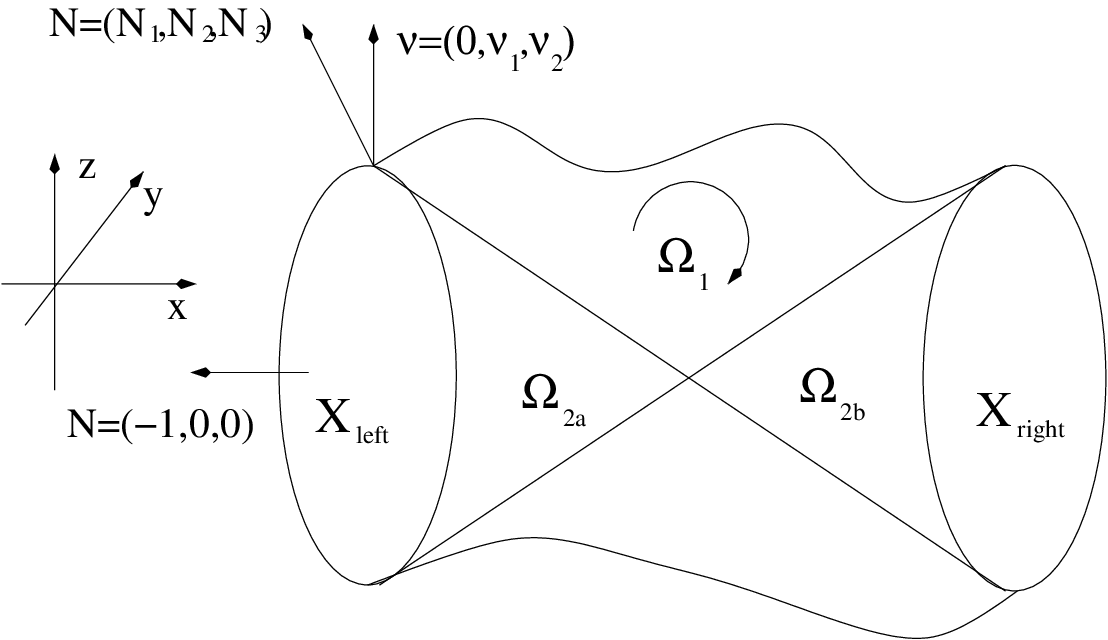}
\ece 
\vs{-0mm}
\caption{\small{Sketch for volume computation.}\label{fvol}}
\end{figure}

\paragraph{On the volume in the numerics}
To compute the volume of (not necessarily straight) 
cylinders $\CC$ we split them into two domains $\Om_1$ and $\Om_2$, 
see Fig.\ref{fvol}, where $\Om_2=\Om_{2,a}\cup\Om_{2,b}$ consists of 
two equal cones, from $(x,y,z)=(0,0,0)$ to the (in general not circular) 
left and right disk at $x=\pm L/2$. 
The idea is to compute $V=\vol(\Om_1)+2\vol(\Om_{2a})$ by the divergence theorem 
$\vol(\Om)=\frac 1 3 \int_\Om \mathrm{ div} (x) \dd x = \frac 1 3 \int_{\pa\Om} 
\spr{x,N}\dd S$. In the splitting $\CC=\Om_1\cup\Om_2$, we have 
$\spr{x,N}=0$ along the inner ``conical'' boundaries 
$\pa{\Om_1}\cap\pa{\Om_2}$, such that 
$\vol(\Om_1)= \frac 1 3 \int_{X} \spr{(x,y,z),N}\dd S$, and 
$\vol(\Om_{2a})=\frac 1 3 \int_{X_\text{left}}\spr{(-L/2,0,0)^T,
(-1,0,0)^T}\dd (y,z)=\frac L 2|X_\text{left}|$, 
and for the latter we can apply Gauss' theorem again to obtain 
\huga{ \vol(\Om_{2a}) = \frac L 6 \frac 1 2 \int_{\pa X_\text{left}}
\spr{(y,z),(\nu_1,\nu_2)}\dd s
=\frac L{12} \int_{\pa X_\text{left}}
\spr{(X_2,X_3),(N_2,N_3)/\|(N_2,N_3)\|}\dd s. }
Thus, altogether, 
\huga{\label{volf}
\vol(\CC)=\frac 1 3 \int_{X} \spr{(x,y,z),N}\dd S+\frac L{6}
\int_{\pa X_\text{left}}
\spr{(X_2,X_3),(N_2,N_3)/\|(N_2,N_3)\|}\dd s.
}

In the numerics we also need the derivative of $V(\CC+u N)$ wrt $u$. 
For $\Om=\Om_0+uN$ and $\vol(\Om)=\frac 1 3 \int_{\pa\Om} \spr{X,N}\dd S$ 
we have $\pa_u \vol(\Om)v=\int_\Om v\dd S$, which can be used 
for $\pa_u \vol(\Om_1)$, and for $\pa_u\vol(\Om_{2a})$ we proceed as follows. 
We choose a positive oriented arc length parametrization 
$\ga:[0,2\pi)\ra \pa X_{\text{left}}$. 
Then the (unit) tangent and (unit) normal vectors are 
$T(s)=\bpm \ga_1'(s), \ga_2'(s)\epm^T$ and $\nu(s)=\bpm -\ga_2'(s), \ga_1'(s)\epm^T$. 
For a perturbation $u$ we consider $\tilde \ga(s) =\ga(s) +\del u(s)\nu(s)=\bpm \ga_1(s) -\del u(s)\ga'_2(s)\\ \ga_2(s)+\del u(s) \ga'_1(s)\epm $, with unit normal vector $\tilde \nu =\bpm -\tilde \ga'_2(s)\\ \tilde \ga'_1(s)\epm / \| \tilde \ga'(s)\|$. Then 
\hualst{ \frac \dd {\dd \del} \bigg|_{\del=0}&\int_0^{2\pi} \spr{\tilde \ga(s), \tilde \nu(s) }\|\tilde \ga'(s)\| \dd s= 
\frac \dd {\dd \del}\bigg|_{\del=0} \int_0^{2\pi} \tilde \ga_2(s)\tilde \ga'_1(s)-\tilde \ga_1(s)\tilde\ga'_2(s) \dd s\\
&= \int_0^{2\pi}  \frac \dd {\dd \del}\bigg|_{\del=0} \Big[\big( \ga_2(s)+\del u(s)\ga_1'(s)\big)\Big(\ga_1'(s) - \del \big(u(s)\ga_2'(s)\big)'\Big) \\
&\hs{15mm}- \big( \ga_1(s)-\del u(s)\ga_2'(s)\big)\Big(\ga_2'(s) + \del \big(u(s)\ga_1'(s)\big)'\Big) \Big]\dd s \\
&=\int_0^{2\pi} u(s) \ga_1'(s)^2 -\ga_2(s) (u(s)\ga_2'(s))' + u(s)\ga_2'(s)^2 - \ga_1(s)(u(s)\ga_1'(s))' \dd s\\
&=2\int_0^{2\pi} u(s)\dd s,
}
which can be easily implemented using the (1D) mass matrix of the 
discretization of $\pa X_{\text{left}}$. 

\paragraph{Numerical flows} We give a brief description of 
our ad hoc method for numerical flows (aka direct numerical simulation DNS) 
for \reff{hflow}, 
and refer to \cite{BMN05,BGN08,BGN20,BGN21} and the references therein for 
more sophisticated and provably stable FEM schemes 
for geometric flows, including the Helfrich flow \reff{hflow} and extensions, and in 
particular also dealing with tangential motion of mesh points. 
Here we only consider normal displacements, and 
approximate the normal velocity  as 
$$\ds V_t=\spr{\pa_t X,N}\approx \spr{\frac 1 h M_{n}(X_{n+1}{-}X_n),N_n}
=\frac 1 h \spr{M_{n}(X_n{+}u_{n+1}N_{n}{-}X_n),N_n}=\frac 1 h 
M_{n}u_{n+1},$$
where $h$ is the current stepsize, $n$ the time stepping index, 
and $M_{n}$ the mass matrix. 
Setting $U_n=(u_n,\Lam_n)$, the PDE $V_t=G(X(t),\Lam(t))$, 
together with the (volume and area) constraints $q(u_{n+1})=0$ (translational 
and rotational phase conditions play no role in the DNS),  
yields the nonlinear system 
\huga{\label{tstep1} 
\CG(U_{n+1}):=\bpm \frac 1 h M_{n} u_{n+1}+G(U_{n+1})\\
q(u_{n+1})\epm\stackrel!=\bpm 0\\0\epm\in \CH\times \R^{n_q} 
} 
for $U_{n+1}$, where $\CH$ stands for the FEM space used. 
Given $h$ and $X_n$, we solve \reff{tstep1} for $U_{n+1}$, 
up to a given tolerance \tol, usually $\tol=10^{-6}$ in the sup--norm, via 
Newton loops. 
For simplicity we do not apply an extrapolation 
step, i.e., the initial guesses for the Newton loops are 
$u_{n+1}^{(0)}=0$ (because $X_n$ has already been updated) and 
$\Lam_{n+1}^{(0)}=\Lam_n$.  Note that there are no 
initial conditions for $\Lam$ (only for $X_0$, which 
corresponds to $u_0=0$), 
and we use $\Lam_1^{(0)}=\Lam_0$, the values of $\Lam$ at the 
steady state which we perturb. 
For the Newton loops we essentially have all 
tools (in particular Jacobians) available from the continuation setup. 
This also suggest some time stepping control (decrease $h$ if many Newton 
iterations are required or no convergence, increase $h$ if fast convergence). 

Additionally we remark that
Newton loops for \reff{tstep1} often fail to converge near straight cylinders 
(or more generally CMC surfaces), where $\pa_u q_1$ and 
$\pa_uq_2$ become linearly dependent ($\CD$ from \reff{ddef} becomes singular),  
but this is not a major concern for us as we are 
more interested in the flow towards nontrivial steady states.  
See, e.g., \cite[Scheme (2.8)]{BGN08} for dealing with 
the constraints when approaching CMC surfaces, which essentially 
consists in setting $\lam_{1,n+1}=\lam_{1,n}$ and only solving 
for $\lam_{2,n+1}$ with $q(u_{n+1})=q_2(u_{n+1})$. 
Moreover, the DNS in \S\ref{numex} are all uncritical wrt to mesh 
quality. 
However, for stronger initial perturbations we expect that allowing 
and controlling tangential motion of mesh points (or other ways of 
mesh--adaptation) will be crucial, see also, e.g., the numerical flows 
for closed vesicles in \cite{geomtut}.

\bibliographystyle{alpha}
\bibliography{./hubib.bib} 

\end{document}